\documentclass[preprint,11pt]{article}
\usepackage{color}
\usepackage[centertags]{amsmath}
\usepackage{amsfonts}
\usepackage{amssymb}
\usepackage{amscd}
\usepackage{amsthm}
\usepackage{color}
\usepackage{graphicx}
\usepackage{fullpage}
\usepackage{subfigure}
\usepackage{epstopdf}
\usepackage{alltt}
\usepackage{float}
\usepackage{url}

\usepackage{lipsum}
\usepackage{amsfonts}
\usepackage{graphicx}
\usepackage{epstopdf}
\usepackage{algorithmic}
\usepackage{a4wide} 
\usepackage{amsmath,amsfonts,amssymb,latexsym}
\usepackage{graphicx}
\usepackage[latin1]{inputenc}
\usepackage{color}
\usepackage{siunitx,tabularx,etoolbox,xcolor,booktabs} 

\usepackage{tikz}
\usepackage[outline]{contour}

\contourlength{1.2pt}

\ifpdf
  \DeclareGraphicsExtensions{.eps,.pdf,.png,.jpg}
\else
  \DeclareGraphicsExtensions{.eps}
\fi

\usepackage{amsopn}

\newcommand{\Hone}{H^{1}(\Omega_{d})}

\newcommand{\Honel}{H^{1}_0(\Omega_{\pml})}
\newcommand{\Honeln}{H^{1}(\Omega_{\pml})}
\newcommand{\Ltwol}{L^{2}(\Omega_{\pml})}

\newcommand{\R}{\mathbb{R}}
\newcommand{\C}{\mathbb{C}}

\newcommand{\HoneR}{H^{1}(\R)}

\providecommand*{\mrm}[1]{\mathrm{#1}}
\newcommand{\iu}{\mrm{i}}

\newcommand{\diff}{\,\mrm{dx}}
\DeclareMathOperator{\Ran}{Ran}
\DeclareMathOperator{\Ker}{Ker}

\def\sa{\mathfrak{a}}
\def\sb{\mathfrak{b}}
\def\sq{\mathfrak{q}}
\def\bu{\tilde{u}}

\def\cL{\mathcal{A}}

\newcommand{\HH}{\mathcal{V}}
\newcommand{\Hh}{\mathcal{V}^{\nu}}

\newcommand{\E}{\mathcal{E}}

\newcommand{\Ph}{{P}^{\nu}}
\def\H{\Hone}

\newcommand{\W}{\mathcal{W}}
\newcommand{\Wh}{\mathcal{W}^{\nu}}

\newcommand{\fem}{\nu} 

\newcommand{\dis}{\displaystyle}

\newcommand{\Om}{\Omega}

\newcommand{\nG}[1]{{\boldsymbol #1}}

\newcommand{\fr}[2]{\displaystyle\frac{#1}{#2}}

\newcommand{\mx}[1]{\displaystyle\mathbb{#1}} 
\newcommand{\mc}[1]{\mathcal{#1}} 

\newcommand{\dtot}[2]{\displaystyle\frac{d #1}{d #2}}

\newcommand{\re}[1]{\Re #1}
\newcommand{\im}[1]{\Im #1}

\newcommand{\dtn}{d} 
\newcommand{\pml}{\ell} 
\newcommand{\kett}{a} 
\newcommand{\res}{a} 
 
\numberwithin{equation}{section}
\numberwithin{figure}{section}

\definecolor{cadmiumgreen}{rgb}{0.0, 0.42, 0.24}
\definecolor{cardinal}{rgb}{0.77, 0.12, 0.23}

\def\cL{\mathcal{L}}

\providecommand{\keywords}[1]{\textbf{\textit{Keywords:}} #1}

\begin{document}

\title{On spurious solutions in finite element approximations of resonances in open systems}

\author{Juan Carlos Araujo-Cabarcas,\and Christian Engstr\"om \\ \and
Department of Mathematics and Mathematical Statistics, Ume\aa\, University, Ume\aa, Sweden}

\maketitle

\begin{abstract}
In this paper, we discuss problems arising when computing resonances with a finite element method. In the pre-asymptotic regime, we detect for the one dimensional case, spurious solutions in finite element computations of resonances when the computational domain is truncated with a perfectly matched layer (PML) as well as with a Dirichlet-to-Neumann map (DtN). The new test is based on the Lippmann-Schwinger equation and we use computations of the pseudospectrum to show that this is a suitable choice. 
Numerical simulations indicate that the presented test can distinguish between spurious eigenvalues and true eigenvalues also in difficult cases.  
\end{abstract}


\section{Introduction}
Open resonators are common in many applications including acoustic properties of musical instruments, laser
cavities, and multilayer x-ray resonators \cite{ELEJABARRIETA2002584,Karman1999,schenk2011optimization}. 
A widely applicable technique to terminate the computational domain in resonance 
as well as in scattering problems is a perfectly matched layer (PML). It is possible to prove that the Galerkin method converges (in gap), which implies that in the asymptotic regime there are no Galerkin eigenvalues that are unrelated to the spectrum of the original operator. However, under realistic conditions it is very costly to use sufficiently fine meshes and most computations are therefore done in the pre-asymptotic regime. The PML problem is highly non-normal and there are frequently numerous eigenvalues that are unrelated to the spectrum of the original operator. These eigenvalues are called spurious eigenvalues and they are a major challenge in engineering applications.

A common approach that aims to detect spurious solutions is to compute numerical approximations with several sets of PML parameters. Then a perturbation argument is used to distinguish true resonances from spurious solutions \cite{gopal08,KettnerSchmidt2011,6632358}.
This require several computations of the eigenvalues for different parameters. The basic assumption in this approach is that spurious eigenvalues react stronger to perturbations than approximations to resonances. However, eigenvalues may also react strongly due to an insufficient approximation and it is unclear how much an eigenvalue should move to be marked as a spurious eigenvalue.


In this paper, we propose a new test based on a volume integral formulation of the problem called the Lippmann-Schwinger equation and argue that this is a suitable choice.  This test is applied to numerical solutions obtained with a finite element method where the computational domain is truncated with a PML as well as with a 
DtN map. Numerical simulations indicate that the presented test can distinguish between spurious eigenpairs and true eigenpairs also in complicated cases when the spurious solutions mix with true approximations of resonances. 
The new test determines if a computed eigenvalue numerically is in an $\epsilon$-pseudospectrum of the integral operator and the corresponding vector is an $\epsilon$-pseudomode. The Lippmann-Schwinger equation is frequently used in acoustic and electromagnetic scattering theory \cite{Colton+Kress1983} and it has previously been used directly to determine resonances \cite{Kao2008412,gopal08,osting13}. The resonant modes grow exponentially at infinity but the integration in the Lippmann-Schwinger equation is only over the resonator where the solutions are well behaved. 
This is a computationally significant advantage in particular for large structures that contain many air holes. However, the direct approach with an integral equation is demanding since it results in a nonlinear eigenvalue problem (NEP), matrices are full, and each evaluation (e.g solver iteration) requires a matrix assembly.

An advantage with the DtN map in one dimension is that the resulting eigenvalue problem only has a quadratic nonlinearity and the formulation contains no free parameters. The PML has the advantage that the resulting eigenvalue problem is linear. However, the method contains several parameters, which influence the computational result.

For the DtN formulation, we use a block operator representation and apply standard techniques to prove an estimate for the gap between the discrete and continuous eigenspaces. Moreover, based on the results \cite{MR0383117,kim09} we state the corresponding estimate for the gap in the PML setting. Convergence in gap proves for the DtN formulation that no spurious solutions exist when the finite element space is large enough but the rate of convergence depend critically on resolvent norms that may be very large. The same conclusions hold in the PML formulation with the additional requirement that the PML layer is thick enough \cite{kim09}. The eigenvalue problem with a truncated PML will have more eigenvalues inside a given region in the complex plane compared with a DtN formulation of the same problem. We introduce a DtN for the truncated PML and determine a region in the complex plane where it is possible to obtain convergence for given PML parameters. Then, we derive a new estimate of the difference between a resonance and an eigenvalue of the finite PML problem. 
Finally, the integral equation based test is then used to detect spurious solutions in the DtN and PML formulations. The numerical examples indicate that the test can detect spurious eigenvalues in solutions computed with relatively coarse discretizations (computed with $h$-FEM) as well as for fine discretizations obtained using $p$-FEM.

The results in this paper are stated for the one-dimensional case, which for the considered test problems, it is possible to reach the asymptotic regime and compute resonances without spurious solutions on a standard computer. This is an advantage, since in this case it is possible to evaluate numerically the performance of the new filtering technique. The filtering process can be extended to higher dimensions and it will then provide a new practical tool for many challenging applications in physics and engineering. 

A procedure to efficiently compute resonances with the PML or DtN formulation is then: (i) Use a course discretization and e. g. ARPACK with several shifts to compute a selection of Galerkin eigenvalues. (ii) Use the new filtering process to sort out approximations of interest and reduce the number of shifts. (iii) Use hp-adaptivity, e.g. a technique similar to \cite{MR2326196} for PML and \cite{MR3529053} for DtN, to reduce the errors in the target eigenvalues. (iv) Check that the new eigenpars significantly reduce the residual in the Lippmann-Schwinger equation.

\section{Convergence of Galerkin spectral approximations for the DtN and PML formulations}\label{sec:convergence}

In this section we introduce the DtN and PML formulations used to truncate the exterior domains and prove convergence of the Galerkin method. Our approach to analyze the DtN formulation follow \cite{MR482047,MR3164142} and we state convergence results that was proved by Bramble and Osborn, et al. \cite{MR1115240}. The results stated for the infinite PML formulation is contained in \cite{kim09}. For the finite PML problem we introduce a DtN map that is used to derive a new error estimate and reference solutions.

\subsection{Computing Resonances with the DtN map}\label{sec:DtN_1d}
Resonance problems are closely related to the underlying scattering problem and we begin therefore with the Helmholtz scattering problem on $\R$ \cite{Colton+Kress1992}. Consider the scattering of a given incoming wave $u_{i}$ by an obstacle $n$ with support in  $\text{supp}\, (n^2-n_0^2)\subset (-\dtn,\dtn)$. Then the outgoing radiation condition on the scattered wave $u_{s}$ is
\begin{equation}\label{eq:outgoing}
	u_s'(-x_0)=-i k \,n_0 \,u_s(-x_0), \quad u_s'(x_0)  =i k\,n_0 \, u_s(x_0),\quad x_0\geq d
\end{equation}
and a function that satisfies \eqref{eq:outgoing} is called \emph{outgoing} \cite{Colton+Kress1992}. 
The scattering problem is then: Find for given $k^2$ with $\Im k^2\geq 0$ the total wave $H^2 (\R)\ni u=u_{i}+u_{s}$ with $u_s$ outgoing, that satisfies 
\begin{equation}\label{eq:HelmBase}
 	-u''-k^2 n^2u=0.
 \end{equation}
The condition \eqref{eq:outgoing} on $u_{s}$ ensures uniqueness of the solution \cite[p. 348]{MR996423}, \cite{Colton+Kress1992}.

For $x\not\in\Omega_\dtn$ and given non-zero $k\in\C$, equation \eqref{eq:HelmBase} has the linearly independent solutions $e^{\pm ik n_0 x}\in H^2_{\hbox{\tiny loc}}(\mx R)$. For $x\geq d$, we have that $u_s(x)=e^{ik n_0 x}$ is the outgoing solution and $u_s(x)=e^{-ik n_0 x}$ is called the \emph{incoming} solution. Similarly,  for $x\leq -d$  the function $u_s(x)=e^{-ik n_0 x}$ is outgoing and $u_s(x)=e^{ik n_0 x}$ is the incoming solution. 

The standard definition of resonances as the poles of the analytical continuation of the resolvent operator is discussed in Sect. \ref{sec:resolvent}. However, these resonances can also be determined by solving a nonlinear eigenvalue problem, where the nonlinearity comes from the  Dirichlet-to-Neumann (DtN) map \cite{lenoir92,schenk2011optimization}. In one space dimension, the problem formulation is (formerly) given by \eqref{eq:outgoing}, \eqref{eq:HelmBase} with $u_{i}=0$. Hence, the resonance problem restricted to $\Omega_\dtn:=(-\dtn,\dtn)$ is: Find a non-zero $u\in H^2(\Omega_\dtn)$ and $k\in\C$ such that
\begin{equation}\label{eq:Helmholtz}
	-u''-k^2n^2 u=0 \,\,\,\hbox{for}\,\, x\in\Omega_\dtn,
\end{equation}
where the DtN map at $x=\pm d$ is
\begin{equation}\label{eq:formalDtN}
	u'(-d)=-i k \,n_0 \,u(-d), \quad u'(d)  =i k\,n_0 \, u(d).
\end{equation}
Note that the eigenvalues $k$ of the resonance problem  \eqref{eq:Helmholtz},  \eqref{eq:formalDtN} will have negative imaginary part, which for the resonance problem posed on $\R$ implies that an outgoing solution grows exponentially at infinity.

Below, we write \eqref{eq:Helmholtz}-\eqref{eq:formalDtN} on a variational form and state results for a conforming Galerkin finite element discretization of a reformulation of the problem as a linear pencil. Let $\sb$ denote a bounded sesquilinear form on $\Hone\times\Hone$. Then $\sb$ is called compact on $\Hone$ if
\begin{equation}\label{eq:compform}
	\sup_{\|v\|_{\Hone}\leq 1}|\sb [u_n-u,v]|\to 0,
\end{equation}
for every weakly convergent sequence $u_n\to u$ \cite{DEMKOWICZ199469}. In the analysis, we use that a compact sesqulinear form \eqref{eq:compform} corresponds to a compact operator on $\Hone$ \cite{DEMKOWICZ199469}. 
Define the continuous sesquilinear form $\sa_{0}: \Hone\times\Hone\rightarrow\C$,
\begin{equation}\label{a0}
		\sa_{0}[u,v] :=\int_{-\dtn}^{\dtn}u'\overline{v}'\diff.
\end{equation}
The trace operators $\tau_{\pm}:\Hone\rightarrow \C$, $\tau_{\pm}u=u(\pm \dtn)$ have finite rank and the sesqulinear form $n_0\left(z_1\overline{v}(\dtn)+z_2\overline{v}(-\dtn)\right)$ is bounded for all $z=(z_1,z_2)\in\C^2$ and $v\in\Hone$. Hence, from the compactness of $\tau_{\pm}$ follows that
\begin{equation}\label{a1}
	\sa_{1} : \Hone\times\Hone\rightarrow\C,\quad \sa_{1}[u,v] :=
	n_0\left( u(\dtn)\overline{v}(\dtn)+u(-\dtn)\overline{v}(-\dtn)\right),
\end{equation}
is compact. For $r>s$ the embedding $H^{s}(\Omega_\dtn)\subset  H^{r}(\Omega_\dtn)$ is compact \cite[Theorem 7.2]{MR895589}. Hence, the bounded form
\begin{equation}\label{a2}
	    	\sa_{2} : L^2(\Omega_d)\times\Hone\rightarrow\mathbb{C},\quad\sa_{2}[u,v]:=\int_{-\dtn}^{\dtn}n^2u\overline{v}\diff
\end{equation}
is compact on $\Hone\times\Hone$. Set $\lambda=-\iu k$ and define for $u,v \in\Hone$ and $\lambda\in\C$ the form-valued function
\begin{equation}\label{eq:quadratic_form}
	\sq(\lambda)[u,v]:=\lambda^{2} \sa_{2}[u,v]+\lambda \sa_{1}[u,v]+\sa_{0}[u,v].
\end{equation}
The form $\sa_{0}$ has a non-empty kernel and we define therefore for fixed $\alpha>0$ the shifted form $\hat\sq(\lambda)=\sq(\lambda+\alpha)$,
\begin{equation}\label{weak_dtn}
	\hat\sq(\lambda)[u,v]:=\lambda^{2} \hat\sa_{2}[u,v]+\lambda \hat\sa_{1}[u,v]+\hat\sa_{0}[u,v],
\end{equation}
where
\begin{equation}\label{weakTM2}
	 \begin{aligned}
		\hat\sa_{0}[u,v] & :=\sa_{0}[u,v]+\alpha\sa_{1}[u,v]+\alpha^2\sa_{2}[u,v],\\
		\hat \sa_{1}[u,v] & :=\sa_{1}[u,v] +2\alpha\sa_{2}[u,v],\quad  \hat \sa_{2}[u,v] :=\sa_{2}[u,v] .
	\end{aligned}
\end{equation}
The quadratic eigenvalue problem is then as follows: Find vectors $u\in\Hone\backslash\{0\}$ and complex numbers $\lambda$ satisfying 
\begin{equation}
	\hat\sq(\lambda)[u,v]=0
	\label{TMii}
\end{equation}
for all $v\in \Hone$.
The shifted sesquilinear form $\hat\sa_{0}$ is coercive: $\hat\sa_{0}[u,u]\geq C\|u\|_{\Hone}$, $C>0$ and it will be used as the inner product on $\Hone$.
Since the forms  $\hat \sa_{n}$, $n=1,2$ are compact on $\Hone\times\Hone$ the operators $Q_{n}:\Hone\rightarrow\Hone$, $n=1,2$ defined by
\begin{equation}
	\hat\sa_{0}[Q_n u,v]:=\hat \sa_{n}[u,v] \quad \text{for all}\,\, u,v\in\Hone,
\end{equation}	
are compact \cite{DEMKOWICZ199469}. Define the operator polynomial
\begin{equation}\label{eq:Q}
		Q(\lambda):=I+\lambda Q_{1} +\lambda^{2}Q_{2}
\end{equation}
in $\Hone$. The operator polynomial \eqref{eq:Q} is a compact perturbation of the identity with an analytic dependence of $\lambda$. Hence, it follows from the analytic Fredholm theorem \cite[Theorem 1.3.1]{Mennicken03} that all eigenvalues, given by $Q(\lambda)u=0$, are isolated and of finite multiplicity. Numerical analysis of this eigenvalue problem can be based on the general theory for analytic Fredholm operator functions; see \cite{MR1393166} and the references therein. However, an alternative approach is to study a corresponding block operator matrix formulation \cite{Markus1988}. Then, approximation theory of linear non-selfadjoint operators can be applied \cite{MR482047,Descloux1978a}. In particular, the results in this section show that the a-posteriori error estimations in \cite{MR3529053} can be applied to the DtN-formulation of the resonance problem.

Let $\W:=\Hone\oplus\Hone$  and assume that $(u_1\, u_2)^t\in\W\setminus \{0\}$, $\lambda\in \C$ is a solution of the generalized eigenvalue problem 
\begin{equation}\label{eq:Lgen}
	\begin{bmatrix}
			I & Q_1 \\
			0 & I
		\end{bmatrix}
	\begin{bmatrix}
			u_1 \\
			u_2
		\end{bmatrix}	
 =\lambda
	\begin{bmatrix}
		\phantom{-}0 & -Q_2 \\
		\phantom{-}I & \phantom{-}0
	\end{bmatrix}
	\begin{bmatrix}
			u_1 \\
			u_2
		\end{bmatrix}.
\end{equation}
Then $Q(\lambda)u_1=0$ follows, and it can be shown that $Q$ and \eqref{eq:Lgen} have the same eigenvalues and that they have the same multiplicities \cite[Lemma 12.5]{Markus1988}.
The operator on the left hand side of \eqref{eq:Lgen} is invertible and we define the operator $\cL:\W\to\W$,
\begin{equation}\label{op:L}
	\begin{bmatrix}
			I & Q_1 \\
			0 & I
		\end{bmatrix}^{-1}	
	\begin{bmatrix}
		\phantom{-}0 & -Q_2 \\
		\phantom{-}I & \phantom{-}0
	\end{bmatrix}
=\begin{bmatrix}
			-Q_1 & -Q_2 \\
			\phantom{-}I & 0
		\end{bmatrix}=:\cL. 
\end{equation}
Consequently, the spectrum of the quadratic operator function $Q$ coincide with the set of numbers
\begin{equation}
	\sigma (Q)=\{\lambda\in\C\,:\,\lambda=1/\mu,\,\mu\in \sigma (\cL)\}.
\end{equation}

Let $\Hh\subset\H$ denote a sequence of conforming finite element spaces with the approximation property
\begin{equation}
\label{eq:P1}
	\lim_{\text{dim} (\Hh)\rightarrow \infty}\inf_{u^{\nu}\in\Hh}||u-u^{\nu}||_{\H}=0, \quad \text{for all}\, u\in\H.
\end{equation}
Let $\|u\|_{\hat\sa_{0}}:=\sqrt{\hat\sa_{0}[u,u]}$ and let $\Ph:\H\rightarrow \Hh$ denote the projection of $\H$ into $\Hh$ defined by the inner product 
$\hat\sa_{0}[\Ph u,v^{\nu}] =\hat\sa_{0}[u,v^{\nu}]$ for all $v^{\nu}\in\Hh$. For the operator polynomial \eqref{eq:Q} define the projected operator function $Q^{\nu}:\HH\rightarrow \Hh$ by $Q^{\nu}=\Ph Q$. The Galerkin eigenvalue problem is to find vectors $u^{\nu}\in\Hh\setminus\{0\}$ and values $\lambda^{\nu}\in\C$ such that
\begin{equation}
\label{Galerkin}
	Q^{\nu}(\lambda^{\nu})u^{\nu}=0.
\end{equation}
Let $\Wh:=\Hh\oplus\Hh$. The corresponding Galerkin eigenvalue problem for \eqref{eq:Lgen} is to find vectors $\bu^{\nu}\in\Wh\setminus\{0\}$ and values $\lambda^{\nu}\in\C$ such that
\begin{equation}\label{eq:PLgen}
	\begin{bmatrix}
			\Ph & \Ph Q_1 \\
			0 & \Ph
		\end{bmatrix}
	\begin{bmatrix}
			u_1^{\nu} \\
			u_2^{\nu}
		\end{bmatrix}	
 =\lambda^{\nu}
	\begin{bmatrix}
		\phantom{-}0 & -\Ph Q_2 \\
		\phantom{-}\Ph & \phantom{-}0
	\end{bmatrix}
	\begin{bmatrix}
			u_1^{\nu} \\
			u_2^{\nu}
		\end{bmatrix}.
\end{equation}
We define as in \cite{MR482047,MR3164142} an auxiliary problem, which has the same eigenvalues and generalized eigenvectors as \eqref{eq:PLgen} on $\Hh\oplus\Hh$ but it is possible to show convergence in norm. Define for $\bu^{\nu}\in\W^h\setminus\{0\}$, $\lambda^{\nu}\in\C$ the generalized eigenvalue problem 
\begin{equation}\label{eq:PLgen2}
	\begin{bmatrix}
			I & \Ph Q_1 \\
			0 & I
		\end{bmatrix}
	\begin{bmatrix}
			u_1^{\nu} \\
			u_2^{\nu}
		\end{bmatrix}	
 =\lambda^{\nu}
	\begin{bmatrix}
		\phantom{-}0 & -\Ph Q_2 \\
		\phantom{-}I & \phantom{-}0
	\end{bmatrix}
	\begin{bmatrix}
			u_1^{\nu} \\
			u_2^{\nu}
		\end{bmatrix}
\end{equation}
and the block operator matrix formulation
\begin{equation}\label{GeqL}
	\cL^{\nu} \bu^{\nu} =\mu^{\nu} \bu^{\nu},\quad \cL^{\nu} =
		\begin{bmatrix}
			-\Ph Q_1 & -\Ph Q_2 \\
			\phantom{-}I & 0
		\end{bmatrix}.
\end{equation}
In  \cite{MR482047} it has been shown that \eqref{eq:PLgen} and \eqref{eq:PLgen2} have the same spectrum and $\cL^{\nu}\rightarrow \cL$ in norm. 

Given a circle $\gamma_{\mu}\in\rho(\cL)$ which encloses $\mu\in\sigma(\cL)$ and no other elements of $\sigma(\cL)$ the Riesz projections $E(\mu)$ and $E^{\nu}(\mu)$ are defined by
\begin{equation}\label{NonOrtPrj}
\begin{aligned}
	E(\mu;\cL) &=\frac{1}{2\pi\iu}\int_{\gamma_{\mu}}(z-\cL)^{-1}dz,\\
	E^{\nu}(\mu;\cL^{\nu}) &=\frac{1}{2\pi\iu}\int_{\gamma_{\mu}}(z-\cL^{\nu})^{-1}dz.
	\end{aligned}
\end{equation}
The range $\E_{\mu}$ of the operator $E(\mu;\cL):\W\rightarrow \W$,
\begin{equation}
	\E_{\mu}:=\Ran(E (\mu;\cL))=\Ker(\mu-\cL)^{\alpha},
\end{equation}
is the corresponding generalized eigenspace and $\alpha$ is the smallest positive integer such that $\Ker(\mu-\cL)^{\alpha}=\Ker(\mu-\cL)^{\alpha+1}$. Moreover, let $\E_{\mu}^{\nu}$ denote the range of the projection $E^{\nu}(\mu;\cL^{\nu})$. Define for closed subspaces $V_1$ and $V_2$ of a Hilbert space the gap $\hat\delta$ between $V_1$ and $V_2$ as
\begin{equation}
	\delta (V_1,V_2):=\sup_{v_1\in V_1,\,\|v_1\|_{V_1}=1}\text{dist}\,(v_1,V_2),\quad \hat{\delta} (V_1,V_2)=\max\left (\delta (V_1,V_2),\delta (V_2,V_1) \right ).
\end{equation}
Let $(\cL-\cL^{\nu})|_{\E_{\mu}}$ denote the restriction of $\cL-\cL^{\nu}$ to $\E_{\mu}$ and assume that $\text{dim}\,\E_{\mu}=\text{dim}\,\E_{\mu}^{\nu}$. Then follows
\begin{equation}\label{eq:OsbornEst}
	\hat{\delta} (\E_{\mu},\E_{\mu}^{\nu})\leq C\,\text{length}(\gamma_{\mu})\sup_{z\in\gamma_{\mu}}\|(z-\cL)^{-1}\|\sup_{z\in\gamma_{\mu}}\|(z-\cL^{\nu})^{-1}\|\|(\cL-\cL^{\nu})|_{\E_{\mu}}\|
\end{equation}
for some $C>0$, \cite[Theorem 1]{MR0383117}. We know that $\hat{\delta} (\E_{\mu},\E_{\mu}^{\nu})\rightarrow 0$ since $\cL^{\nu}\rightarrow \cL$ in norm, but the condition $\hat{\delta} (\E_{\mu},\E_{\mu}^{\nu})< 1$ should be satisfied for a finite dimensional $\W^{\nu}$ to guarantee that $\text{dim}\,\E_{\mu}=\text{dim}\,\E_{\mu}^{\nu}$ \cite[IV.2 Corollary 2.6]{Kato1980}. Note that $\delta (\E_{\mu},\E_{\mu}^{\nu})< 1$ only implies $\text{dim}\,\E_{\mu}\leq\text{dim}\,\E_{\mu}^{\nu}$ and in particular $\delta (\E_{\mu},\E_{\mu}^{\nu})=0$ implies that $\E_{\mu}\subset\E_{\mu}^{\nu}$. The estimate  \eqref{eq:OsbornEst} depends on the norm of the resolvents $(z-\cL)^{-1}$ and $(z-\cL^{\nu})^{-1}$ over $\gamma_{\mu}$, which can be very large for non-normal operators \cite{MR2359869}. The inequality
\begin{equation}\label{est:approx}
	\|(\cL-\cL^{\nu})|_{\E_{\mu}}\|\leq C\rho^{\nu}, \quad \rho^{\nu}:=\sup_{u\in \E_{\mu}}\inf_{u^{\nu}\in\W^{\nu}}\|u-u^{\nu}\|_{\W}
\end{equation}
holds and $\rho^{\nu}$ depends on the approximation properties of the finite element space \cite{Babuska+Buo+Osborn1989}. This indicates that very good approximation properties of $\W^{\nu}$ may be necessary to clear a given region in the complex plane from spurious eigenvalues. 
Define
\[
	\hat\E_{\mu}:=\Ran(E (\mu;\cL^*))=\Ker(\mu-\cL^*)^{\alpha}, \quad\hat\rho^{\nu}:=\sup_{\hat u\in \hat\E_{\mu}}\inf_{u^{\nu}\in\W^{\nu}}\|\hat u-\hat u^{\nu}\|_{\W}.
\]
For large enough $\dim(\W^{\nu})$ Kolata \cite{MR482047} proved the following estimates
\begin{equation}
	|\mu-\mu^{\nu}|\leq C\rho^{\nu}\hat\rho^{\nu},\quad \|u-u^{\nu}\|_{\W}\leq C\rho^{\nu}.
\end{equation}
In our case, the eigenfunctions of $\cL$ and of $\cL^*$ will have the same regularity. Hence, we expect that the convergence of the eigenvalues are $\mathcal{O}((\delta^{\nu})^2)$. 

\subsection{Computing Resonances with the PML}\label{sec:PML_1d}

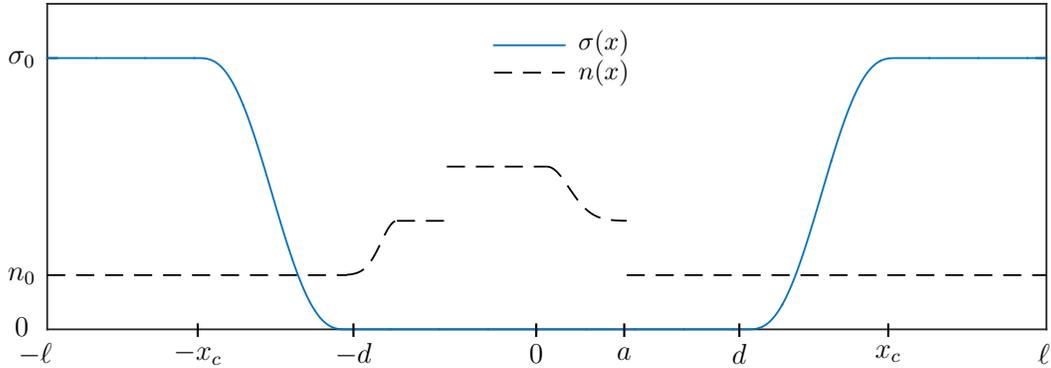
\begin{figure}[!h]
	\centering
	\begin{tikzpicture}[thick,scale=0.9, every node/.style={scale=0.9}]
		\tikzstyle{lab}  = [fill=none,font=\large,inner sep=4pt]
		\tikzstyle{form} = [fill=none,font=\normalsize,inner sep=4pt]
		
		\draw( 0.00, 0.0) node { \includegraphics[scale=1.0]{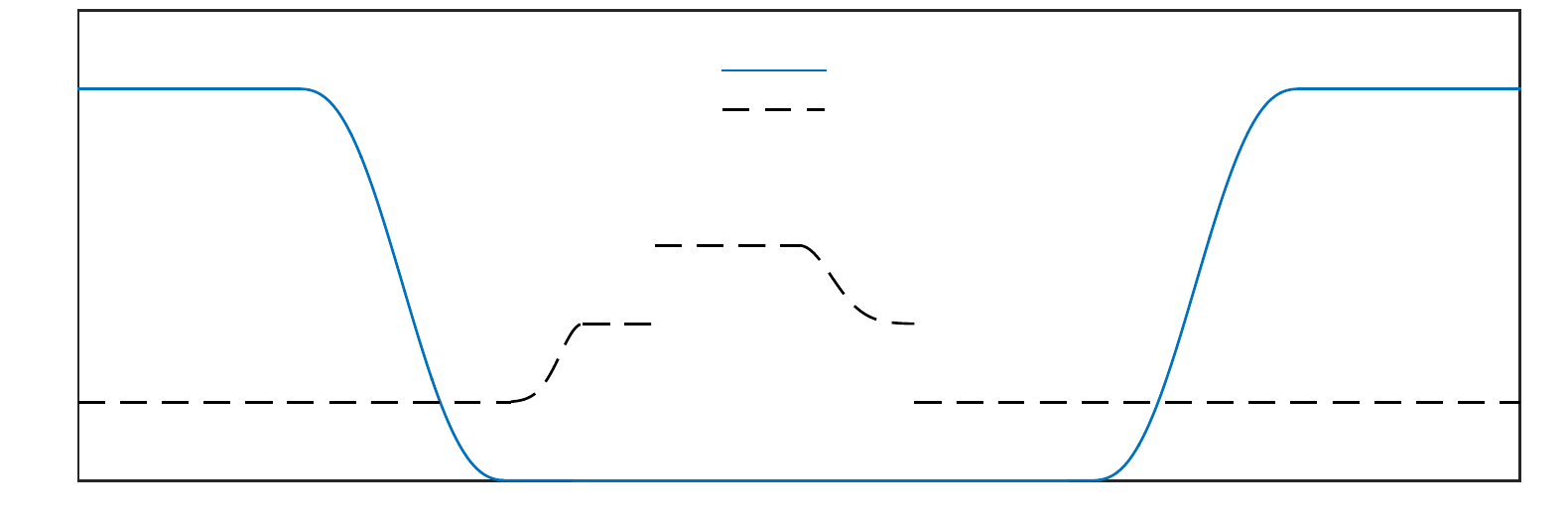} };
		
		\node[lab] at ( 7.5,-2.6) {$ \pml$};
		\node[lab] at (-7.4,-2.6) {$-\pml$};
		\node[lab] at (-7.6,-2.2) {$0$};
		\node[lab] at (-7.6,-1.45) {$n_0$};
		
		\node[lab] at ( 0.0,-2.6) {$0$};
		\draw[thick,-] (0.0,-2.15) -- (0.0,-2.35);   
		
		\node[lab] at ( 1.3,-2.6) {$a$};
		\draw[thick,-] (1.3,-2.15) -- (1.3,-2.35);   
		
		\node[lab] at (-2.7,-2.6) {$-d$};
		\draw[thick,-] (-2.7,-2.15) -- (-2.7,-2.35);
		
		\node[lab] at ( 3.0,-2.6) {$ d$};
		\draw[thick,-] (3.0,-2.15) -- (3.0,-2.35);
		
		\node[lab] at (-5.0,-2.6) {$-x_c$};
		\draw[thick,-] (-5.0,-2.15) -- (-5.0,-2.35);
		
		\node[lab] at ( 5.2,-2.6) {$ x_c$};
		\draw[thick,-] (5.2,-2.15) -- (5.2,-2.35);   
		
		\node[lab] at (-7.6,1.75) {$\sigma_0$};
		
		\node[form] at ( 1.0,2.0) {$\sigma(x)$};
		\node[form] at ( 1.0,1.55) {$n(x)$};
		
	\end{tikzpicture}
	\vspace*{-3mm}
	\caption{\emph{Finite PML strength function in solid line, and a typical refractive index profile in dashed line.}}
	\label{fig:pml_strength}
\end{figure}
In the previous section, a DtN-map was used to reduce the exterior Helmholtz problem to a bounded domain. 
Then, resonances were computed by solving the quadratic eigenvalue problem \eqref{TMii}. In this section, we consider an alternative approach based on a complex coordinate stretching (the PML method), which results in a linear eigenvalue problem \cite{kim09}. 
Take $x_c>d$ and let $P$ be the third order polynomial satisfying: $P(d)=0$, $P'(d)=0$, $P(x_c)=\sigma_0$, and $P'(x_c)=0$. Then, we define the \emph{PML strength function} $\sigma\in C^1(\mx R)$ as
\begin{equation}
	\sigma(x):=\left\{
	\begin{array}{lcl}
		0, & \hbox{for} & |x|\leq d \\
		P(|x|), & \hbox{for} & d<|x|\leq x_c \\
		\sigma_0, & \hbox{for} & x_c<|x|
	\end{array} \right. .
	\label{eq:sigma_fun}
\end{equation}
The chosen PML strength function \eqref{eq:sigma_fun} is increasing for $d<|x|\leq x_c$ and satisfies $\sigma(x)\geq 0$.  The PML problem is in Sect. \ref{sub:finitePML} restricted to $(-\pml,\pml)$ and the PML strength function has then  the profile shown in Fig. \ref{fig:pml_strength}. In the following sections, we consider the complex change of variable and transformation rule
\begin{equation}
	\tilde x = \dis\int_{x_0}^x \alpha(y)\, dy,\,\, \dtot{}{\tilde x}=\frac{1}{\alpha( x)}\dtot{}{ x},
	\,\,\hbox{with}\,\,\alpha(x)=1+i\sigma(x)
	\label{eq:coordstr}
\end{equation} 
where $x_0=-\infty$ in $\R^{-}$ and  $x_0=d$ in $\R^{+}$.

\subsubsection{The infinite PML problem}

Formally, applying \eqref{eq:coordstr} to $-u''-k^2n^2 u=0$ results in the \emph{infinite PML problem}: Find a non-zero $u$ and $k\in\C$ such that
\begin{equation}
	\begin{array}{rc}
		-\dtot{}{x} \left(\frac{1}{\alpha}\dtot{u}{x}\right)-k^2n^2\alpha \,u=0.
	\end{array}
		\label{eq:inf_PML}
\end{equation}
A variational formulation of \eqref{eq:inf_PML}, is written by defining the bounded sesquilinear forms $\sa:\HoneR\times\HoneR\rightarrow\C$, $\sb:\HoneR\times\HoneR\rightarrow\C$,
\begin{equation}\label{eq:defPML}
		\sa[u,v] :=\int_{\R} \frac{1}{\alpha} u' \bar v'  \,dx,\quad \sb[u,v] :=\int_{\R} n^2\alpha\,u \bar v  \,dx,
\end{equation}
and the shifted sesquilinear form $\tilde\sa[u,v] :=\sa[u,v]+\sb[u,v]$. S. Kim and J. E. Pasciak \cite{kim09} proved that the operator $A:\HoneR\rightarrow \HoneR$ defined by
\begin{equation}\label{eq:inf_pml_operator}
	\tilde\sa[A u,v]=\sb[u,v],
\end{equation}	
is well defined and bounded. Note that $Au=\lambda u$, $\lambda=1/(k^2+1)$ implies that $\sa[u,v]=k^2\sb[u,v]$ for all $v\in\HoneR$.

\subsubsection{The finite PML problem}\label{sub:finitePML}

For finite element computations we restrict the domain to $\Omega_{\pml}:=(-\pml,\pml)$ and choose similarly as in  \cite{kim09} homogeneous Dirichlet boundary conditions. Formally, the \emph{finite PML problem} is then: Find the eigenpairs $(u,k)$ such that
\begin{equation}\label{eq:trun_PML}
		-\dtot{}{x} \left(\frac{1}{\alpha}\dtot{u}{x}\right)-k^2n^2\alpha \,u=0,\,\, x\in \Omega_{\pml},\,\,\,
		u(\pml)=0 \,\,\,\hbox{and}\,\,\, 	u(-\pml)=0. 
\end{equation}
In the following, we consider a variational formulation of \eqref{eq:trun_PML} when $0<n_{\min}\leq n(x)\leq n_{\max}$ for all $x\in\Omega_{\pml}$. Define the bounded sesquilinear forms $\sa_\pml:\Honel\times\Honel\rightarrow\C$, $\sb_\pml:\Honel\times\Honel\rightarrow\C$,
\begin{equation}\label{eq:aPML1}
		\sa_{\pml}[u,v] :=\int_{-\pml}^\pml \frac{1}{\alpha} u' \bar v'  \,dx,\quad \sb_{\pml}[u,v] :=\int_{-\pml}^\pml n^2\alpha\,u \bar v  \,dx,
\end{equation}
and the shifted sesquilinear form
\begin{equation}\label{eq:aPML}
		\tilde\sa_{\pml}[u,v] :=\sa_{\pml}[u,v]+\sb_{\pml}[u,v].
\end{equation}
The eigenvalue problem is then as follows: Find vectors $u\in\Honel\backslash\{0\}$ and complex numbers $\lambda$ satisfying 
\begin{equation}\label{eq:pml_fem}
	\tilde\sa_{\pml}[u,v] =\lambda \sb_{\pml}[u,v],\quad \lambda=k^2+1,\,\,\hbox{for all }\,\,v\in \Honel.
\end{equation}
A straightforward calculation shows that $\tilde\sa_{\pml}$ is coercive: $\re{\tilde\sa_{\pml}[u,u]}\geq C\|u\|^{2}_{\Honeln}$,\\
$C= \min\{n_{\min}^2,1/(1+\sigma_0^2)\}$. Hence,  $\tilde\sa_{\pml}$ satisfies the inf-sup condition and $\sb_{\pml}$ is compact since it is continuous on $\Ltwol\times\Honel$ and the embedding $\Ltwol\subset \Honeln$ is compact. Then, it follows that the operator $A_{\pml}:\Honel\rightarrow \Honel$ defined by
\begin{equation}\label{eq:finite_pml_operator}
	\tilde\sa_{\pml}[A_\pml u,v]=\sb_{\pml}[u,v],\quad\text{for all}\ v\in\Honel
\end{equation}	
is compact \cite{MR482047}. Moreover, it can be shown \cite[Theorem 3.1]{kim09} that $A_{\pml}$ can be extended to a bounded operator with domain $\HoneR$.

Let $V:=\Ran(E (\lambda; A_\pml))$ denote the generalized eigenspace associated with $\lambda$, where $E (\lambda; A_\pml)$ is the Riesz projection
\begin{equation}\label{PML:Riesz}
	E(\lambda;A_\pml)=\frac{1}{2\pi\iu}\int_{\gamma_{\lambda}}(z-A_\pml)^{-1}dz.
\end{equation}
Let $S^{\fem}_0\subset\Honel$ denote a sequence of conforming finite element spaces and define the projection $P^{\nu}: \Honel\rightarrow S^{\fem}_0$ by
\[
	\tilde\sa_\pml[P^{\nu} u,v]=\tilde\sa_\pml[u,v], u\in\Honel,\quad v\in S^{\fem}_0.
\]
The operator on $S_0^{\fem}$ can then be written $A^{\nu}_{\pml}=P^{\nu}A_{\pml}$ and $A^{\nu}_{\pml}\rightarrow A_{\pml}$ in norm \cite{MR0383117,kim09}. Let $V_{\pml}^{\nu}:=\Ran(E (\lambda; A^{\nu}_\pml))$ denote the generalized eigenspace associated with $\lambda$. In the following, we consider the generalized eigenspaces $V_{\pml}$, $V_{\pml}^{\nu}$, as subspaces of $\HoneR\supset\Honel$. Then
\begin{equation}
	\hat{\delta} (V_{\pml},V_{\pml}^{\nu})\leq C\,\text{length}(\gamma_{\lambda})\sup_{\lambda\in\gamma_{\lambda}}\|(\lambda-A_{\pml})^{-1}\|\sup_{\lambda\in\gamma_{\lambda}}\|(\lambda-A_{\pml}^{\nu})^{-1}\|\|(A_{\pml}-A_{\pml}^{\nu})|_{V_{\pml}}\|,
\end{equation}
\[
	\|(A_{\pml}-A_{\pml}^{\nu})|_{V_{\pml}}\|\leq C\hat\rho^{\fem}_0,\, \quad\hat\rho^{\fem}_0=\sup_{u\in V_{\pml}, \|u\|=1}\inf_{v\in S^{\fem}_0}\|u-v\|_{\Honel},
\]
where $C$ depends on the inverse of the inf-sup constant for the discrete problem \cite[Theorem 1]{MR0383117}, \cite[p. 696]{MR1115240}. Assume $\text{dim}\,V_{\mu}=\text{dim}\,V_{\mu}^{\nu}$ and that $\pml$ is large enough. Then from the proof of \cite[Theorem 4.1]{kim09} follows
\begin{equation}
	\hat{\delta} (V,V_{\pml})\leq C\,\text{length}(\gamma_{\lambda})\sup_{\lambda\in\gamma_{\lambda}}\|(\lambda-A)^{-1}\|\sup_{\lambda\in\gamma_{\lambda}}\|(\lambda-A_{\pml})^{-1}\|e^{-\alpha_1 \pml}
\end{equation}
for some $\alpha_1 >0$. Hence, the inequality $\hat{\delta} (V,V_{\pml}^{\nu})\leq 2(\hat{\delta} (V,V_{\pml})+\hat{\delta} (V_{\pml},V_{\pml}^{\nu}))$ \cite[IV.2]{Kato1980} implies
\begin{equation}\label{eq:gapPML}
	\hat{\delta} (V,V_{\pml}^{\nu})\leq C_0\,\text{length}(\gamma_{\lambda})\sup_{\lambda\in\gamma_{\lambda}}\|(\lambda-A_\pml)^{-1}\|\left (\sup_{\lambda\in\gamma_{\lambda}}\|(\lambda-A_{\pml}^{\nu})^{-1}\|\hat\rho^{\fem}_0+ C_1\sup_{\lambda\in\gamma_{\lambda}}\|(\lambda-A)^{-1}\| e^{-\alpha_1 \pml}\right).
\end{equation}

The estimates \eqref{eq:OsbornEst} and \eqref{eq:gapPML} show that the gap between the generalized eigenspace and the corresponding approximation depend critically on the norm of the resolvent. The non-normality of the operator and the associated finite element matrix are therefore important. 


\subsubsection{A DtN map for the finite PML problem}\label{sec:dtn_pml}

The eigenvalues of the DtN formulation of the resonance problem coincide with the eigenvalues of the infinite PML formulation provided that the PML strength $\sigma_0$ is large enough. However, the truncated PML problem may have more eigenvalues inside a given region in the complex plane compared with these exact formulations. Here, we derive a DtN-map for the finite PML problem that enables us to compare the standard DtN formulation with the finite PML problem in a new way. Moreover, the DtN map for the PML problem will be used to derive scalar eigenvalue relations that in Sect. \ref{sec:ex_sol_pml} are used to obtain highly accurate reference solutions. 

Let $\Omega_1:=(-\pml,-d)$, $\Omega_2:=(d,\pml)$. Then the problem with a finite PML layer \eqref{eq:trun_PML} can be replaced by the following coupled problem
\begin{subequations}\label{eq:trun_PML_split}
	\begin{align}
		-u''_0-k^2n^2 \,u_0=0 \quad  \text{in} \quad  & \Omega_{\dtn}, \label{eq:trun_PML_u0} \\
		-\dtot{}{x} \left(\frac{1}{\alpha}\dtot{u_j}{x}\right) 
			-k^2n^2_0\,\alpha \,u_j=0 \quad  \text{in} \quad  & \Omega_j,\,j=1,2, \label{eq:trun_PML_u12}
	\end{align}	
\end{subequations}
with boundary conditions $u_1(-\pml)=u_2(\pml)=0$ and the compatibility conditions
\begin{equation}\label{eq:compat_u}
		u_0(-d)=u_{1}(-d),\quad u'_0(-d)=u'_{1}(-d),\quad u_0(d)=u_{2}(d),\quad  u'_0(d)=u'_{2}(d). 
\end{equation}

In the following discussion we use the definitions
\begin{equation}
	\beta:=n_0(\pml-a)(1+i\sigma_\pml),\,\, \sigma_\pml:=\sigma_0\frac{\pml-\hat x}{\pml-a},\,\,
	\hat x=\fr{\dtn+x_c}{2},\,\,\hbox{and}\,\,
\phi(k):=\left(\frac{1+e^{2ik\beta}}{1-e^{2ik\beta}} \right).
	\label{eq:def_param}
\end{equation}
Define the sets $\mathcal{D}^+:=\{k\in\C\,:\,1+e^{2ik\beta}=0 \}$,
$\mathcal{D}^-:=\{k\in\C\,:\,1-e^{2ik\beta}=0 \}$ and $\mathcal{D}_c:=\mathcal{D}^+\cup \mathcal{D}^-$.
From the explicit expressions for $u_1, u_2$, and \eqref{eq:compat_u} we obtain $u_0(-d)=u_0(d)=0$ for $k\in \mathcal{D}^-$
and $u'_0(-d)=u'_0(d)=0$ for $k\in \mathcal{D}^+$. Hence, for all $k\in\mathcal{D}_c$ the equations \eqref{eq:trun_PML_u0} and \eqref  {eq:trun_PML_u12} decouple.
Values $k\in \mathcal{D}_c$ cannot be solutions of \eqref{eq:trun_PML_u0} because the eigenvalue problem \eqref{eq:trun_PML_u0} with homogeneous Dirichlet or Neumann boundary conditions has only real eigenvalues $k$, but $1\pm e^{2ik\beta}=0$ has no other real solutions besides the trivial solution. 

In the case, $k\notin \mathcal{D}_c$, the solutions of the coupled problem 
\eqref{eq:trun_PML_split}-\eqref{eq:compat_u} in $\Omega_{\dtn}$ are equivalent to \eqref{eq:trun_PML_u0} with the DtN-map 
\begin{equation}\label{eq:PML_DtN}
u'_0(-\dtn)=-i k \,n_0\phi(k) u_0(-\dtn) \,\,\,\hbox{and}\,\,\,
u'_0(\dtn) =i k \,n_0\phi(k) u_0(\dtn).
\end{equation}
The condition $|e^{2ik\beta}|< 1$ ensures that waves decay exponentially in the PML region, then we define the \emph{critical line} as the subset of $\mx C$ such that
\begin{equation}\label{eq:CriticalLine}
	\arg {k} =\arg\left (\frac{1}{1+i\sigma_\pml}\right ).
\end{equation}
We use \eqref{eq:CriticalLine} to divide the 4th quadrant of the complex plane into \emph{feasible} and \emph{nonfeasible} searching regions for resonances. In the infinite PML problem it is known \cite{co96,kim09}, that the sector is limited by $\arg {k} =\arg (1/(1+i\sigma_0))$, which is a larger sector compared with the sector defined by \eqref{eq:CriticalLine}.

Set $\lambda=-\iu k$ and define for $u,v \in\Hone$, and $\lambda\in \mathcal{D}:=\{\lambda\in\C\,:\, e^{-2\beta\lambda}\neq 1\}$ the form-valued function
\begin{equation}\label{eq:pmlDtN_form}
	\mathfrak{t}(\lambda)[u,v]:=\lambda^{2} \sa_{2}[u,v]+\lambda \sa_{1}[u,v]+\sa_{0}[u,v]+g(\lambda)\sa_{1}[u,v],
\end{equation}
where $\sa_{n}$, $n=0,1,2$ are defined in \eqref{a0}, \eqref{a1}, \eqref{a2}, and $g(\lambda)=2\lambda e^{-2\beta\lambda}/(1-e^{-2\beta\lambda})$.
Define as in \eqref{weak_dtn} the shifted form $\hat{\mathfrak{t}}(\lambda):=\mathfrak{t}(\lambda+\alpha)$, $\alpha>0$, and let $\hat g(\lambda)=g(\lambda+\alpha)$. Let $\hat\sa_{0}[A_1 u,v]=\sa_{1}[u,v]$ for all 
 $u,v \in\Hone$.  The operator function $T$ corresponding to the shifted problem
\begin{equation}\label{eq:ShiftedpmlDtN_form}
	\hat{\mathfrak{t}}(\lambda)[u,v]:=\lambda^{2} \hat\sa_{2}[u,v]+\lambda \hat\sa_{1}[u,v]+\hat\sa_{0}[u,v]+\hat g(\lambda)\sa_{1}[u,v]
\end{equation}
is then
\begin{equation}\label{eq:ShiftedpmlDtN_form}
	T(\lambda):=Q(\lambda)+\hat g(\lambda)A_1,
\end{equation}
where $Q$ is defined in \eqref{eq:Q}. Hence, $T$ is a finite rank perturbation of $Q$ and $\hat g(\lambda)\rightarrow 0$, when $\sigma_\pml\rightarrow\infty$. We derive below an estimate of the distance between $\lambda\in\sigma (Q)$ and $\tilde\lambda\in\sigma (T)$ for simple eigenvalues.
Assume that $v$ and $\tilde v$ are right eigenvectors of $Q(\lambda)$ and $T(\tilde{\lambda})$, respectively. Let  $w$ and $\tilde w$ denote the corresponding left eigenvectors with the normalizations $w^{*} Q'(\lambda) v=1$ and $\tilde{w}^{*} T'(\tilde\lambda) \tilde{v}=1$.  For a simple eigenvalue $\lambda\in\sigma (Q)$ there exists a neighbourhood $\mathcal{N}$ containing $\lambda$ such that
\begin{equation}\label{eq:Keldysh}
	Q^{-1}(z)=\frac{vw^{*}}{z-\lambda}+R(\lambda),
\end{equation}
where $R$ is analytic on $\mathcal{N}$ \cite{Kel1971}. Let $\gamma\subset \mathcal{N}$ be a Cauchy contour around $\lambda$ and define
\begin{equation}\label{eq:A}
	A^{(0)}:=\frac{1}{2\pi i}\int_\gamma Q^{-1}(z)dz=vw^*,\quad A^{(1)}:=\frac{1}{2\pi i}\int_\gamma zQ^{-1}(z)dz=\lambda vw^* .
\end{equation}
The resolvent of $T$ close to a simple eigenvalue can also be represented in the form \eqref{eq:Keldysh}. Assume that it exists exactly one simple eigenvalue $\tilde\lambda\in\sigma (T)$ in the neighbourhood $\mathcal{N}$ of $\lambda\in\sigma (Q)$ and define 
\begin{equation}\label{eq:tildeA}
	\tilde A^{(0)}:=\frac{1}{2\pi i}\int_\gamma T^{-1}(z)dz=\tilde v\tilde w^*,\quad\tilde  A^{(1)}:=\frac{1}{2\pi i}\int_\gamma zT^{-1}(z)dz=\tilde\lambda \tilde v\tilde w^* .
\end{equation}
Set $A^{(n)}_u:=u^* A^{(n)} u$ and  $\tilde A^{(n)}_u:=u^* \tilde A^{(n)} u$ for $n=0,1$ and let $u\in\Hone$ denote a function such that $A^{0}_u\neq 0$ and  $\tilde A^{0}_u\neq 0$. Then, the following estimate holds
\[
	|\lambda-\tilde\lambda |\leq\frac{1}{| A^{(0)}_u |}\left (| A^{(1)}_u-\tilde A^{(1)}_u |+|\tilde\lambda | | A^{(0)}_u-\tilde A^{(0)}_u | \right ).
\]
From the identity $Q^{-1}-T^{-1}=\hat g(\lambda)T^{-1}A_1Q^{-1}$ and the integral representations \eqref{eq:A}, \eqref{eq:tildeA} follows
\[
	|\lambda-\tilde\lambda |\leq C\left (\max_{z\in\gamma}|z\hat g(z)|+|\tilde\lambda |\max_{z\in\gamma}|\hat g(z)|\right )
\]
for some positive constant $C$. Note that $| \lambda-\tilde\lambda |\rightarrow 0$ when $\sigma_\pml\rightarrow\infty$ but $| \lambda-\tilde\lambda |$ can for $e^{-2\beta\lambda}\approx 1$ be large even if $\sigma_\pml$ is very large.






\section{Reference solutions}\label{sec:ex_sol_pml}
Assume that $\{\psi_1(x,k),\psi_2(x,k)\}$ for a given $n(x)$ are two independent solutions of \eqref{eq:Helmholtz}. Then, with the DtN-map \eqref{eq:formalDtN} we find the implicit eigenvalue relation
\begin{equation}
	 \fr{\psi'_1( d,k) - ik \,\psi_1( d,k) }{\psi'_2( d,k) - ik \,\psi_2( d,k) }
	=\fr{\psi'_1(-d,k) + ik \,\psi_1(-d,k) }{\psi'_2(-d,k) + ik \,\psi_2(-d,k) }.
	\label{eq:general_solution}
\end{equation}
Similarly, the DtN-map for the finite PML formulation \eqref{eq:trun_PML_u0}-\eqref{eq:PML_DtN} and $k\notin \mathcal{D}_c$, gives the implicit eigenvalue relation 
\begin{equation}
	 \fr{\psi'_1( d,k) - ik \,\phi(k)\psi_1( d,k) }{\psi'_2( d,k) - ik \,\phi(k)\psi_2( d,k) }
	=\fr{\psi'_1(-d,k) + ik \,\phi(k)\psi_1(-d,k) }{\psi'_2(-d,k) + ik \,\phi(k)\psi_2(-d,k) }.
	\label{eq:general_PML_solution}
\end{equation}
If $\{\psi_1(x,k),\psi_2(x,k)\}$ are known, the solutions of the scalar equations \eqref{eq:general_solution}, \eqref{eq:general_PML_solution} can be calculated to very high accuracy using a complex Newton root finder \cite{Yau98}. These values will be used as reference solutions for our finite element computations. Note that this procedure is slow since many computations with different initial guesses are necessary to compute all resonances in a given region of the complex plane.

In the remainder of the section we introduce three test cases used to compare results with the discussed formulations.

\subsection{Single slab problem}\label{sec:single_slab}
The following problem has been considered by several authors including \cite{osting13,kim09}. Define for $\eta\geq 1$ the piecewise constant function $n$ as
\begin{equation}\label{eq:ref_idx}
	n \left( x\right) =\left\{ 
	\begin{array}{lcr}
	\eta & \hbox{if} & |x|\leq \res\\
	1& \hbox{if} & |x|> \res
	\end{array} \right.
\end{equation}

\textit{Solution with DtN map:} Equation \eqref{eq:general_solution} reduces to $e^{-4i\eta k \res}=R^2$, where  $R:=(\eta -1)/(\eta +1)$ is called the reflectance. For $\eta\neq 1$ the solutions are 
\begin{equation}\label{eq:eigs_tm_1d}
	k_m=\frac{\pi m}{2\eta \res}-i\frac{\ln\left|1/R\right|}{2\eta \res},\,\,m=0,1,2,\ldots
\end{equation}
and for $\eta=1$ the equation has no solutions.

\textit{Solution with finite PML layer:}
Using $n(x)$ as in \eqref{eq:ref_idx} with $\eta\neq 1$, equation \eqref{eq:general_PML_solution} becomes
\begin{equation}
	e^{-4i\eta k \res}=\left( \frac{ \eta - \phi(k) }{\eta + \phi(k) }\right)^{2},\,\,\,
	\phi(k):=\left(\frac{1+e^{2ik\beta}}{1-e^{2ik\beta}} \right),
	\label{eq:eig_tm_pml_1d}
\end{equation}
which for given $\eta$ is solved numerically with a complex Newton root finder. 
For $\eta=1$ we obtain $e^{4ik(\beta-\res)}=1$, which has infinite many solutions 
\begin{equation}\label{eq:SlabPML}
	k_m=\tfrac{2m+1}{2(\beta-\res)}\pi,\,\,m=0,1,2,\ldots.
\end{equation}
The eigenvalues $k_m$ in \eqref{eq:SlabPML} are close to the critical line defined by \eqref{eq:CriticalLine}. Note that the \eqref{eq:SlabPML} are exact solutions of the finite PML formulation but the problem with the usual DtN-map has no eigenvalues. 
Moreover,  $|k_{m+1}-k_m|\rightarrow 0$ and $\arg k_m\rightarrow \arg(1/(1+i\sigma_0))$ when  $\pml\rightarrow\infty$.

For a general $n(x)$, we expect that the finite PML formulation has more solutions than the formulation with a DtN-map and these additional solutions are called spurious solutions.

\subsection{Air-filled-cavity problem}\label{sec:air_filled}

The single slab problem was studied in \cite{ket12} and spurious eigenvalues were successfully detected with a perturbation approach. In this section, we present a more demanding problem where the perturbation approach only selects a very small region of the complex plane as free of spurious solutions as indicated in \cite[Figure 5.24]{ket12}. In particular, we derive a scalar equation for the eigenvalues, which enable us to compute highly accurate reference solutions with a Newton root finder.

Define for  $\kett>1$ the refractive index
\begin{equation}
	 n \left( x\right) =\left\{ 
	\begin{array}{lcll}
	1   & \hbox{if} & &\! |x|\leq 1\\
	\gamma & \hbox{if} & 1<&\! |x|\leq \kett\\
	\eta & \hbox{if} & \kett< &\! |x|
	\end{array} \right..
\label{eq:tm_ket_1d}
\end{equation}
From \eqref{eq:general_solution}, we derive the implicit eigenvalue relation
\begin{equation}
\begin{array}{ll}
e^{-4ik} \!\!\!&\!\!\left(\frac{
		(1\!+\!\eta/\gamma)(1\!+\!\gamma)\,e^{ik(\kett(\eta-\gamma)+\gamma)}+
		(1\!-\!\eta/\gamma)(1\!-\!\gamma)\,e^{ik(\kett(\eta+\gamma)-\gamma)} }{
		(1\!+\!\eta/\gamma)(1\!-\!\gamma)\,e^{ik(\kett(\eta-\gamma)+\gamma)}+
		(1\!-\!\eta/\gamma)(1\!+\!\gamma)\,e^{ik(\kett(\eta+\gamma)-\gamma)} 
}\right)\\[3mm]
 & =  \left(\frac{
		(1\!-\!\eta/\gamma)(1\!+\!\gamma)\,e^{ ik(\kett(\gamma-\eta)-\gamma)}+
		(1\!+\!\eta/\gamma)(1\!-\!\gamma)\,e^{-ik(\kett(\gamma+\eta)-\gamma)}}{
		(1\!-\!\eta/\gamma)(1\!-\!\gamma)\,e^{ ik(\kett(\gamma-\eta)-\gamma)}+
		(1\!+\!\eta/\gamma)(1\!+\!\gamma)\,e^{-ik(\kett(\gamma+\eta)-\gamma)}
}\right).
\end{array}
	\label{eq:eigs_tm_kett_1d}
\end{equation}

\begin{table}
\robustify\bfseries
\centering
\begin{tabular}{rS[table-format=2.10, detect-weight]
                 S[table-format=2.10, detect-weight]
                 rS[table-format=2.10, detect-weight]
                 S[table-format=2.10, detect-weight] }
\toprule
{$j$} & {$\re k_j$} & {$\im k_j$} & {$j$} & {$\re k_j$} & {$\im k_j$} \\
\midrule
 0 & 	0.0000000000 & -0.8948801287 &  8 & 	6.6087515863 & -0.8788560394 \\
 1 & 	0.4869949494 & -0.6502632860 &  9 & 	7.0248667636 & -0.7730423533 \\
 2 & 	1.5955486049 & -0.3950551466 & 10 & 	7.9794721839 & -0.4166038034 \\
 3 & 	2.7503593706 & -0.5843773974 & 11 & 	9.1753687526 & -0.4808796847 \\
 4 & 	3.3047923378 & -0.8909296467 & 12 & 	9.9108347715 & -0.8579829521 \\
 5 & 	3.7465666834 & -0.7159810538 & 13 &  10.3153076002 & -0.8180915326 \\
 6 & 	4.7869777032 & -0.4021092410 & 14 &  11.1740110180 & -0.4393352673 \\
 7 & 	5.9689601644 & -0.5268047778 & 15 &  12.3746790920 & -0.4461923754 \\
\bottomrule
\\
\end{tabular}


\vspace*{-5mm}
\caption{\emph{Selected reference eigenvalues for the air-filled-cavity problem (Sect. \ref{sec:air_filled}) ordered by $|\Re k_j|$.}}
\label{tab:air_filled_reference}
\end{table}
In our finite element calculations we used the profile \eqref{eq:tm_ket_1d} with $\kett=1.5,\gamma=\sqrt{3.5}$, and $\eta=\sqrt{2.5}$. A selection of eigenvalue approximations computed from \eqref{eq:eigs_tm_kett_1d}
are given in Table \ref{tab:air_filled_reference}. In these computations, we use the complex Newton root finder \cite{Yau98} with machine precision as stopping-criteria.  

\subsection{Bump problem}\label{sec:bump} 
\begin{table}
	\robustify\bfseries
	\centering
	\begin{tabular}{rS[table-format=2.10, detect-weight]
                 S[table-format=2.10, detect-weight]
                 rS[table-format=2.10, detect-weight]
                 S[table-format=2.10, detect-weight] }
\toprule
{$j$} & {$\re k_j$} & {$\im k_j$} & {$j$} & {$\re k_j$} & {$\im k_j$} \\
\midrule
 0 &  0.0000000000 & -0.4271986734 &  6 &  6.0034893253 & -0.7920181369  \\
 1 &  1.1402018812 & -0.4825101535 &  7 &  6.9572111153 & -0.8281487827  \\
 2 &  2.1432843061 & -0.5771518110 &  8 &  7.9089927230 & -0.8604952505  \\
 3 &  3.1204984325 & -0.6473255266 &  9 &  8.8593105049 & -0.8897868318  \\
 4 &  4.0868340691 & -0.7036943333 & 10 &  9.8084919100 & -0.9165558262  \\
 5 &  5.0470974941 & -0.7510601464 & 11 & 10.7567710490 & -0.9412039599  \\
\bottomrule
\\
\end{tabular}


	\vspace*{-5mm}
	\caption{\emph{Selected reference eigenvalues for the refractive index \eqref{eq:tm_bump} ordered by $|\Re k_j|$.}}
	\label{tab:bump_reference}
\end{table}

Problems with continuous refractive index are interesting from the application point of view \cite{BAYLISS1985,Wadbro2008,Gillman2015,Shiri201623}. Motivated by the application in \cite{Shiri201623}, we introduce a refractive index profile resembling a continuous \emph{bump}, and define
\begin{equation}
	 n \left( x\right) =\left\{ 
	\begin{array}{lcll}
	2-x^2  & \hbox{if} & &\! |x|\leq 1\\
	1  & \hbox{if} & &\! |x| > 1
	\end{array} \right..
\label{eq:tm_bump}
\end{equation}

The solutions $\psi_1,\psi_2$ in 
\eqref{eq:general_solution}, are not available for this problem, and we therefore compute reference solutions by using a very fine FE discretization.

\section{The Lippmann-Schwinger equation and pseudospectrum}\label{sec:resolvent}

For a closed linear operator $T$ on a Banach space we denote by $\text{Ker}\, T$, $\text{Ran}\, T$, $\rho (T)$, and $\sigma (T)$, its kernel, range, resolvent set, and spectrum, respectively. Let $A_0:L^2(\R)\to L^2(\R)$, $A_0=-n_0^{-2}D_x^2$, denote the Laplacian times a constant $-n_0^{-2}$ with domain $D(A_0)=H^2(\R)$. 

For $\Im k>0$ the resolvent $R_0(k): L^2(\R)\to L^2(\R)$ of $A_0$ is
\[
	R_0(k)u:=(A_0-k^2)^{-1}u=\frac{i n_0}{2k}\int_{\R}e^{i n_0 k|x-y|}u(y)dy.
\]
Let $L^2_c(\R)$ denote the space of $L^2$-functions with compact support and let $L^2_{\mrm{loc}}(\R)$ denote the space of functions such that the restriction to every bounded subset $\Omega$ of $\R$ lies in $L^2(\Omega)$. The resolvent operator $R_0$ extends for $\Im k<0$ to a meromorphic family of operators $R_0(k): L^2_c(\R)\to L^2_{\mrm{loc}}(\R)$ \cite{LectureZworski}. Assume that $\chi_{\dtn}$ is a $C^{\infty}$ function with support in $\Omega_\dtn:=(-\dtn,\dtn)$. Then
\begin{equation}\label{eq:resAzero}
	\|\chi_{\dtn} R_0(k)\chi_{\dtn}\| \leq C\frac{e^{2 n_0 \dtn(\Im k)_{-}}}{|k|},
\end{equation}
where $x_{-}:=\max (0,-x)$ and $C>0$  \cite[Theorem 2.1]{LectureZworski}. Let $n\in L^{\infty}(\R)$ denote the refractive index and assume that $n(x)>0$ and $n_0>0$. Define as above the operator $A=-n^{-2}D_x^2$ and its resolvent $R(k)=(A-k^2)^{-1}$. The resolvent operator $R(k): L^2(\R)\to L^2(\R)$, $\Im k^2>0$ extends to a meromorphic family of operators $R(k): L^2_c(\R)\to L^2_{\mrm{loc}}(\R)$, $k\in\C$ and the poles of $R$ are called resonances \cite{MR1119196}, see also \cite{MR1037774,MR1350074}. 
The identity $a^{-1}-b^{-1}=b^{-1}(b-a)a^{-1}$ gives
\[
	R(k)-R_0(k)\frac{n^{2}}{n_0^2}=k^2R_0(k)\frac{n^{2}-n_0^2}{n_0^2}R(k),\quad k\in\rho (A)\cap\rho (A_0).
\]
Hence
\begin{equation}\label{eq:LS}
	R(k)=T^{-1}(k)R_0(k)\frac{n^2}{n_0^2}, \quad T(k):=1-k^2 R_0(k)\frac{n^{2}-n_0^2}{n_0^2}
\end{equation}
and $k$ is a pole of $R$ if the pair $(u,k)$ satisfies the Lippmann-Schwinger equation  
\begin{equation}
	T(k)u=0.
	\label{eq:lipp_sch_eq}
\end{equation}
The non-linear eigenvalue problem \eqref{eq:lipp_sch_eq} has successfully been used to compute resonances \cite{osting13} and the eigenvalues of $T$ equal the resonances of $A$ (comp. e.g. \cite[Chapter 2]{LectureZworski}). Let $\chi_{d}$ denote a $C^{\infty}$ function with support in  $\text{supp}\, (n^2-n_0^2)\subset (-\dtn,\dtn)$. From \eqref{eq:LS} follows
\[
	T^{-1}(k)=I+k^2R(k)\frac{n^{2}-n_0^2}{n^2}
\]
and the following resolvent estimate holds:
\begin{equation}\label{eq:ResEst}
	\|\chi_{d}T^{-1}(k)\chi_{d}\|\leq 1+|k|^2\max \left | \frac{n^{2}-n_0^2}{n^2} \right | \| \chi_{d} R(k)\chi_{d}\|.
\end{equation}
We define for $\epsilon>0$ the restricted $\epsilon$-pseudospectrum $\sigma_{\epsilon}(T)$ as the set of all $k\in\C$ such that $\|\chi_{d}T^{-1}(k)\chi_{d}\|>\epsilon^{-1}$.Then, the restricted $\epsilon$-pseudospectrum $\sigma_{\epsilon}(A)$ is the set of all $k\in\C$ such that
\begin{equation}\label{RNorm}
	\|\chi_{\dtn} R(k)\chi_{\dtn} \|>\epsilon^{-1},
\end{equation}
where the norm \eqref{RNorm} increases with $\dtn$ and we expect exponential growth. In particular, for $n=n_0$ the norm of $\chi_{\dtn} R(k)\chi_{\dtn}$ grows exponentially with $d$  \cite[Theorem 2.2]{LectureZworski} but $\|\chi_{d}T^{-1}(k)\chi_{d}\|=1$. Hence, it is plausible that the norm of $\chi_{d}T^{-1}(k)\chi_{d}$ is much smaller than the norm of $\chi_{\dtn} R(k)\chi_{\dtn}$. Moreover, in Sect. \ref{sec:res_air_filled} we consider the pseudospectrum of a discretized Lippmann-Schwinger operator and the numerical calculations suggest that the resolvent norm is very well behaved.

An equivalent condition for $k\in\sigma_{\epsilon}(T)$ is that it exist a normalized function $u$ for which 
\begin{equation}\label{IntegralTest}
	\|\chi_{d}T(k)\chi_{d}u\|<\epsilon.
\end{equation} 
Such $u$ is called an approximate eigenvector or $\epsilon$-pseudomode \cite[p. 255]{MR2359869} and we will in the paper check if the numerically computed functions are $\epsilon$-pseudomodes of $T(k)$. The main result of the paper is that \eqref{IntegralTest} can be used as a suitable 
measure to separate true eigenvalues from the spurious eigenvalues.



\begin{figure}
	\begin{tikzpicture}[thick,scale=1.0, every node/.style={scale=0.9}]
		\draw(  0.00, 0.0) node {
			\includegraphics[scale=0.43]{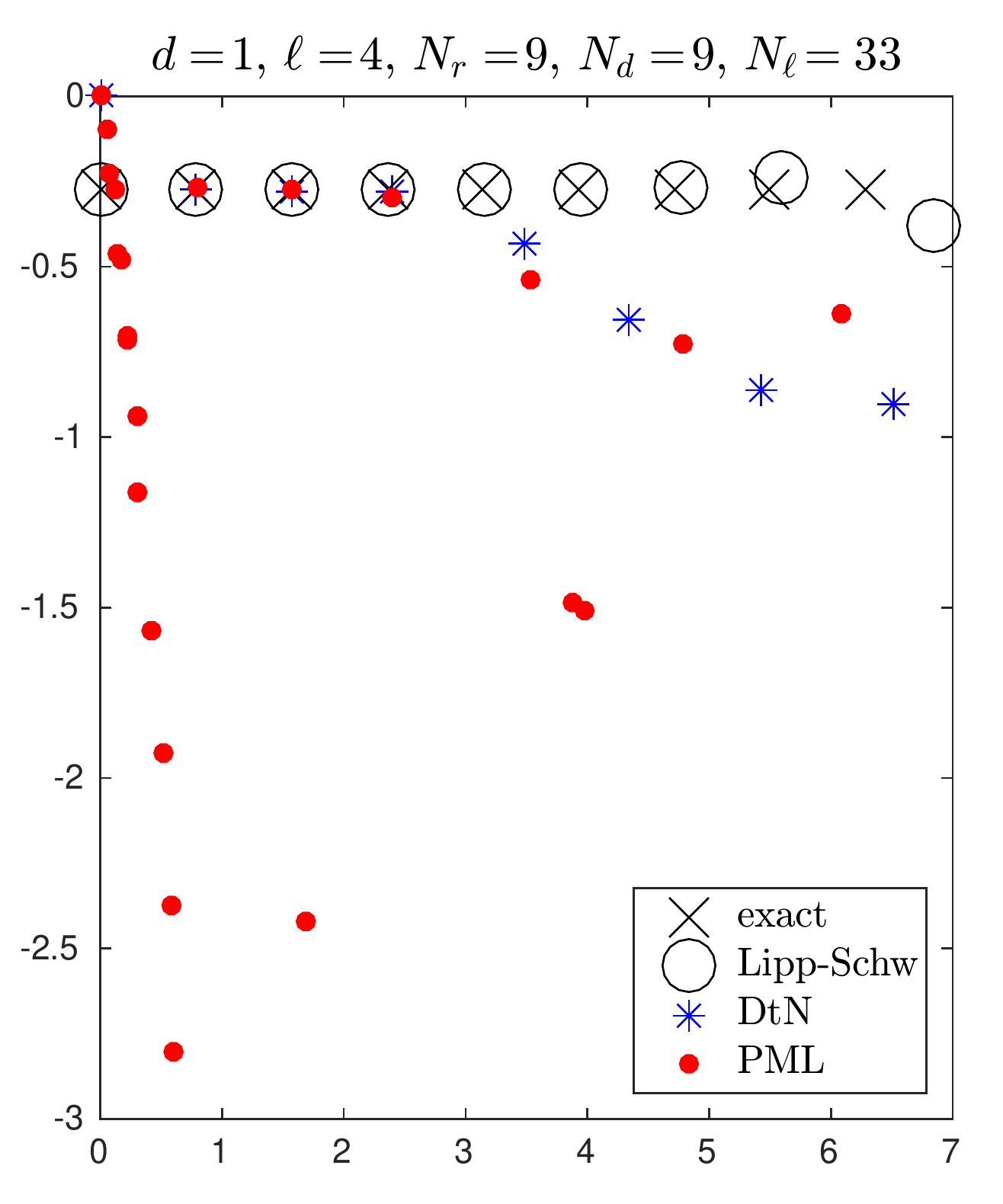} };
		\draw(  5.00, 0.0) node {	\includegraphics[scale=0.43]{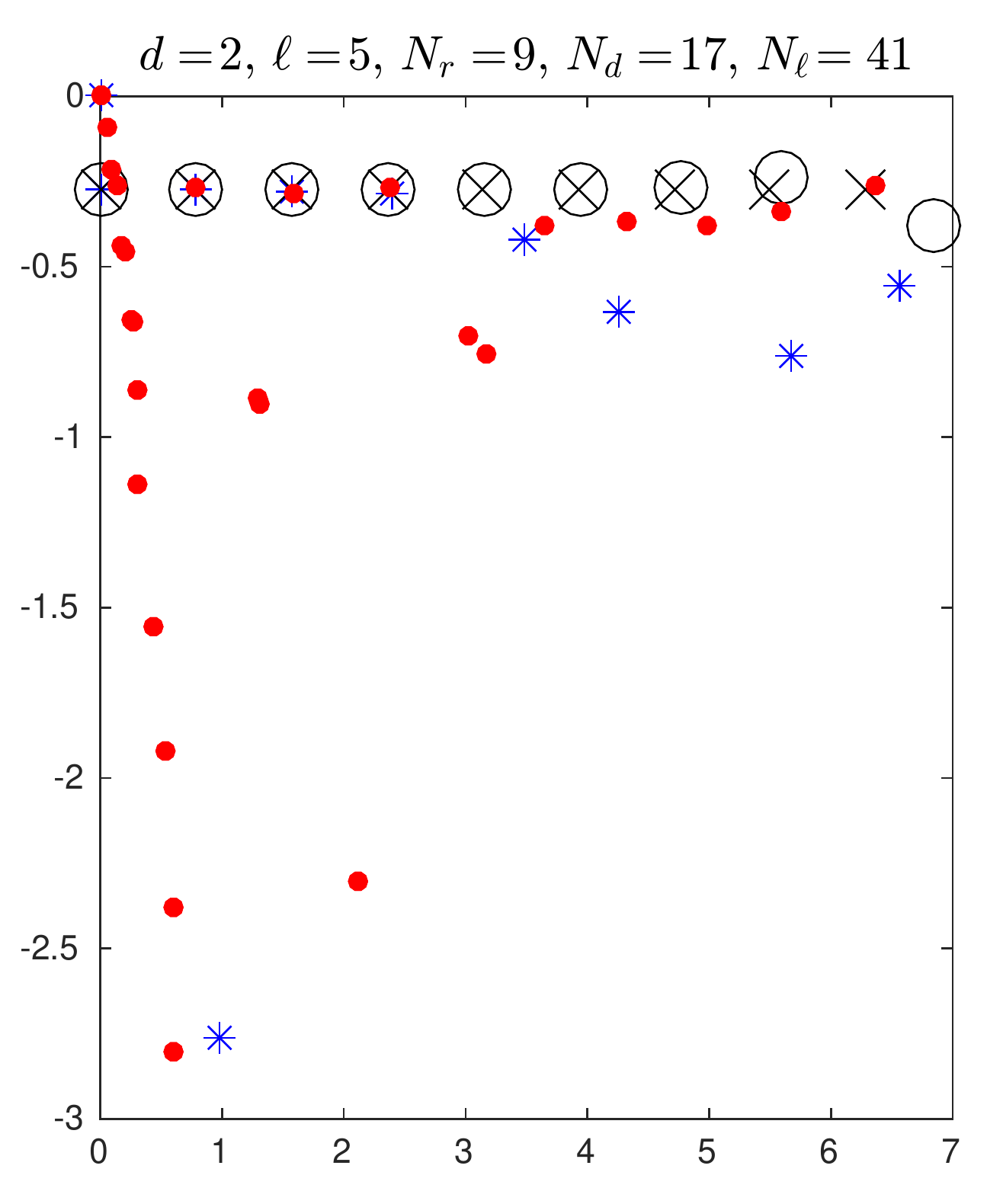} };
		\draw( 10.00, 0.0) node {	\includegraphics[scale=0.43]{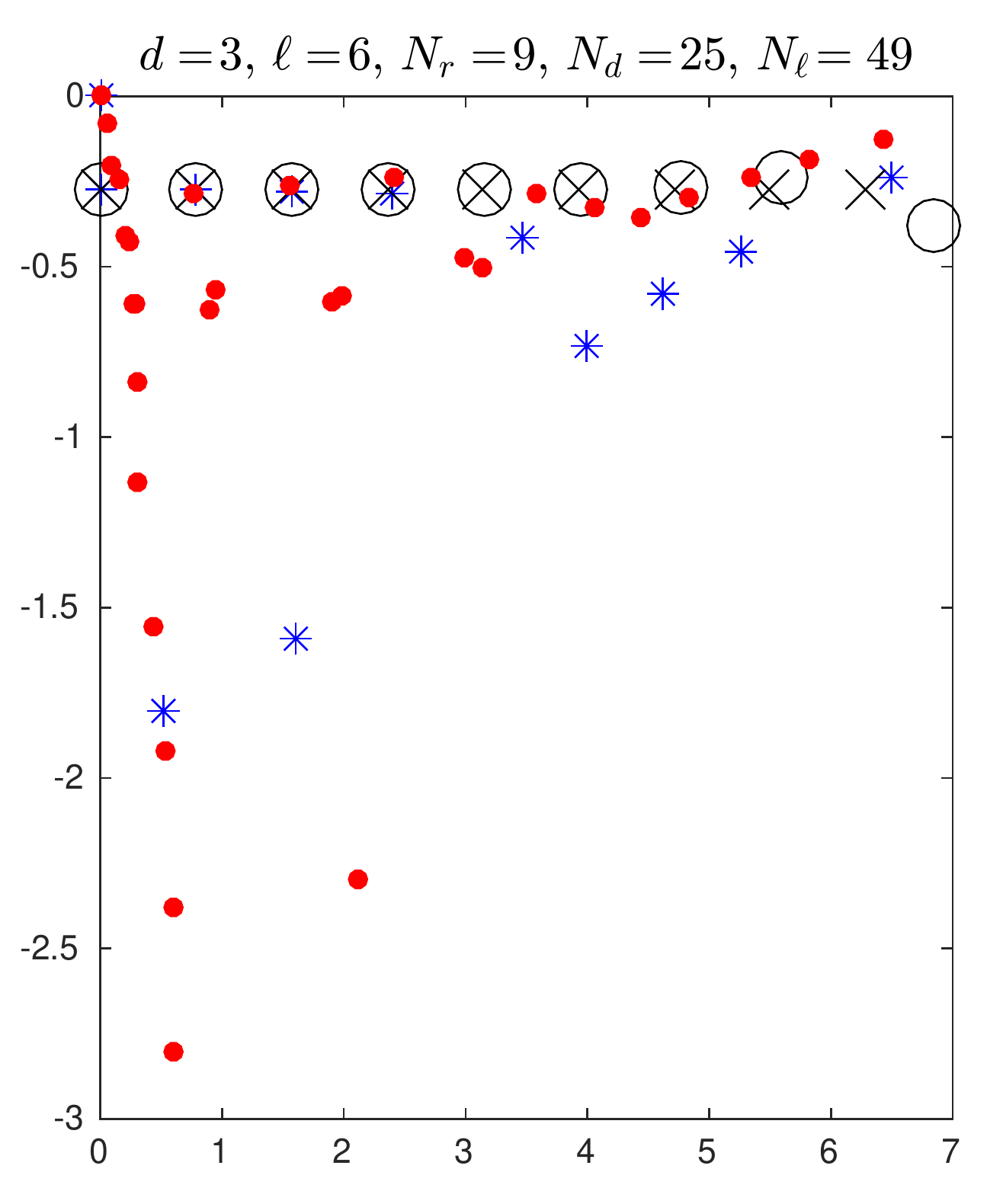} };
		
		\node[text width=3cm] at ( -1.00, 2.95) {$a)$};
		\node[text width=3cm] at (  4.00, 2.95) {$b)$};
		\node[text width=3cm] at (  9.00, 2.95) {$c)$};
	\end{tikzpicture}
	\vspace*{-3mm}
	\caption{\emph{Comparison of exact eigenvalues \eqref{eq:eigs_tm_1d} and approximations computed with the DtN-FEM \eqref{eq:dtn_discrete}, finite PML-FEM \eqref{eq:mat_pml} and Lippmann-Schwinger \eqref{eq:lipp_collocation} for a coarse discretization. In Panel $a$) both the DtN and PML are placed in the boundary of the resonator: $a=1$. In $b$) these are placed at $d=2$ and in $c$) at $d=3$, allowing some air in between the resonator and the truncation of the domain. In the computation we used $p=2$ and $\sigma_0=5$. }}
	\label{fig:comp_eigs_single_slab}
\end{figure}
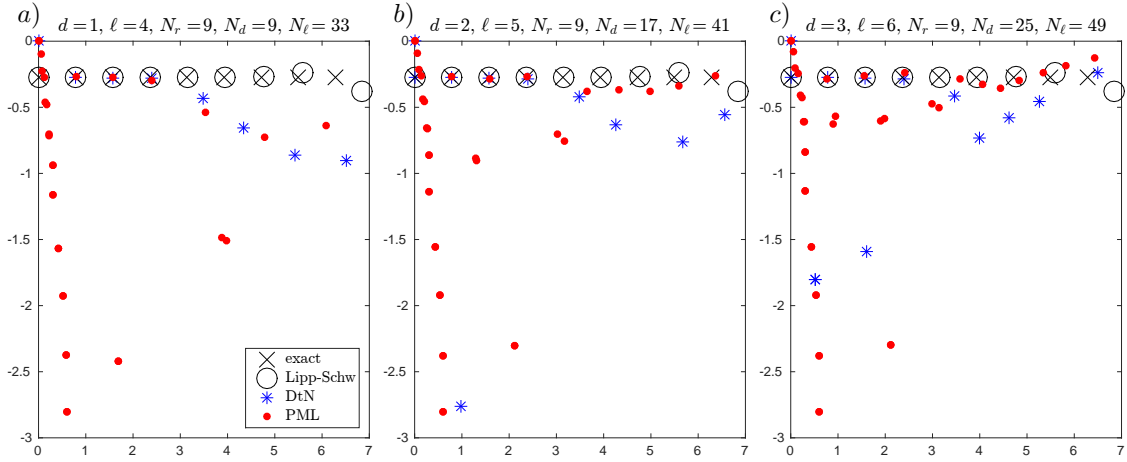

\section{Finite element discretization}\label{sec:fem}

A conforming finite element method is used to discretize the DtN and PML based formulations of the resonance problem. 

\textit{The DtN based formulation:} Let the open interval $\Omega_\dtn:=(-\dtn,\dtn)$, be covered with a regular and quasi uniform finite element mesh $\mc T$ consisting of elements $\{K_i\}^N_{i=1}$. Let $\rho_i$ be the length of the interval $K_i$ and denote by $h$ the maximum mesh size $h:=\max{\rho_i}$. Let $\mc P_p$ denote the space of polynomials on $\mx R$ of degree $\leq p$ and set $\fem:=\{h,p\}$. We define the finite element space $S^{\fem}(\Omega_\dtn):=\{u\in H^1(\Omega_\dtn):\left. u\right|_{K_i} \in\mc P_p(K_i)\,\,\hbox{for}\,\,K_i\in\mc T\}$, and $N_{\dtn}:=\dim(S^{\fem}(\Omega_\dtn))$ \cite[Ch 2]{Schwab1998}. Let $\{\varphi_1,\dots,\varphi_N\}$ be a basis of $S^{\fem}(\Omega_\dtn)$. From \eqref{eq:quadratic_form} we form a companion linearization similar to \eqref{eq:PLgen} and state the corresponding matrix problem: Find $\xi,\eta\in \mx C^{N_\dtn}$ and $\lambda^\fem\in\mx C$ such that
\begin{equation}\label{eq:dtn_discrete}
	\begin{bmatrix}
			A & E \\
			0 & I
		\end{bmatrix}
	\begin{bmatrix}
			\xi \\
			\eta
		\end{bmatrix}	
 =\lambda^{\nu}
	\begin{bmatrix}
		0 &-M \\
		I & \phantom{-}0
	\end{bmatrix}
	\begin{bmatrix}
			\xi \\
			\eta
		\end{bmatrix},
\end{equation}
with
\begin{equation}
	A_{ji}=\int_{-d}^d \varphi'_j\varphi'_i\,dx,\,\,M_{ji}=\int_{-d}^d n^2\varphi_j\varphi_i\,dx,\,\,
	E_{ji}=n_0 \left(\varphi_j(-d)\varphi_i(-d)+\varphi_j(d)\varphi_i(d)\right).
	\label{eq:mat_dtn}
\end{equation}
Then, we recover $k^\fem$ by the rule $k^\fem=-i\lambda^\fem$. 

\textit{The PML based formulation:}
Similarly to the DtN formulation, we define the finite element space $S_0^{\fem}(\Omega_\pml):=\{u\in H^1_0(\Omega_\pml):\left. u\right|_{K_i} \in\mc P_p(K_i)\,\,\hbox{for}\,\,K_i\in\mc T\}$, and $N_{\pml}:=\dim(S_0^{\fem}(\Omega_\pml))$. Let $\{\varphi_1,\dots,\varphi_N\}$ be a basis of $S_0^{\fem}(\Omega_\pml)$.
From $\sa_{\pml}[u,v] =\lambda \sb_{\pml}[u,v]$ as defined in \eqref{eq:defPML}, we formulate the corresponding matrix problem: Find $\xi\in \mx C^{N_\pml}$ and $\lambda^\fem\in\mx C$ such that
\begin{equation}
	\tilde A\xi =\lambda^\fem \tilde M\xi,\,\,\,\hbox{with}\,\,\, 
	\tilde A_{ji}=\int_{-d}^d \frac{1}{\alpha}\varphi'_j\varphi'_i\,dx,\,\,
	\tilde M_{ji}=\int_{-d}^d n^2\alpha \,\varphi_j\varphi_i\,dx.
	\label{eq:mat_pml}
\end{equation}
Then, we recover $k^\fem$ by the rule $k^\fem=\sqrt{\lambda^\fem}$. 

We can evaluate convergence rates for those eigenvalues that approximate resonances. All considered equations have piecewise analytic coefficients $n$, $\alpha$ in $\Omega_{\dtn}$ and $\Omega_{\pml}$. Hence, the eigenfunctions of the reference problems in Sect. \ref{sec:ex_sol_pml}, and the eigenfunctions of the corresponding adjoint problems are all piecewise analytic \cite{Babuska+Buo+Osborn1989}. Then, since meshes that respect the non-smoothness of the coefficients are used, exponential convergence is expected with $p$-FEM (fix $h$ and increase $p$), and optimal converge rates are expected with $h$-FEM (fix $p$ and decrease $h$) \cite{Babuska+Buo+Osborn1989,Schwab1998}. Approximation properties of spaces of piecewise polynomials has been extensively discussed in the literature, and we refer to \cite{Schwab1998,brenner+scott2002} for overviews and further details.

All meshes in our computations are uniform with coarse cells starting from a cell length $h=0.5$. In the presented simulations we refer to $q$ as an index for the formulation used: $q=d$ for the DtN, $q=\pml$ for the PML and $q=r$ for the Lippmann-Schwinger formulations. We use $N_q$ for the number of degrees of freedom, $p$ the polynomial degree of the shape functions $\varphi_j$, \emph{cells} the number of initial coarse cells $K_i$, and \emph{ref} the number of uniform refinements. 
Then, the relationship: $N_q=p\times$\emph{cells}$\times 2^{\hbox{\emph{ref}}}+1$, holds for each FE discretization in use. We approximate eigenpairs of 
the two problems \ref{eq:dtn_discrete} and \ref{eq:mat_pml} using \emph{deal.II} \cite{dealII82} for FE, \emph{PETSc} \cite{petscref} for the linear algebra routines and \emph{SLEPc} \cite{slepc05} for the eigenvalue solvers. The shape functions are based on Gauss-Lobatto shape functions \cite{solin04}. \\
\begin{figure}[!h]
	\begin{tikzpicture}[thick,scale=0.9, every node/.style={scale=0.9}]
		\draw( 0.00, 7.2)  node { \includegraphics[scale=0.45]{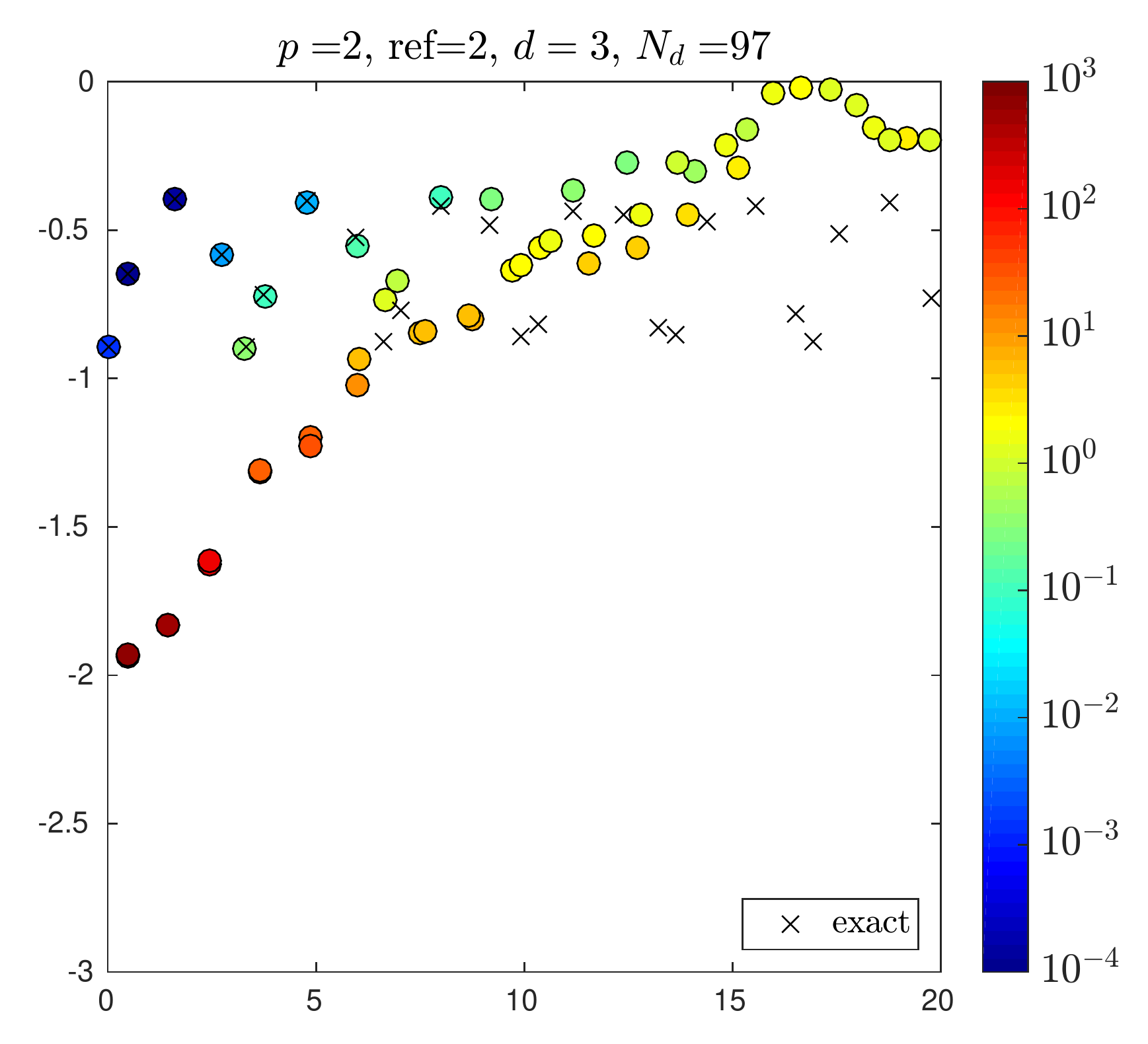} };
		\draw( 8.00, 7.2)  node { \includegraphics[scale=0.45]{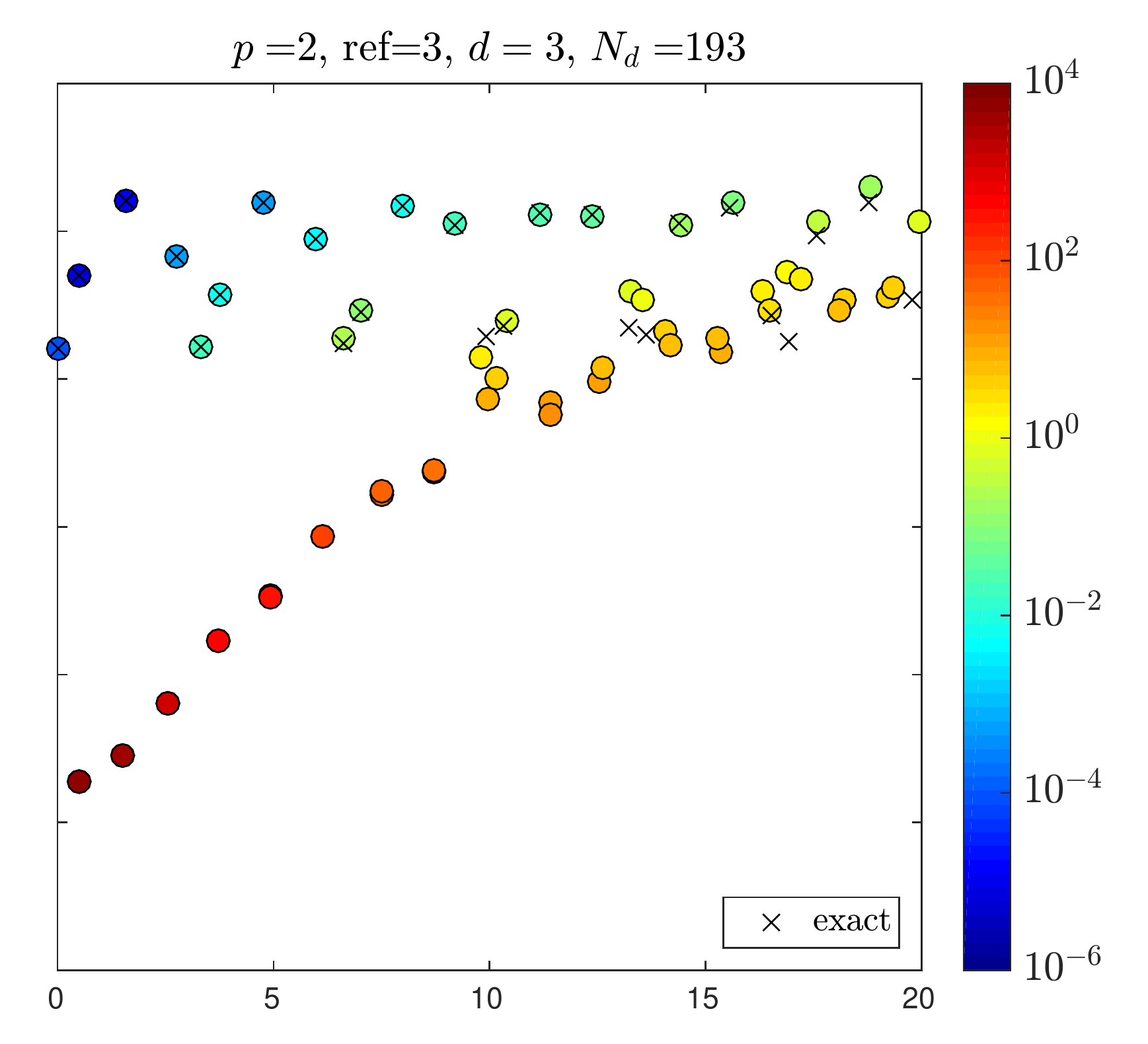} };
		
		\draw( 0.00, 0.0)  node { \includegraphics[scale=0.45]{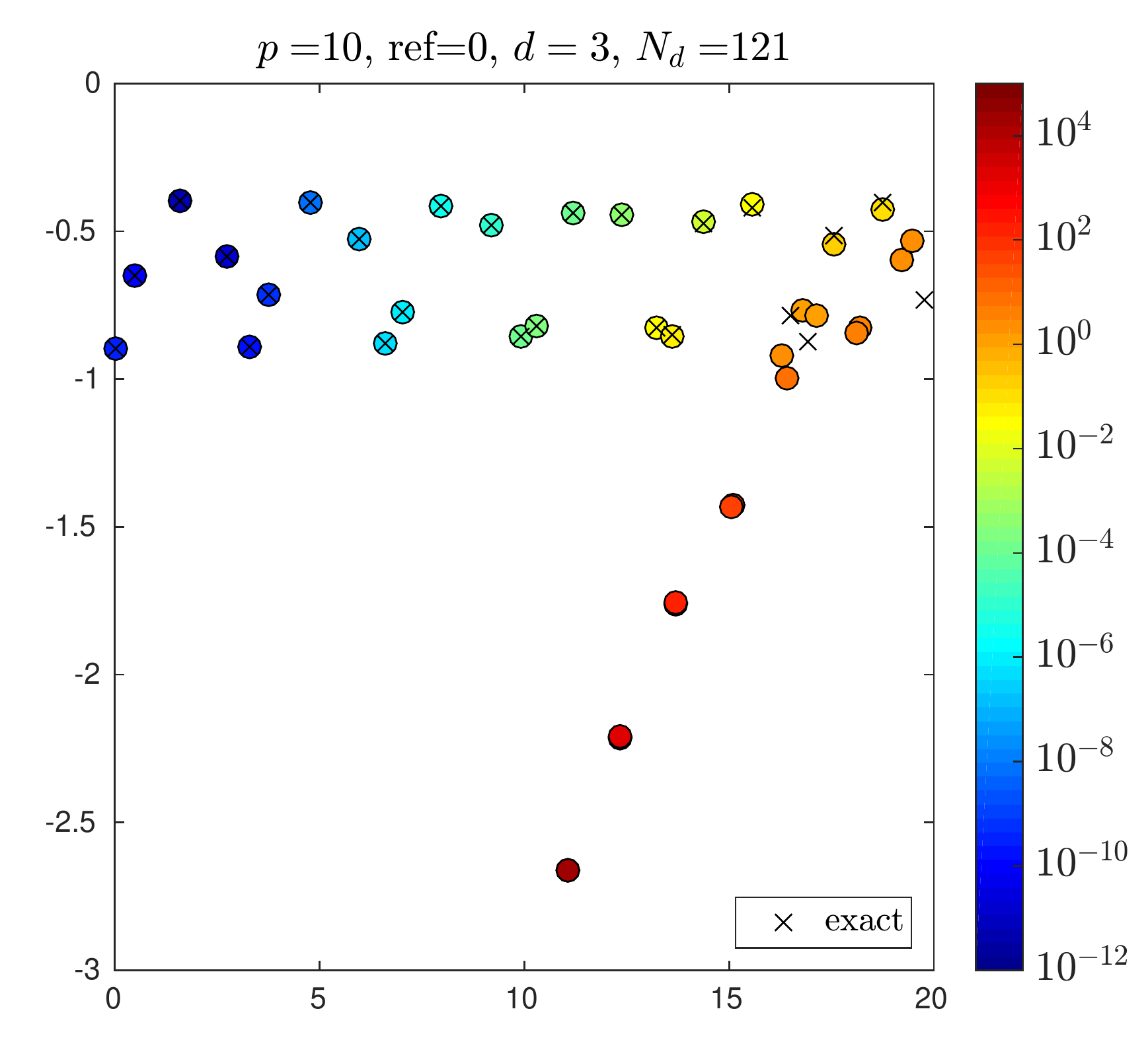} };
		\draw( 8.00, 0.0)  node { \includegraphics[scale=0.45]{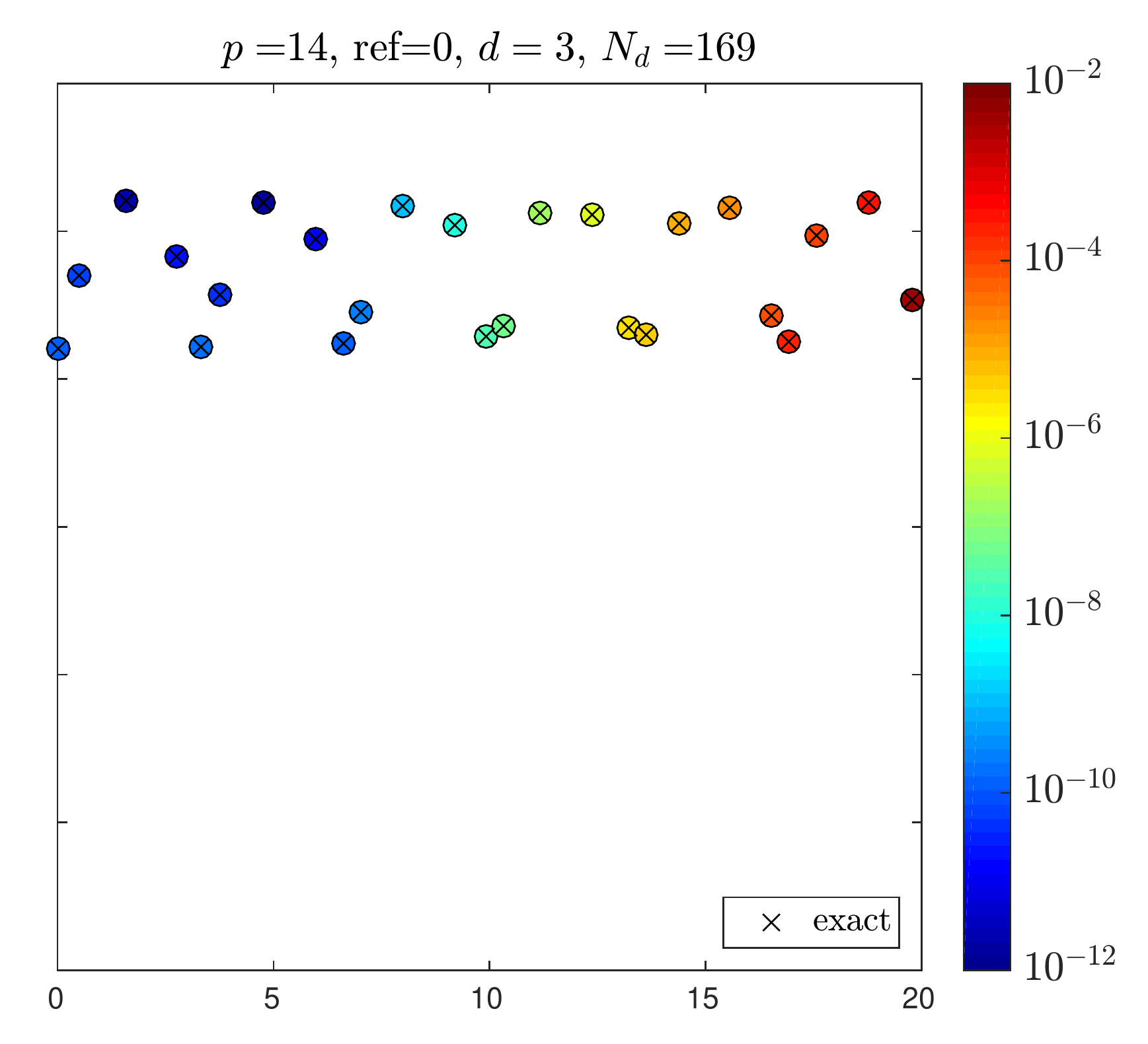} };
		
		\draw (0.4, 9.5) ellipse (0.21cm and 0.33cm);
		
		\node[text width=3cm] at (-1.60, 3.25) {$c)$};
		\node[text width=3cm] at ( 6.00, 3.25) {$d)$};
		
		\node[text width=3cm] at (-1.60, 10.45) {$a)$};
		\node[text width=3cm] at ( 6.00, 10.45) {$b)$};
		
	\end{tikzpicture}
	\vspace*{-5mm}
	\caption{\emph{Computed Eigenvalues $k_j^\fem$ of the DtN formulation \eqref{eq:dtn_discrete}, corresponding to the air-filled-cavity problem described in Sect.	\ref{sec:air_filled}. In colors we give $\epsilon_j$ computed with \eqref{eq:eps_check} for each $k_j^\fem$.} }
	\label{fig:i_check_DtN_air_filled}
\end{figure}
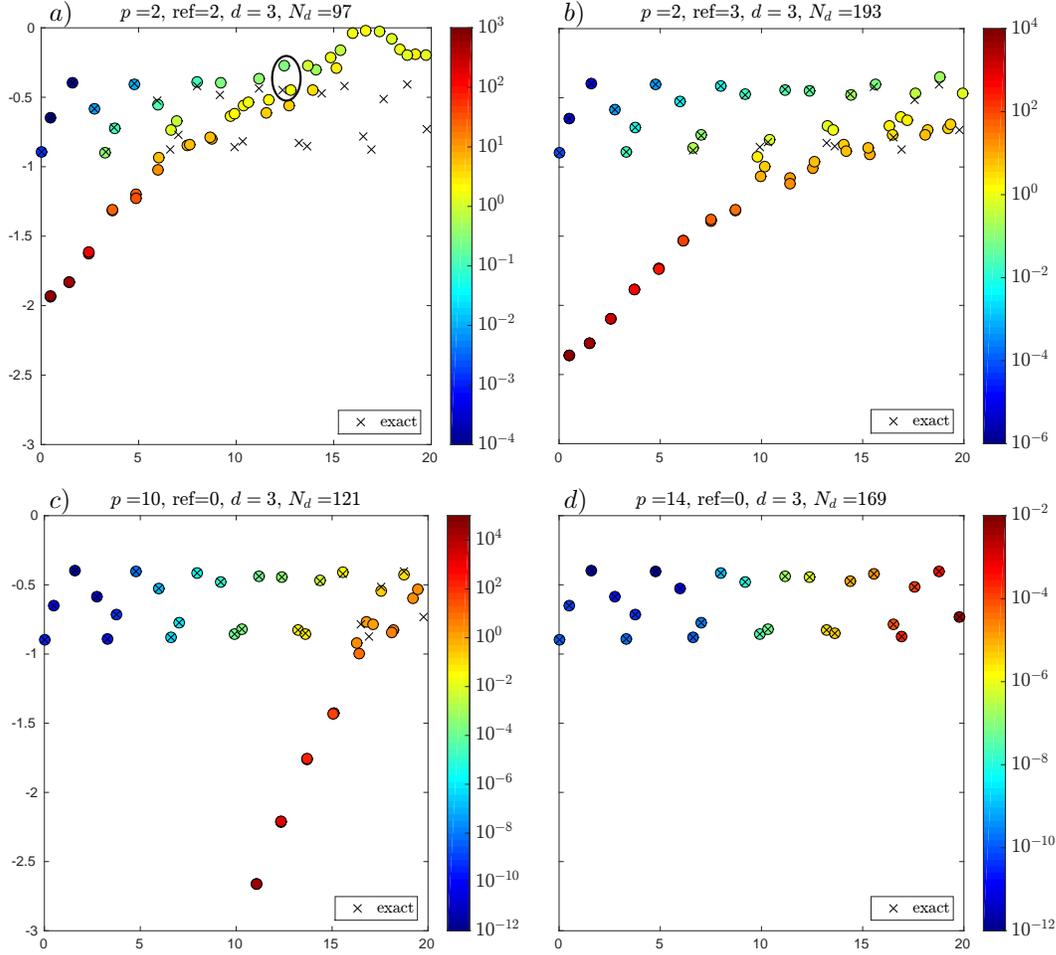
\noindent \textit{The Lippmann-Schwinger based formulation:}

The estimates (2.26) and (3.14) show that the gap between the generalized eigenspace and the corresponding approximation depend critically on the norm of the resolvent. Hence, the non-normality of the operator is important. However, the derivations assumed that the finite element space is large enough.

In Sect. \ref{sec:resolvent}, we argued that the Lippmann-Schwinger operator function $T$ has a well behaved resolvent and we expect therefore no spurious eigenvalues when resonances are computed from a numerical scheme based on \eqref{eq:lipp_sch_eq}.
To illustrate this, we present a collocation scheme referred to as case (A) of the Galerkin methods in \cite{ikebe72}, and used in \cite[Sect. 3.2]{gopal08} for resonance computations.
Let $\Omega_r:=\text{supp}\, (n^2-n_0^2)$, and let $\{\varphi_j\}$ be a basis for $S^\fem(\Omega_r)$ with the property $\varphi_j(x_i)=\delta_{ji}$, $\{x_i\}_{i=1}^{N_r}\in \Om_r$. 
Then, by plugging the ansatz $u^\fem=\sum^{N_r}_{j} \xi_j\varphi_j$ into \eqref{eq:lipp_sch_eq}, we obtain the nonlinear eigenvalue problem: Find $\xi\in\mx C^{N_r}$ and $k^\fem\in\mx C$ such that 
\begin{equation}
	T^\fem(k^\fem)\xi=(I-K(k^\fem))\xi=0\,\,\,\hbox{with}\,\,\, K_{ij}(k^\fem):=\frac{ik^\fem}{2n_0}\int_{\Omega_r}e^{in_0k^\fem |x_i-y|}(n^2(y)-n^2_0)\varphi_j(y)\, dy.
	\label{eq:lipp_collocation}
\end{equation}

This type of nonlinear matrix eigenvalue problems can be solved using a contour integration based method \cite{Sakurai09,MR2914550} and we apply NEPCISS to \eqref{eq:lipp_collocation}.
Perturbation estimates for eigenpairs of Fredholm valued functions show how the convergence rates of this type of methods are connected with the spectral properties of $T$ \cite{MR3399077,MR3423605}. Our numerical computations in Sect. \ref{sec:I_check} illustrate that no spurious eigenvalues are present in the Lippmann-Schwinger based formulation. However, the proposed filtering process of the finite element approximations described in Sect. \ref{sec:I_check} is in general a much more efficient way to compute resonances.

Computations of the \emph{pseudospectra} provide insight into the behavior of the resolvent $(T^\fem)^{-1}$. In these computations, we use that $\sigma_\epsilon(T^\fem)$ is the set of all $z\in\mx C$ such that 
\begin{equation}
	s_{\scriptsize \hbox{min}}\, T^\fem(z)<\epsilon,
	\label{eq:pseudo_mat}
\end{equation} 
where $s_{\scriptsize \hbox{min}}\, T^\fem(z)$ denotes the smallest singular value of $T^\fem(z)$ \cite[Def. 2.10]{trefethen07}. For the singular value computations we used SLEPc \cite{slepc05}.

\subsection{Numerical detection of spurious solutions}\label{sec:I_check}

In this section we derive a discrete form of \eqref{IntegralTest} that allow us to identify resonances from spurious solutions once we have computed FE solutions $(u^\nu_j,k^\nu_j)$ to \eqref{eq:dtn_discrete} or to \eqref{eq:mat_pml}.
The resulting expression for the filter is a discrete form of the condition $\|\chi_{d}T(k)\chi_{d}u\|<\epsilon$, where $u$ is a FE solution restricted to $\Omega_d$. In the numerical computations we use the minimal computational domain $\Omega_r:=\text{supp}\, (n^2-n_0^2)$.
Then, the Lippmann-Schwinger equation \eqref{eq:LS} can be written in the form
 \begin{equation}\label{eq:discrete_lipp}
	T(k)u=u-K(k)u\,\,\,\hbox{with}\,\,\, K(k)u:=\frac{ik}{2n_0}\int_{\Omega_r}e^{in_0k |x-y|}(n^2(y)-n^2_0)u(y,k)\, dy.
\end{equation}
Let $\{\varphi_j\}$ be a basis for $S^\fem(\Omega_r)$ and let $P^\fem$ be the $L_2$-projection on $S^\fem(\Omega_r)$. Then, we define
\[
u^\fem:=\sum_{j=1}^{N_r} \xi_j\varphi_j,\quad 
P^\fem K(k^\fem)u^\fem:=\sum_{j=1}^{N_r} \eta_j\varphi_j,\,\,\,
M^r_{ij}=\int_{\Omega_r}\varphi_j \varphi_i\,dx,
\]
with $\|u^\fem\|_{L^2(\Omega_r)}=1$, and compute
\begin{equation}
	T^\fem(k^\fem)u^\fem:=u^\fem-P^\fem K(k^\fem)u^\fem=\sum_{j=1}^{N_r} \left(\xi_j-\eta_j \right)\varphi_j(x).
\end{equation} 
The discrete form of $\|\chi_{r}T(k)\chi_{r}u\|<\epsilon$ is then
\begin{equation}
	\|T^\fem(k^\fem)u^\fem\|_{L^2(\Omega_r)}=
	\left|\left|\sum_{j=1}^N 
	\left(\xi_j-\eta_j \right)\varphi_j\right|\right|_{L^2(\Omega_r)}=\sqrt{(\nG\xi-\nG\eta)^T M^r(\nG\xi-\nG\eta)}<\epsilon.
	\label{eq:eps_check}
\end{equation}
Hence, $k^\fem$ belongs, for given $\epsilon>0$, to the $\epsilon$-\emph{psudospectrum} $\sigma_{\epsilon}(T^\fem)$ if the pair $(u^\fem,k^\fem)$ satisfies \eqref{eq:eps_check}.
The integral in \eqref{eq:discrete_lipp} is assembled as the sum of the contributions per element $K_m$. Each subinterval is split in two sub intervals delimited by $x$, then numerical integration is performed by using \emph{Gauss-Legendre} quadratures of the form
$\int_{K_m} f(x)\,dx\approx \sum_{i=1}^{N_q} w_i f(x_i)$, where $w_i$ are the quadrature weights and $x_i$ are the scaled roots of the Legendre polynomials \cite[Ch 4]{solin04}.

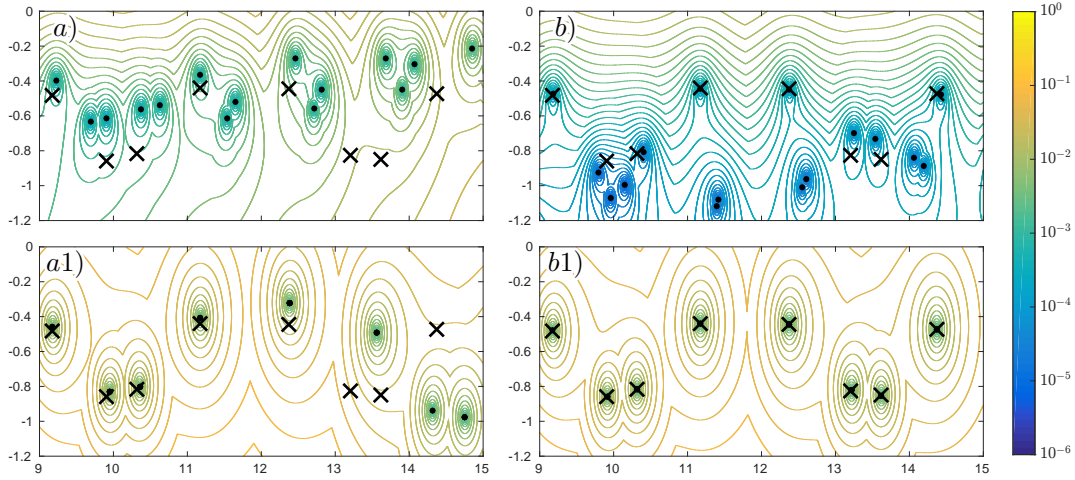
\begin{figure}
	\centering
	\begin{tikzpicture}[thick,scale=0.95, every node/.style={scale=0.95}]
		\draw( 9.80, 1.5) node { \includegraphics[scale=0.425]{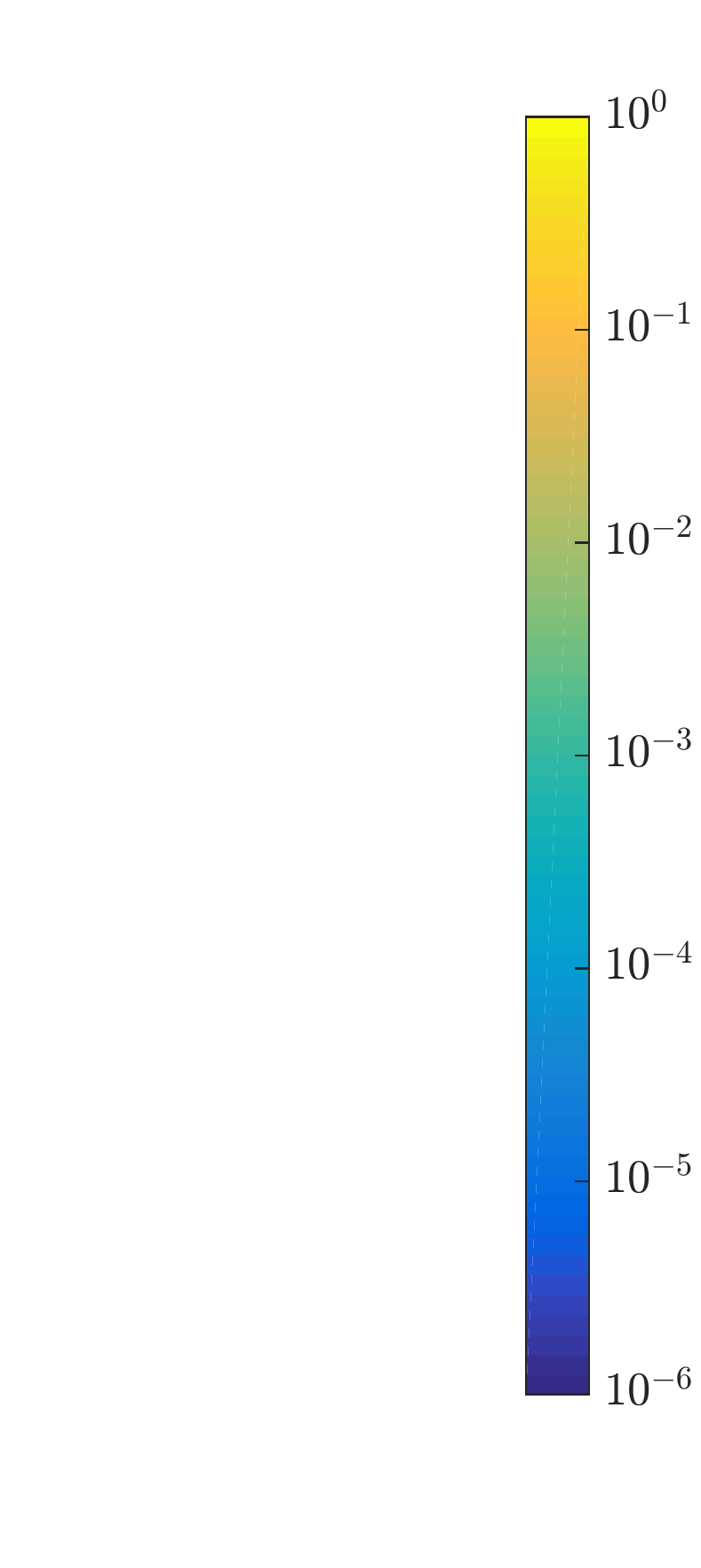} };
		
		\draw( 0.00, 0.0) node { \includegraphics[scale=0.43,clip]{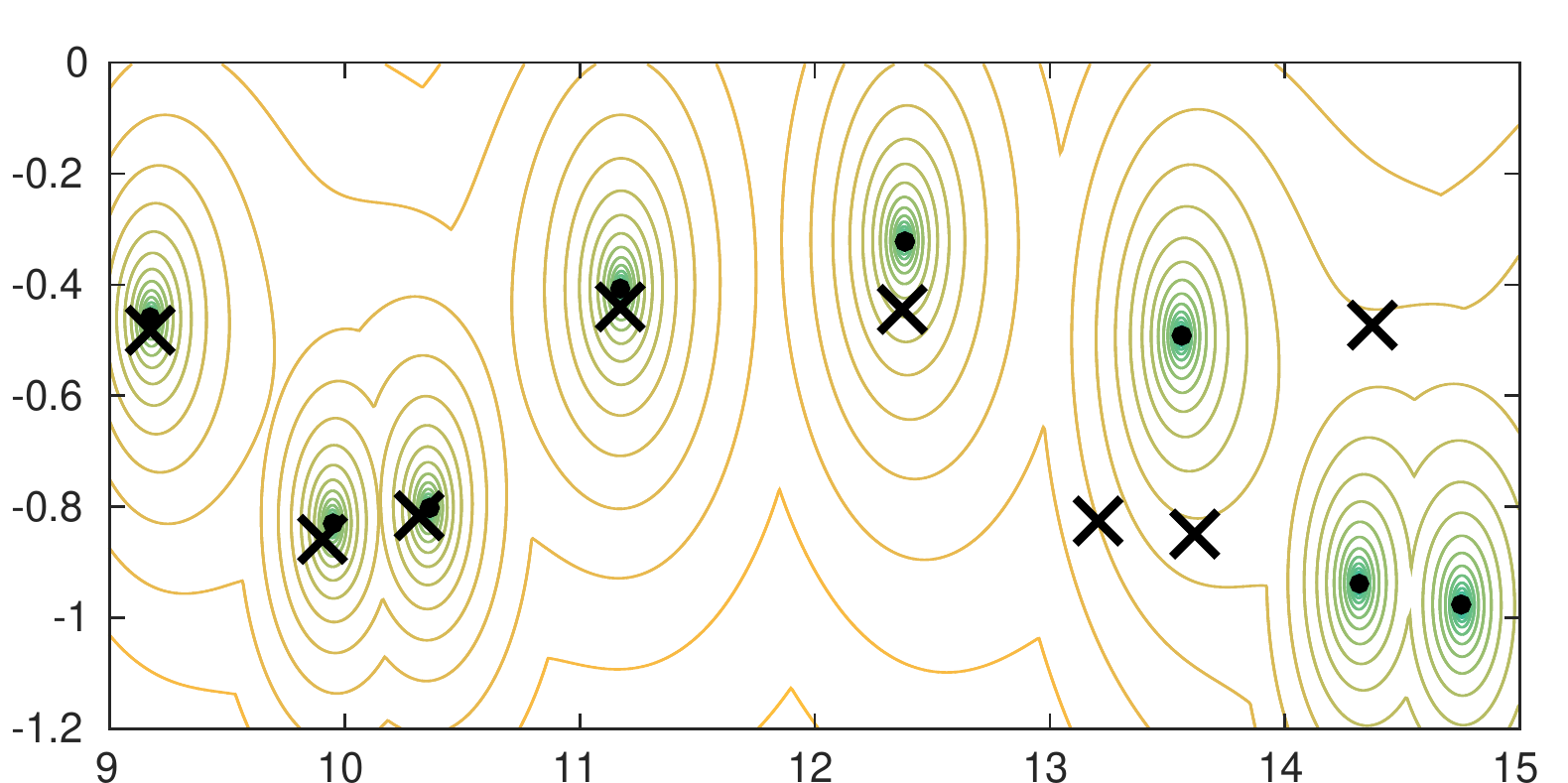} };
		\draw( 7.00, 0.0) node { \includegraphics[scale=0.43,clip]{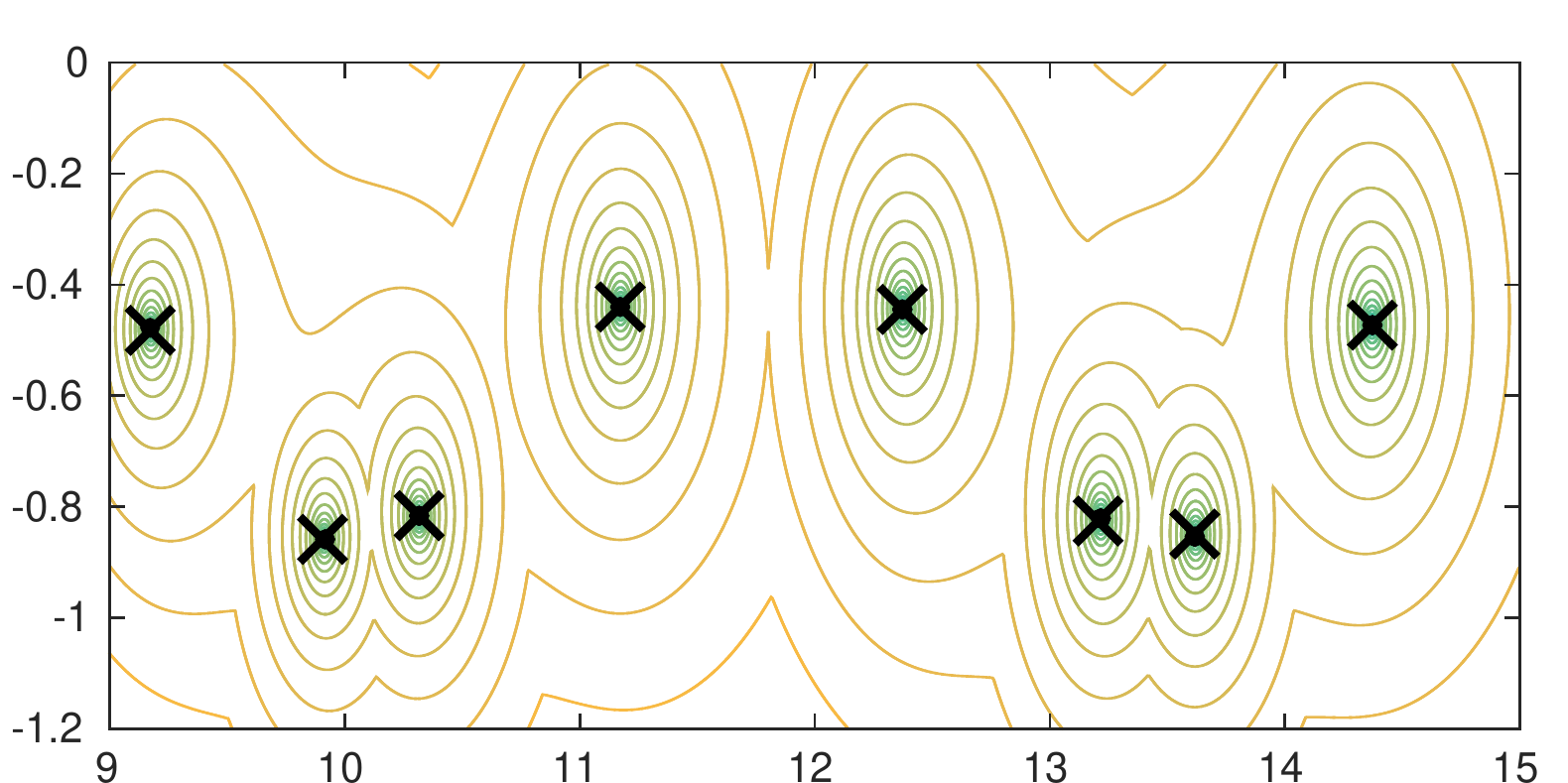} };
		
		\draw( 0.00, 3.3) node { \includegraphics[scale=0.43]{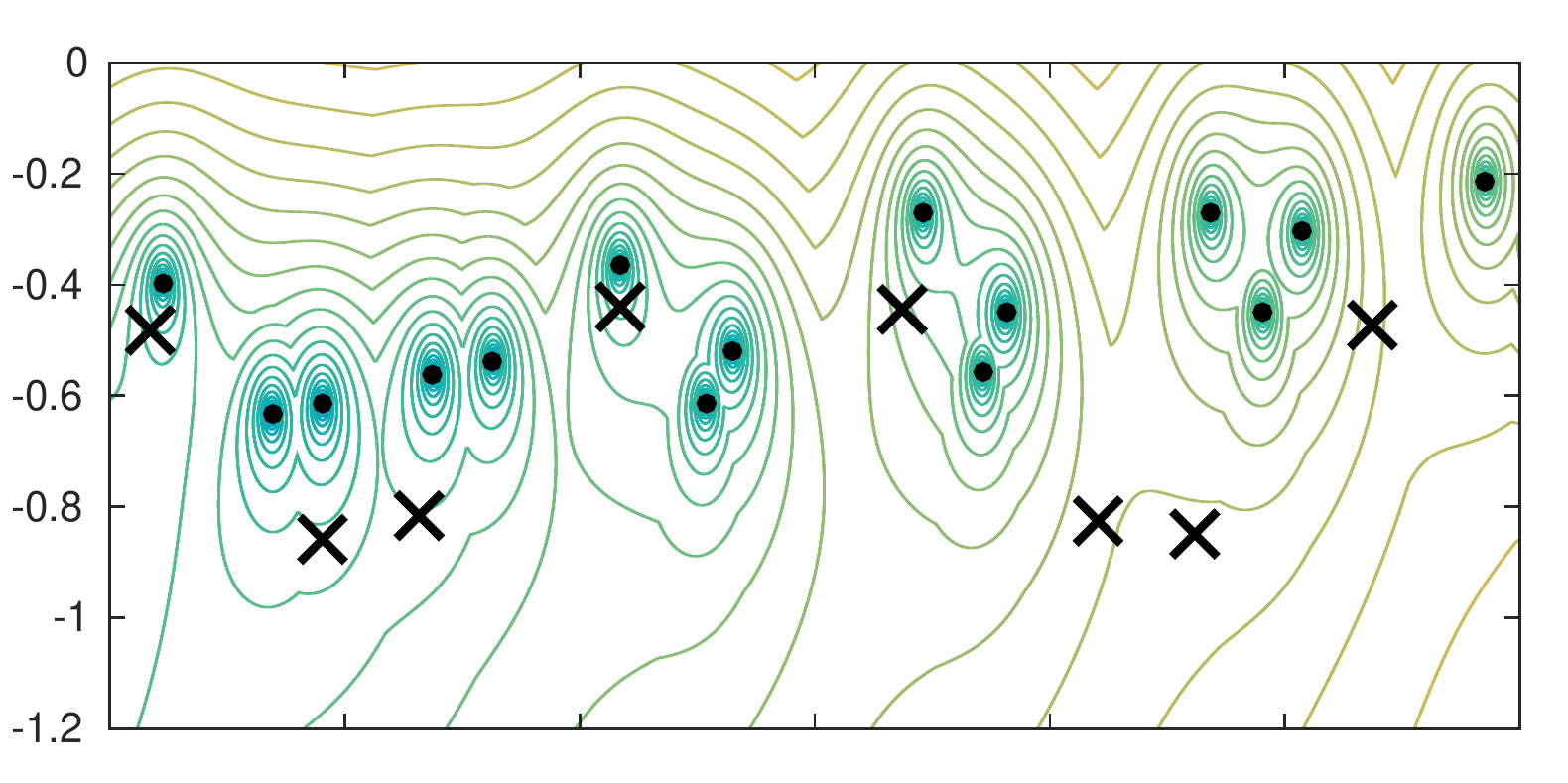} };
		\draw( 7.00, 3.3) node { \includegraphics[scale=0.43]{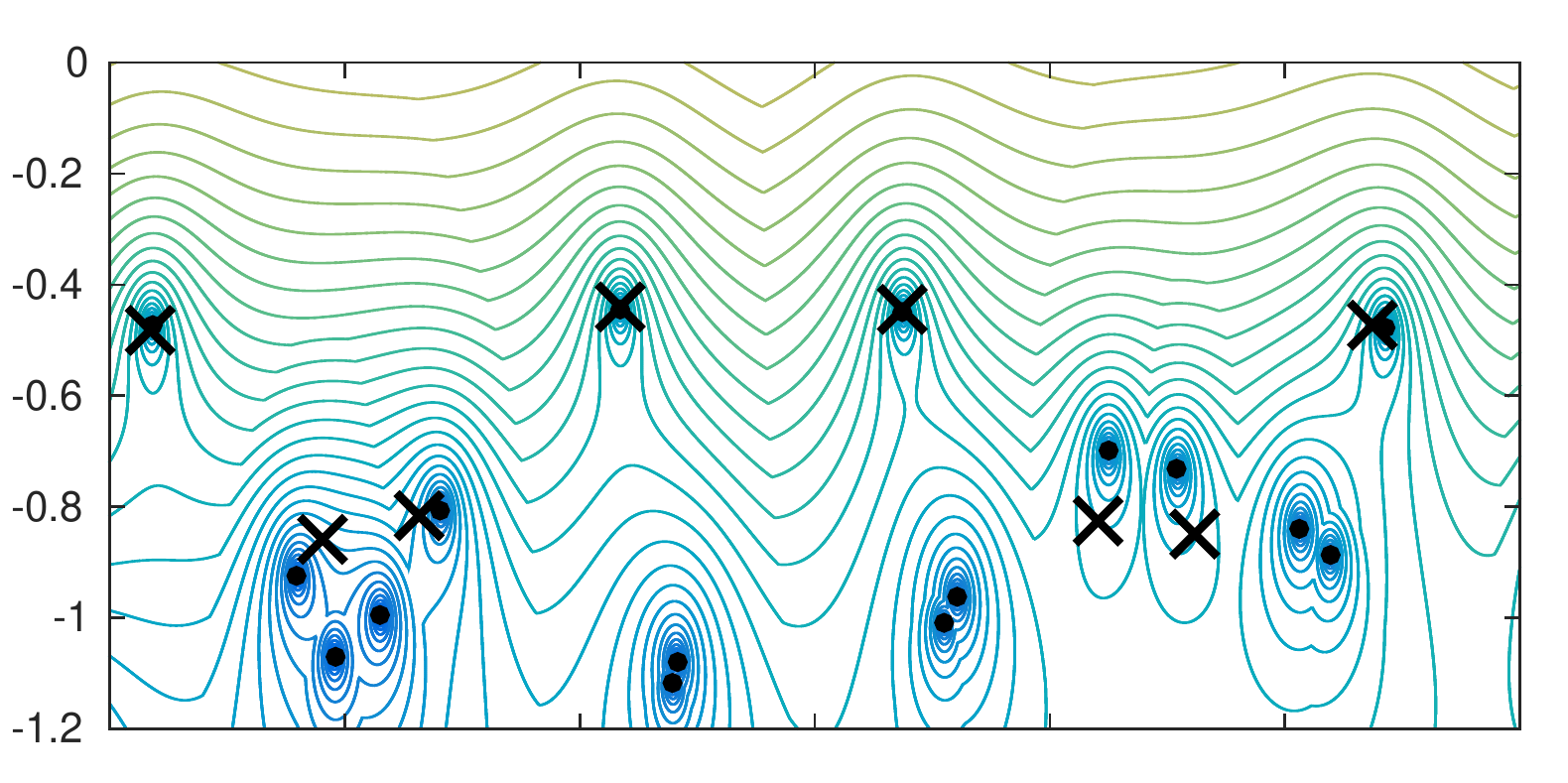} };
		
		\node at (-2.6,4.5) {\contour{white}{\normalsize $a)$}};
		\node at (-2.6,1.2) {\contour{white}{\normalsize $a1)$}};
		
		\node at (4.4,4.5) {\contour{white}{\normalsize $b)$}};
		\node at (4.4,1.2) {\contour{white}{\normalsize $b1)$}};
	\end{tikzpicture}
	\vspace*{-6mm}
	\caption{\emph{Pseudospectrum for the air-filled-cavity problem:  Panels $a)$ and $b)$ are computed using the DtN based formulation and $a1)$ and $b1)$ are computed using the Lippmann-Schwinger formulation. The finite element space used in $a)$, $a1)$ is the same as in Fig. \ref{fig:i_check_DtN_air_filled}, $a)$ and the space used in $b)$, $b1)$ is the same as in Fig. \ref{fig:i_check_DtN_air_filled}, $b)$. For reference, we mark the corresponding eigenvalues \eqref{eq:eigs_tm_kett_1d} with crosses $(\times)$.} }
	\label{fig:psudo_air_filled}
\end{figure}
\subsubsection{Results for the single slab problem}\label{sec:res_single_slab}

In this subsection, the finite element method is used to approximate a selection of eigenpairs $(u_j^\nu,k_j^\nu)$ to the single slab problem described in Sect. \ref{sec:single_slab}. 
The aim is to illustrate the presence of spurious solutions in the different formulations, and not the use of an optimal finite element space. The inclusion of air layers is often unavoidable in applications and we show therefore computations with air-layers in the physical domain.

We set $\Omega_\dtn$ as the physical domain containing the resonator, $\dtn\geq a$, and compare computations on three grids with $a=1$, and $d=1,2,3$. The details of the discretization are $p=2$, $h=0.5$, and the PML parameters are fixed to $\sigma_0=5,\,x_c=d+1,\,\pml=d+3$. 

The results are presented in Fig. \ref{fig:comp_eigs_single_slab}, where we show the exact eigenvalues \eqref{eq:eigs_tm_1d} and approximations computed with the DtN-FEM \eqref{eq:dtn_discrete}, finite PML-FEM \eqref{eq:mat_pml}, and Lippmann-Schwinger \eqref{eq:lipp_collocation}. 
The case with no air-layer is shown in panel ($a$), here the DtN formulation ({\color{blue}\large $*$}) results in no spurious solutions, whereas the PML formulation (\tikz\draw[red,fill=red] (0,0) circle (.5ex);) results in several spurious solutions. The inclusion of air-layers is shown in panels ($b$) and ($c$), where we observe an increased number of spurious eigenvalues in the PML formulation as well as in the DtN based formulation. However, no spurious solutions were computed with the Lippmann-Schwinger formulation (\tikz\draw[black] (0,0) circle (.5ex);) and air-layers are of no concern since the integration is only over $\text{supp}\, (n^2-n_0^2)$. Moreover, the accuracy of the computed set of eigenvalues was superior to DtN-FEM and PML-FEM. This supports the idea introduced in Sect. \ref{sec:resolvent} that a filter based on the Lippmann-Schwinger operator can be used to identify spurious solutions in the PML formulation and in the DtN formulation.

\begin{figure}[!h]
	\begin{tikzpicture}[thick,scale=0.9, every node/.style={scale=0.9}]
		
		\draw( 0.00, 7.2)  node { \includegraphics[scale=0.45]{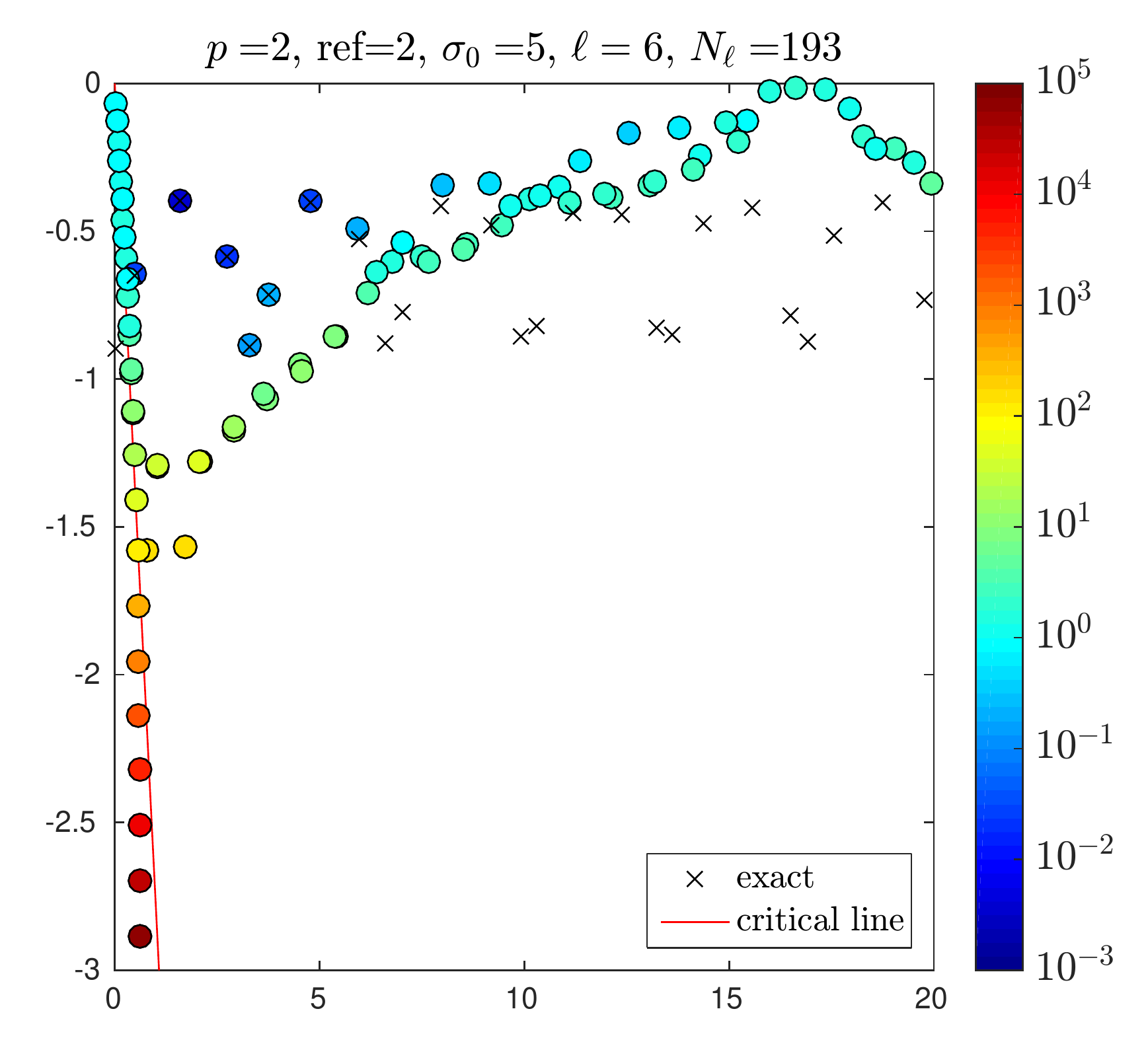} };
		\draw( 8.00, 7.2)  node { \includegraphics[scale=0.45]{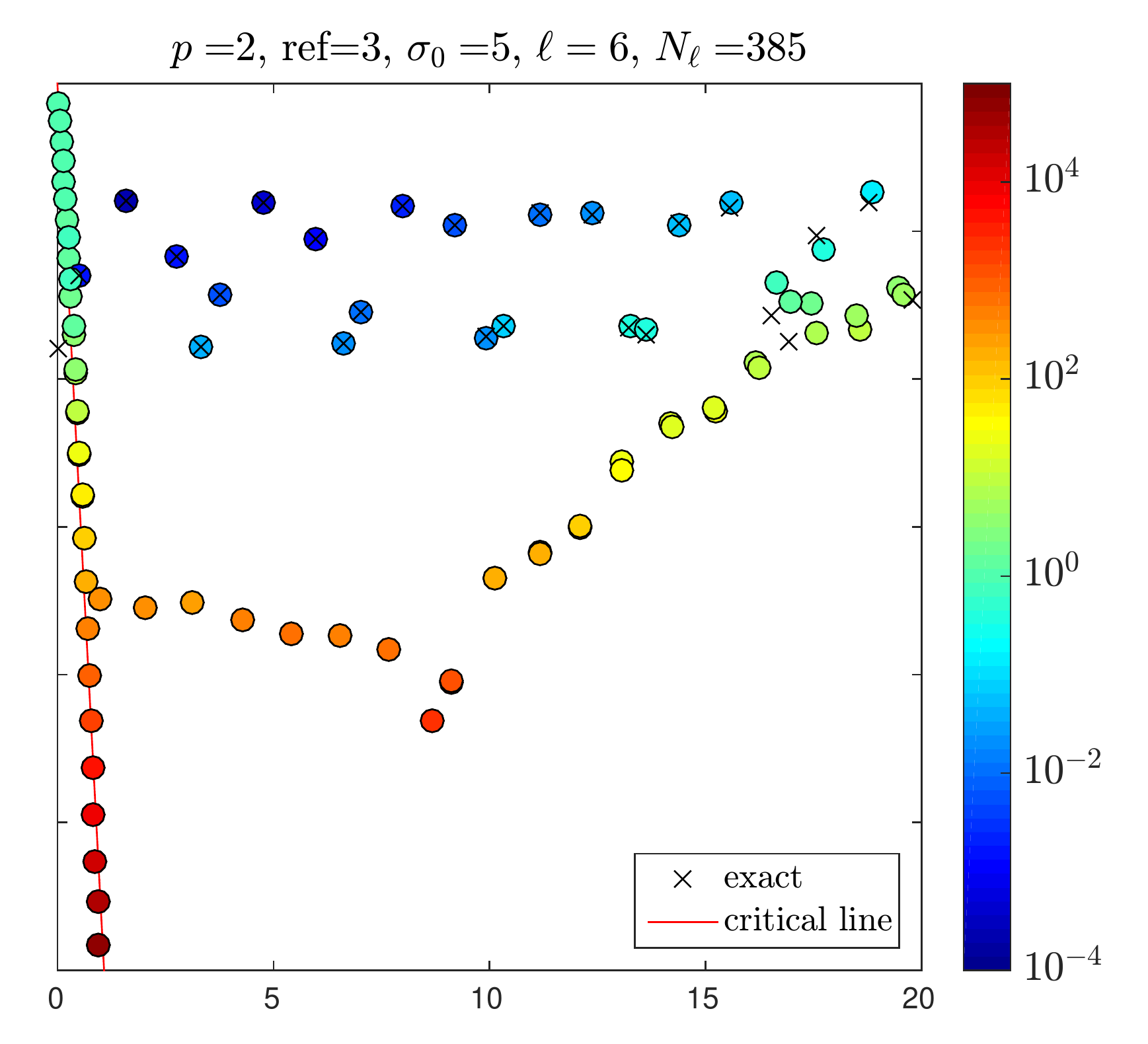} };
		
		\draw( 0.00, 0.0)  node { \includegraphics[scale=0.45]{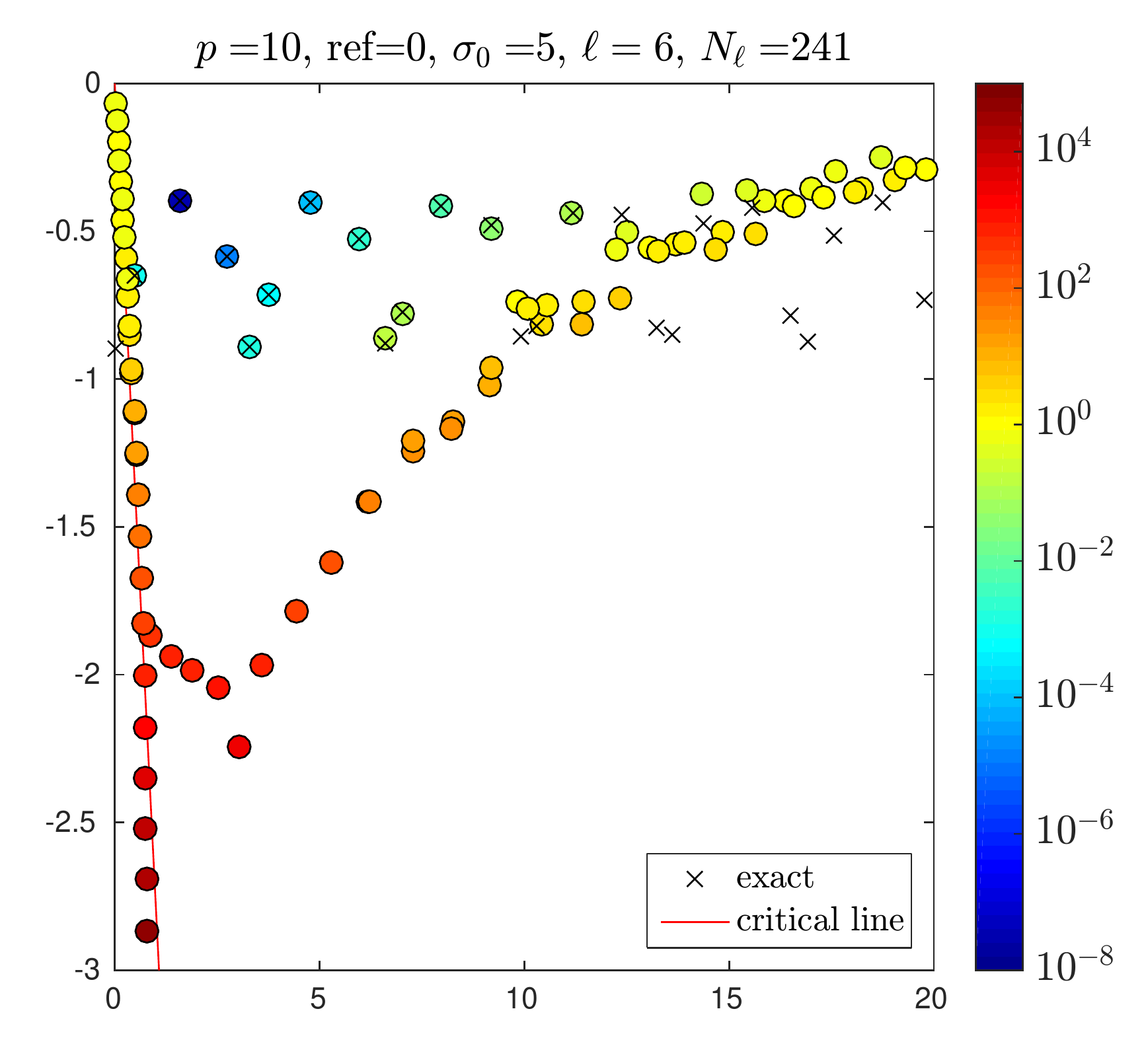} };
		\draw( 8.00, 0.0)  node { \includegraphics[scale=0.45]{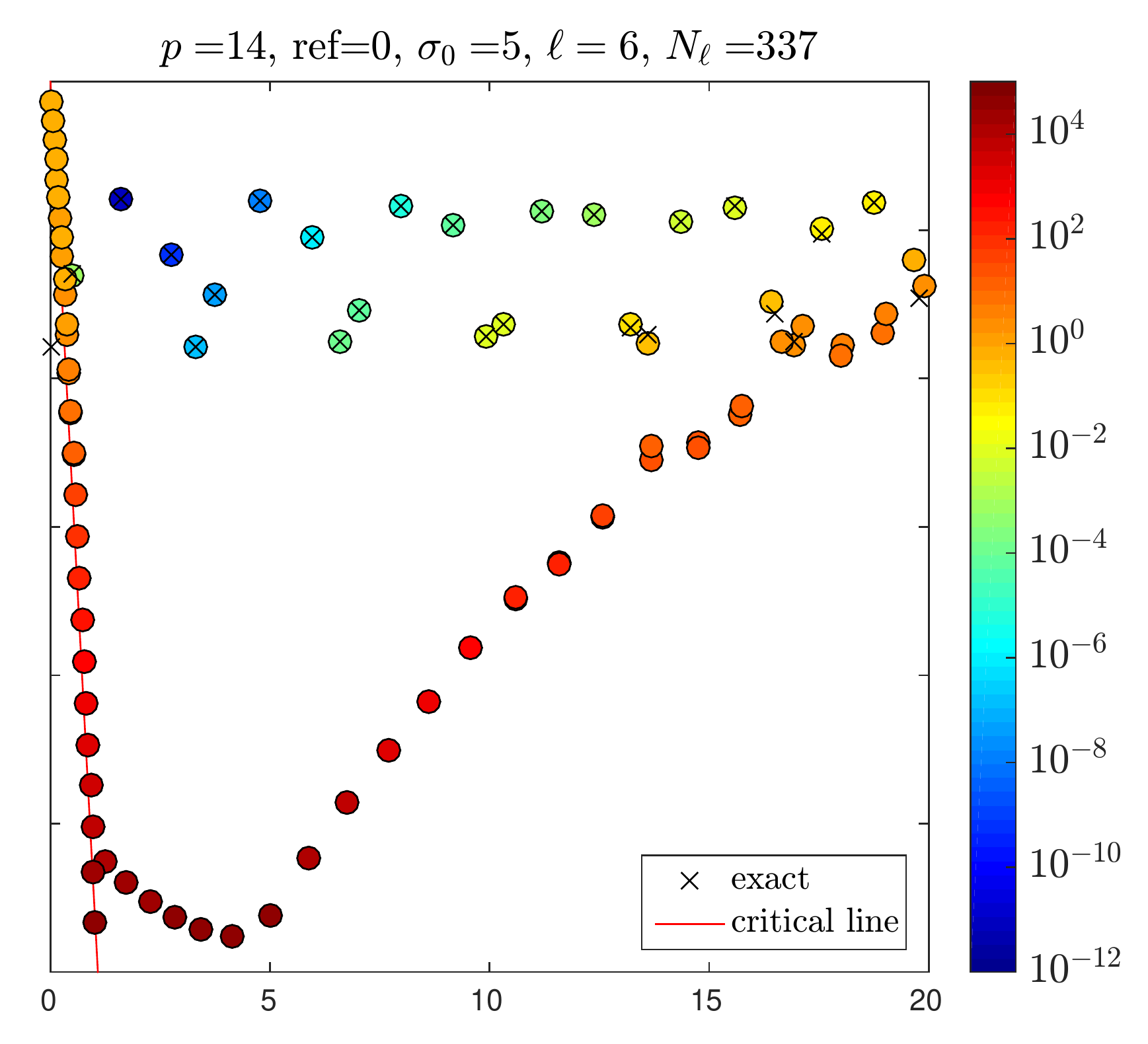} };
		
		\node[text width=3cm] at (-1.60, 3.25) {$c)$};
		\node[text width=3cm] at ( 6.00, 3.25) {$d)$};
		
		\node[text width=3cm] at (-1.60, 10.45) {$a)$};
		\node[text width=3cm] at ( 6.00, 10.45) {$b)$};
		
	\end{tikzpicture}
	\vspace*{-5mm}
	\caption{\emph{Computed Eigenvalues $k_j^\fem$ of \eqref{eq:mat_pml}, corresponding to the air-filled-cavity problem described in Sect. \ref{sec:air_filled}. In colors we give $\epsilon_j$ computed with \eqref{eq:eps_check} for each $k_j^\fem$.} }
	\label{fig:i_check_PML_air_filled}
\end{figure}
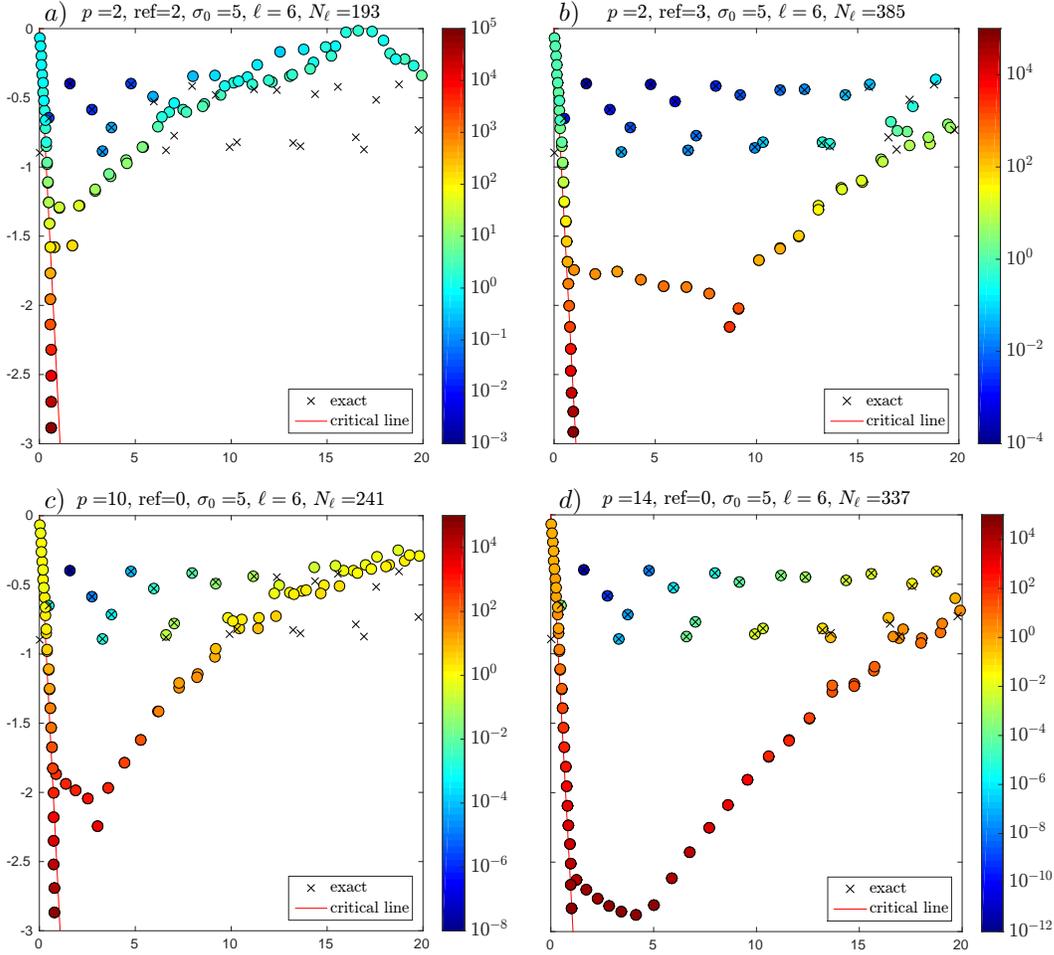

In Sect. \ref{sec:dtn_pml}, we showed that the coupled problem \eqref{eq:trun_PML_split} decouples when $k\in\mathcal{D}_c:=\{k\in\C\,:\,1\pm e^{2ik\beta}=0 \}$ but there are no eigenvalues in the set $\mathcal{D}_c$. However it is plausible that solutions exist close to $\mathcal{D}_c$, and the results presented in Fig. \ref{fig:comp_eigs_single_slab}, Fig. \ref{fig:i_check_PML_air_filled}, Fig. \ref{fig:PML_air_filled_sigma}, and Fig. \ref{fig:psudo_bump} show numerical eigenvalues located close to $\mathcal{D}_c$. For the single slab problem we have used \eqref{eq:eig_tm_pml_1d} to verified that they are indeed approximations of eigenvalues of the finite PML problem. This expression is an exact relationship for all the eigenvalues and by using a Newton root solver with initial guesses in $\mathcal{D}_c$ we find solutions close to the FE approximations. Note that these eigenvalues are not approximations of resonances.


\subsubsection{Results for the air-filled-cavity problem}\label{sec:res_air_filled}

In this subsection, the finite element method is used to approximate a selection of eigenpairs $(u_j^\nu,k_j^\nu)$ to the air-filled-cavity problem with the DtN-map and with the PML formulations. The problem is described in Sect. \ref{sec:air_filled} and reference solutions are listed in Tab. \ref{tab:air_filled_reference}. Then, we use \eqref{eq:eps_check} to determine the smallest $\epsilon>0$ such that $k_j^\nu\in\sigma_{\epsilon}(T^\fem)$.
\begin{figure}[h!]
	\centering
	\begin{tikzpicture}
		\draw( 0.00, 00.0) node { \includegraphics[scale=0.65]{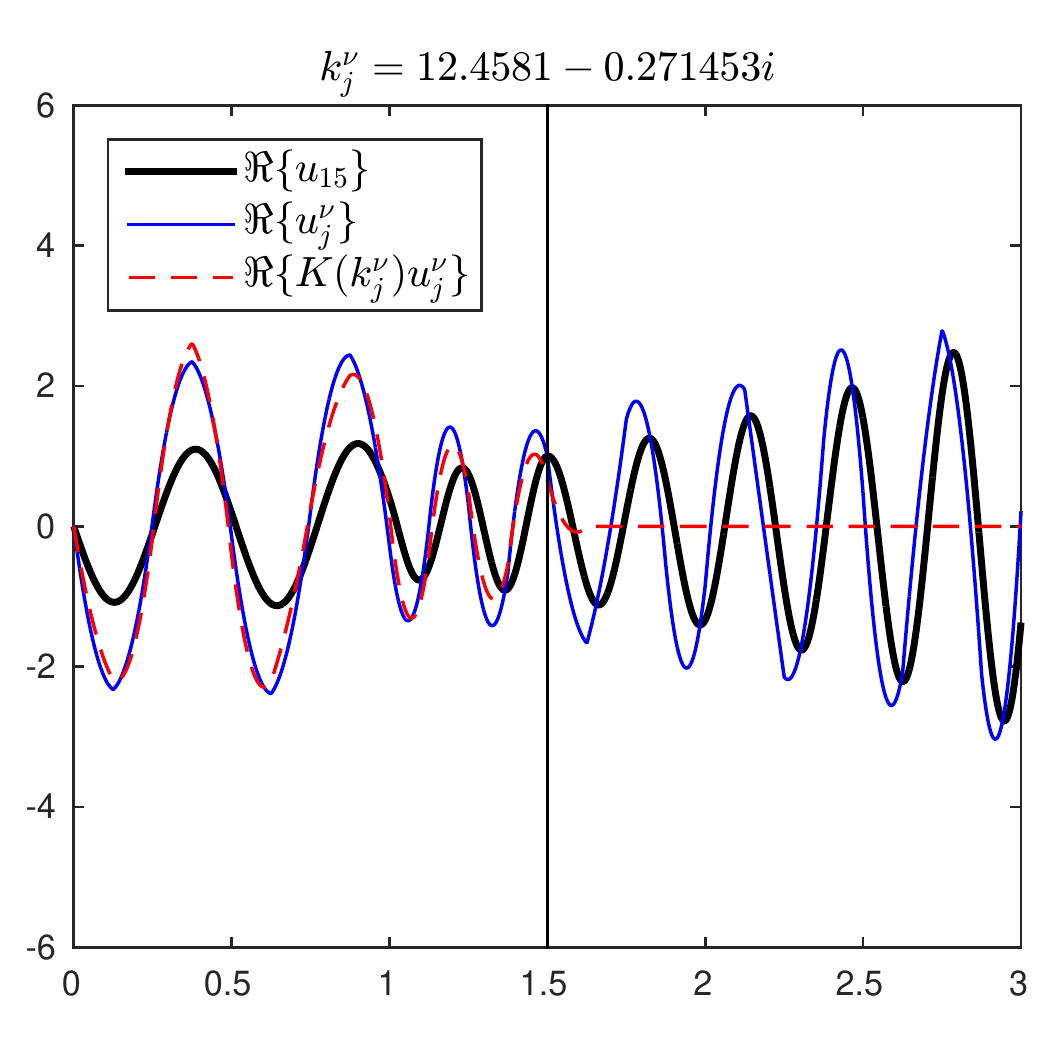} };
		\draw( 7.50, 00.0) node { \includegraphics[scale=0.65]{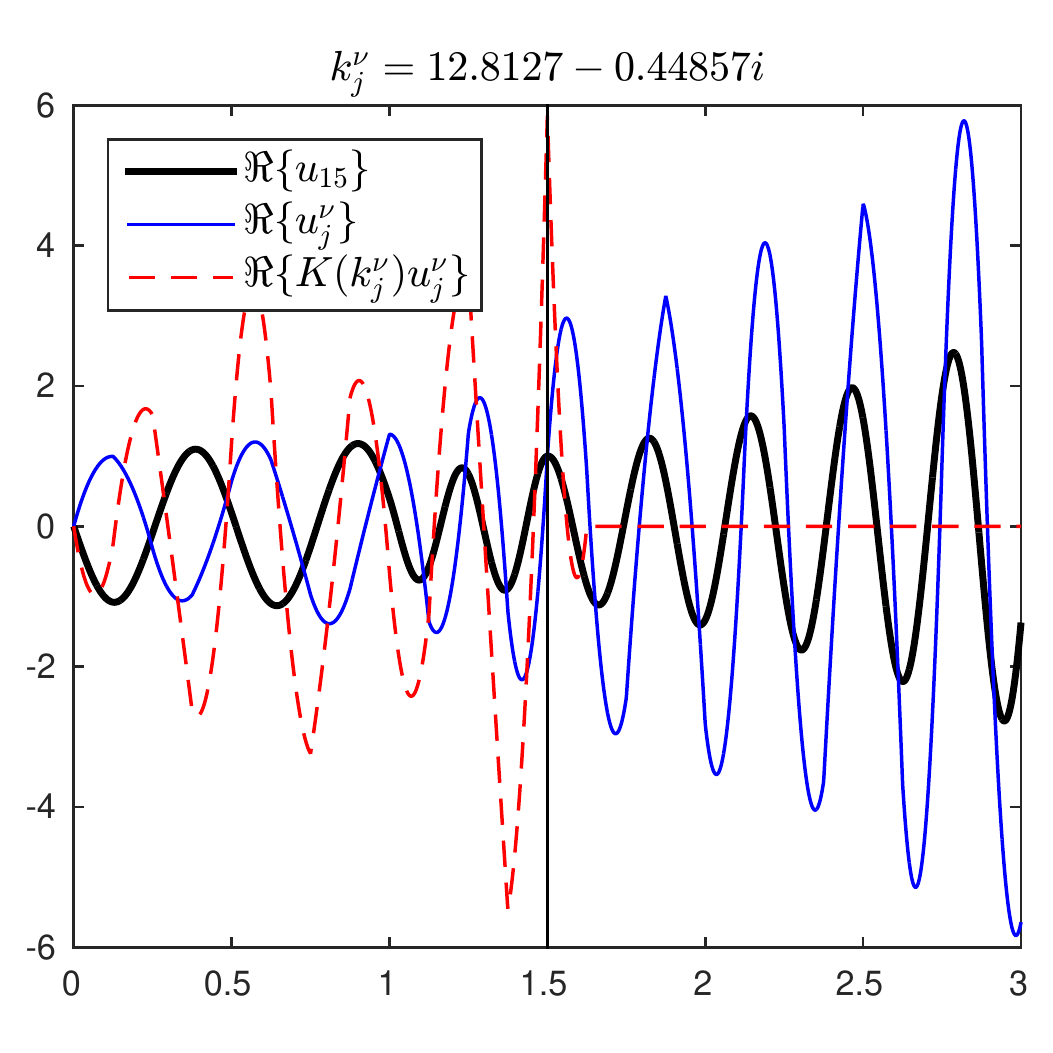} };
		
		\node[text width=3cm] at (-1.30, 3.05) {$a)$};
		\node[text width=3cm] at ( 6.30, 3.05) {$b)$};
	\end{tikzpicture}
	\vspace*{-5mm}
	\caption{\emph{Functions corresponding to eigenpairs enclosed by the black oval 
	in Fig. \ref{fig:i_check_DtN_air_filled}-(a).	We compare the real parts of the exact resonance function $u_{15}$ (thick line), with computed eigenfunctions $u^\fem_j$ (thin line), and $K(k^\fem_j) u^\fem_j$ (dashed line). For the two chosen pairs, $k_j^\fem\in \sigma_{\epsilon}(T^\fem)$
	with $\epsilon=0.325$ in $a$), 
	while $\epsilon=1.99$ in $b$) . In agreement, we see that the pointwise distance $|\Re\{u^\fem_j(x)-K(k^\fem_j) u^\fem_j(x)\}|$ is smaller in $a)$ compared to in $b)$.  }}
	\label{fig:dtn_comp_sp}
\end{figure}
\textit{The DtN based formulation:}
The estimate \eqref{eq:OsbornEst} shows that the gap will not decrease when the resolvent norm of at least one of the block operators \eqref{op:L}, \eqref{GeqL} dominates the approximability of the finite element space. The resolvent norm of \eqref{op:L}, as discussed in Sect. \ref{sec:resolvent}, will grow exponentially when $\im k\rightarrow -\infty$.
In Fig. \ref{fig:i_check_DtN_air_filled} we show eigenvalues $k_j^\fem$ for various discretizations. An eigenpair is numerically close to a resonance pair if 
$k_j^\nu\in\sigma_{\epsilon}(T^\fem)$ for a small $\epsilon$. 
Note that the test \eqref{eq:eps_check} uses not only the eigenvalues but also the corresponding  eigenvectors. As expected $\epsilon$ is in general larger for $k$ away from the real axis and the minimum $\epsilon$ is then found for low values on $\re k_j$. The reason is that the quality of the approximated pairs $(u_j^\fem,k_j^\fem)$ deteriorate with greater oscillatory behavior \cite{Sauter10}. We remark
that the size of the linearized problem \eqref{eq:mat_dtn} is $2N_\dtn$. 

Panels \ref{fig:i_check_DtN_air_filled} a), b) show that for a low polynomial degree ($p=2$) the lowest eigenvalue $k_j^\fem$ is in $\sigma_{\epsilon}(T^\fem)$ for $\epsilon\approx 10^{-2}$. In Fig. \ref{fig:i_check_DtN_air_filled} c) we use $p=10$ giving $\epsilon\approx 10^{-9}$ for the lowest eigenvalue. Fig. \ref{fig:i_check_DtN_air_filled} d) shows the computed eigenvalues for $p=14$. No spurious eigenvalues are computed in the selected region and the lowest eigenvalue is in $\sigma_{\epsilon}(T^\fem)$ for $\epsilon\approx 10^{-11}$. 

\textit{Reliability of the filtering process:} \\
In Fig. \ref{fig:psudo_air_filled}, we show eigenvalues reference ($\times$), for the DtN formulation in panels $a)$ and $b)$ that correspond to the discretizations $a)$ and $b)$ from Fig. \ref{fig:i_check_DtN_air_filled}. In panels $a1)$ and $b1)$ we show eigenvalues for the Lippmann-Schwinger formulation. Discretizations in $a)$ and $a1)$ have the same span of functions $\varphi_j$ covering $\Omega_r$, and similarly for discretizations in $b)$ and $b1)$. As in Sect. \ref{sec:res_single_slab}, we observe that while the Lippmann-Schwinger formulation gives equal number of eigenvalues as exact resonances, the DtN formulation gives too many eigenvalues.

To get further insight on the reliability of the filtering test 
we also include (in contours) the pseudospectrum for the DtN formulation and Lippmann-Schwinger formulation in Fig. \ref{fig:psudo_air_filled}. It is evident that the resolvent norm of the DtN formulation grows with $-\Im k$, while for the Lippmann-Schwinger formulation we observe a well-behaved resolvent norm away from the spectrum.

In order to illustrate the role played by the eigenfunctions in the filtering scheme, we refer to
panel $a)$ of Fig. \ref{fig:i_check_DtN_air_filled}. Two computed eigenvalues that are located near to $k_{15}$ ($\times$) are enclosed by a black ellipse. The corresponding eigenfunctions are visualized in Fig. \ref{fig:dtn_comp_sp} where $a$) corresponds to the case $\epsilon=0.325$ and $b$)
to the case $\epsilon=1.99$. 
Clearly the function $u_j^\fem$ in $a$) is a better approximation to $u_{15}$ than the one plotted in $b$). 

Moreover, the function $K(k^\fem_j)u^\fem_j$ in $a$) as defined in \eqref{eq:discrete_lipp} follows $u^\fem_j$ closely, whereas in $b$) there is no correspondence explaining why $\epsilon_j$ is larger in $b$).
\begin{figure}[h!]
	\centering
	\begin{tikzpicture}[thick,scale=0.95, every node/.style={scale=0.95}]
		
		\draw(10.70, 5.42) node { \includegraphics[scale=0.82]{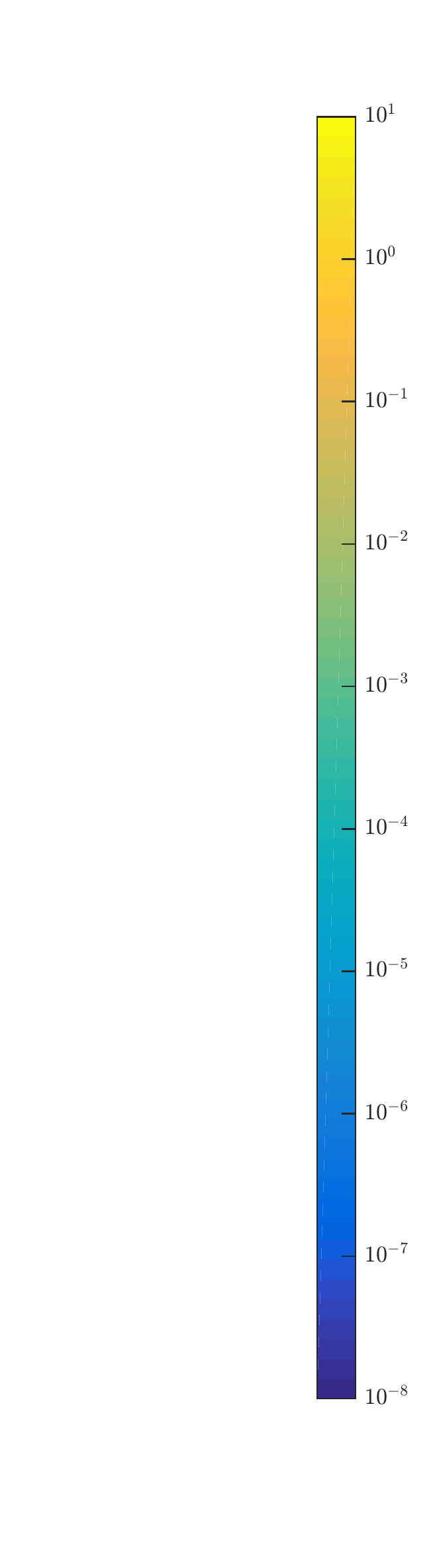} };
		
		\draw( 0.00, 11.4) node { \includegraphics[scale=0.43]{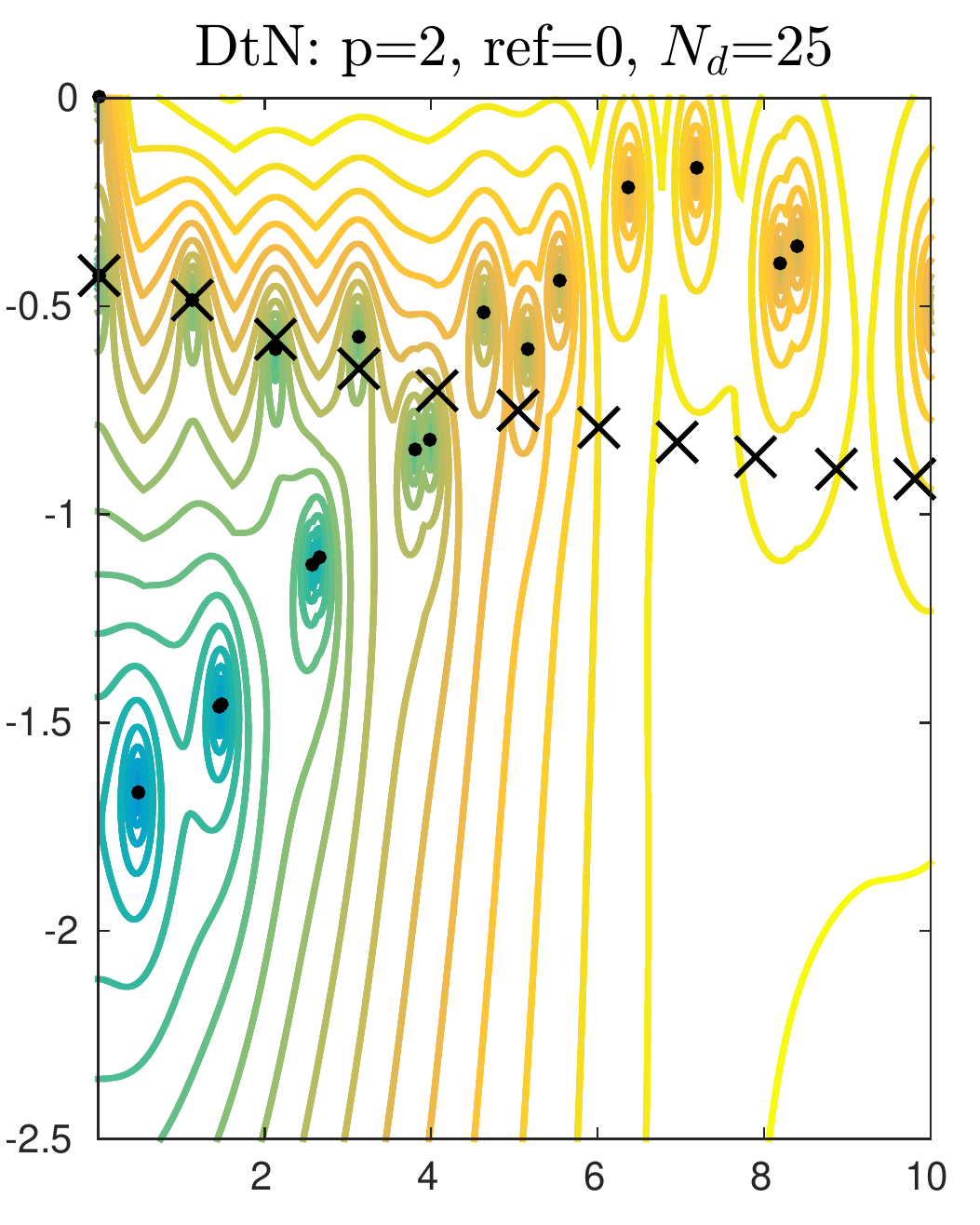} };
		\draw( 4.70, 11.4) node { \includegraphics[scale=0.43]{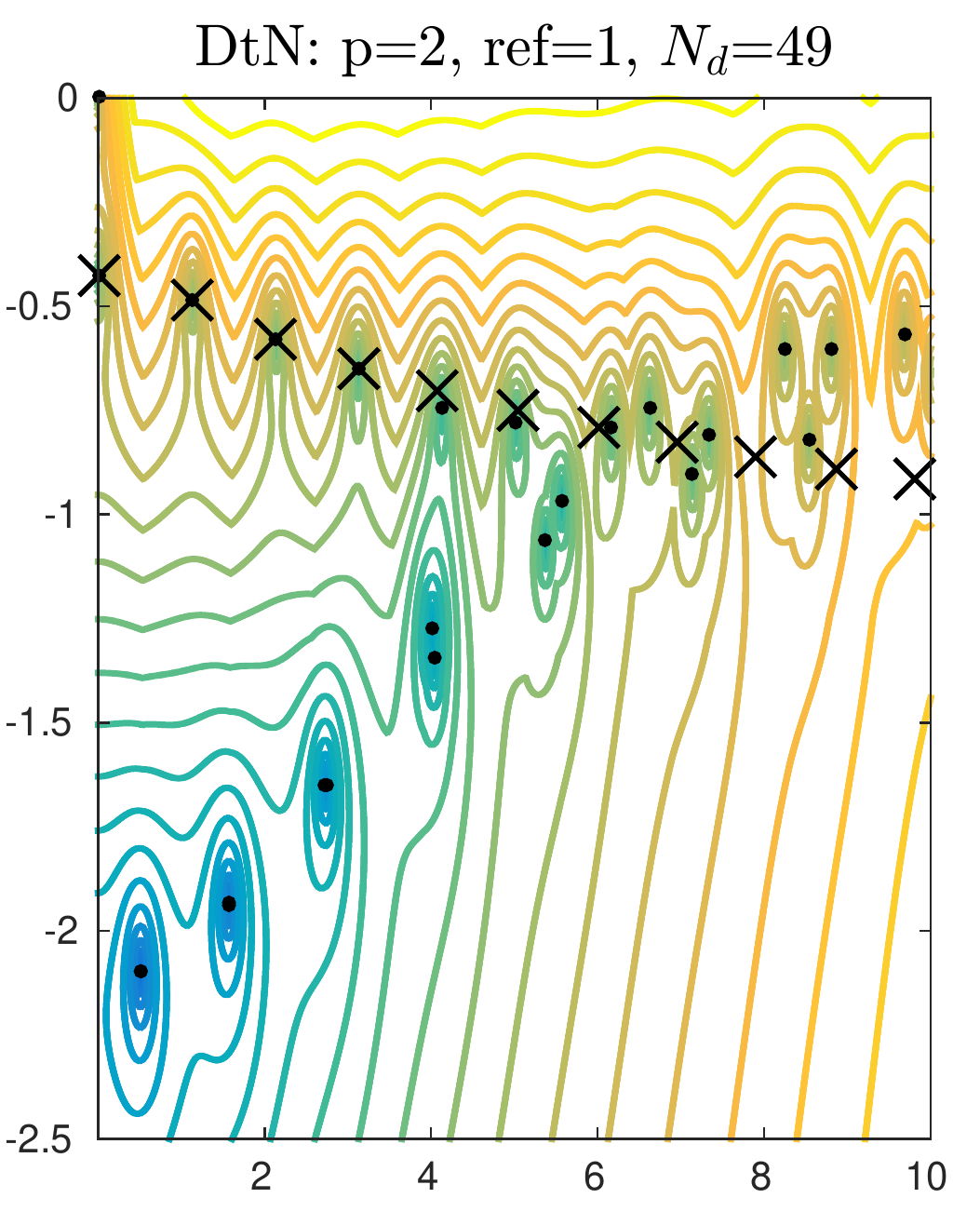} };
		\draw( 9.40, 11.4) node { \includegraphics[scale=0.43]{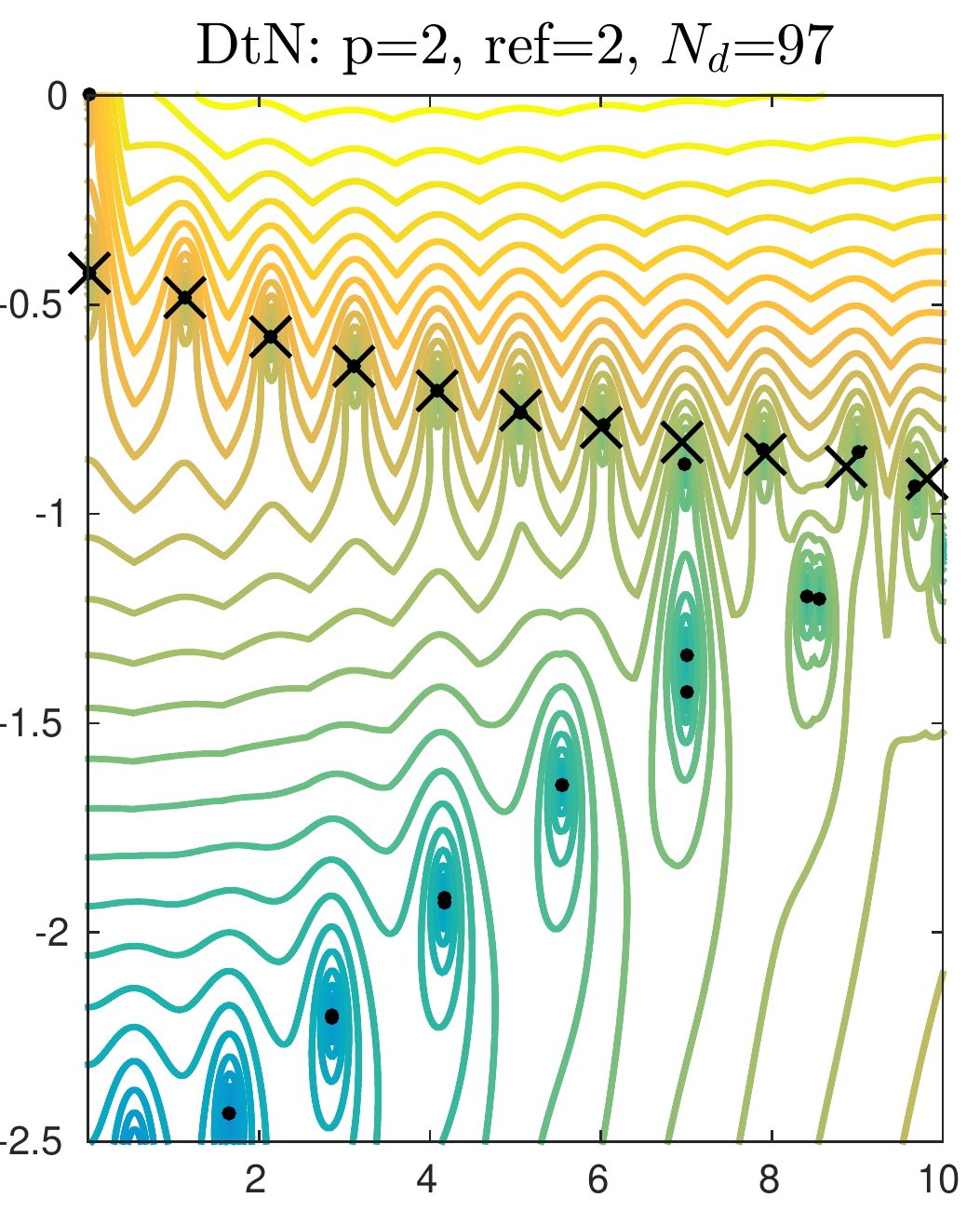} };
		
		\draw( 0.00,  5.7) node { \includegraphics[scale=0.43]{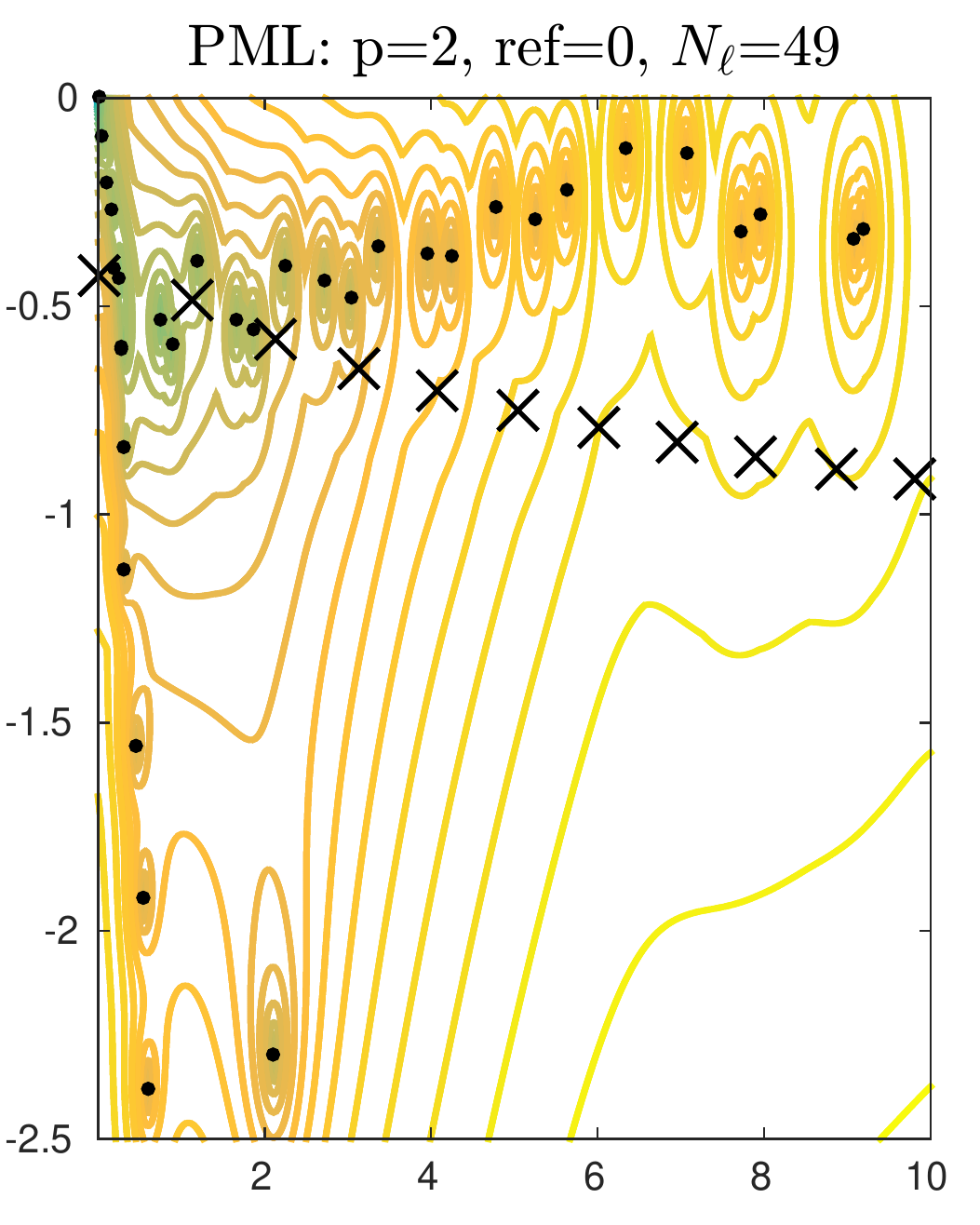} };
		\draw( 4.70,  5.7) node { \includegraphics[scale=0.43]{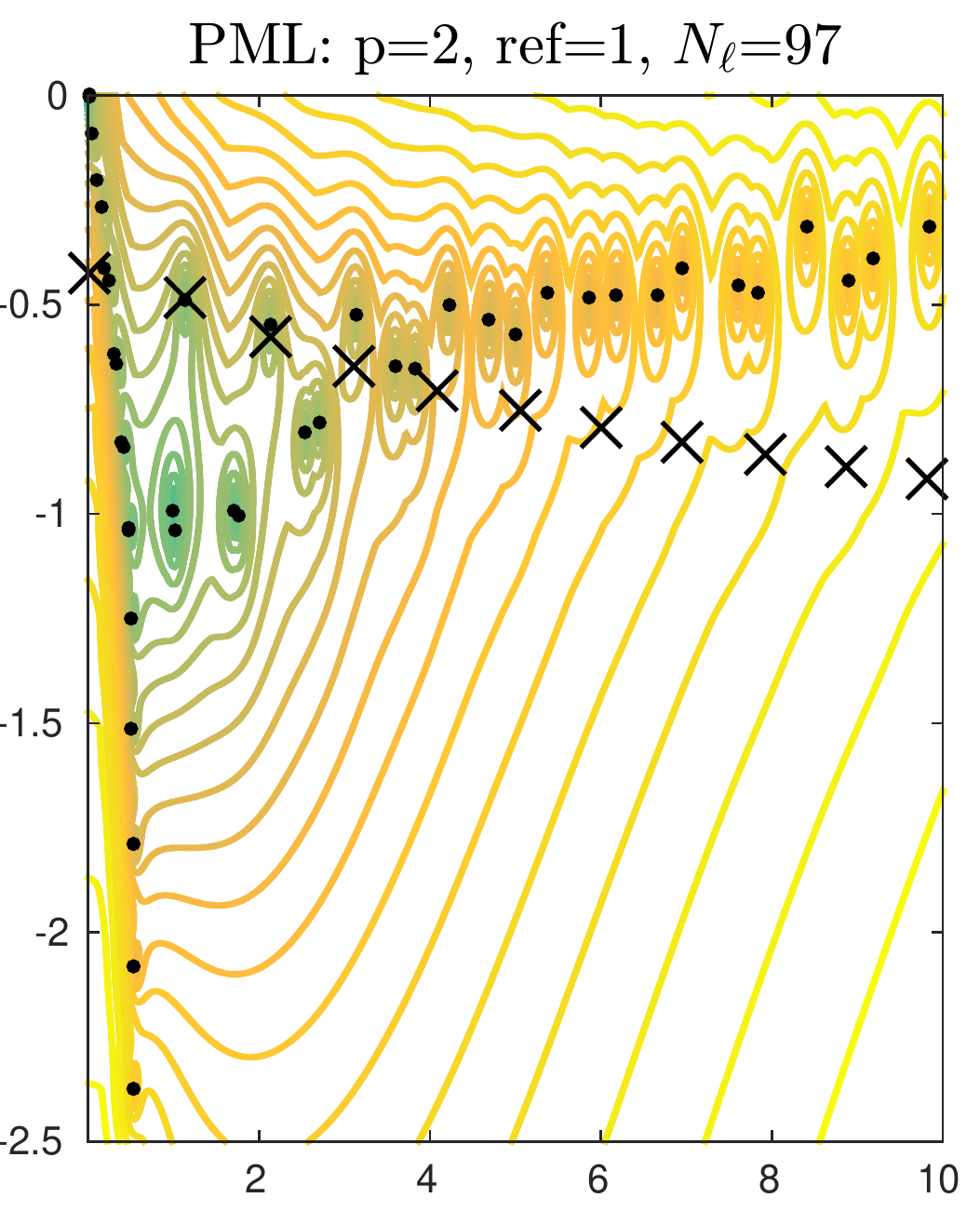} };
		\draw( 9.40,  5.7) node { \includegraphics[scale=0.43]{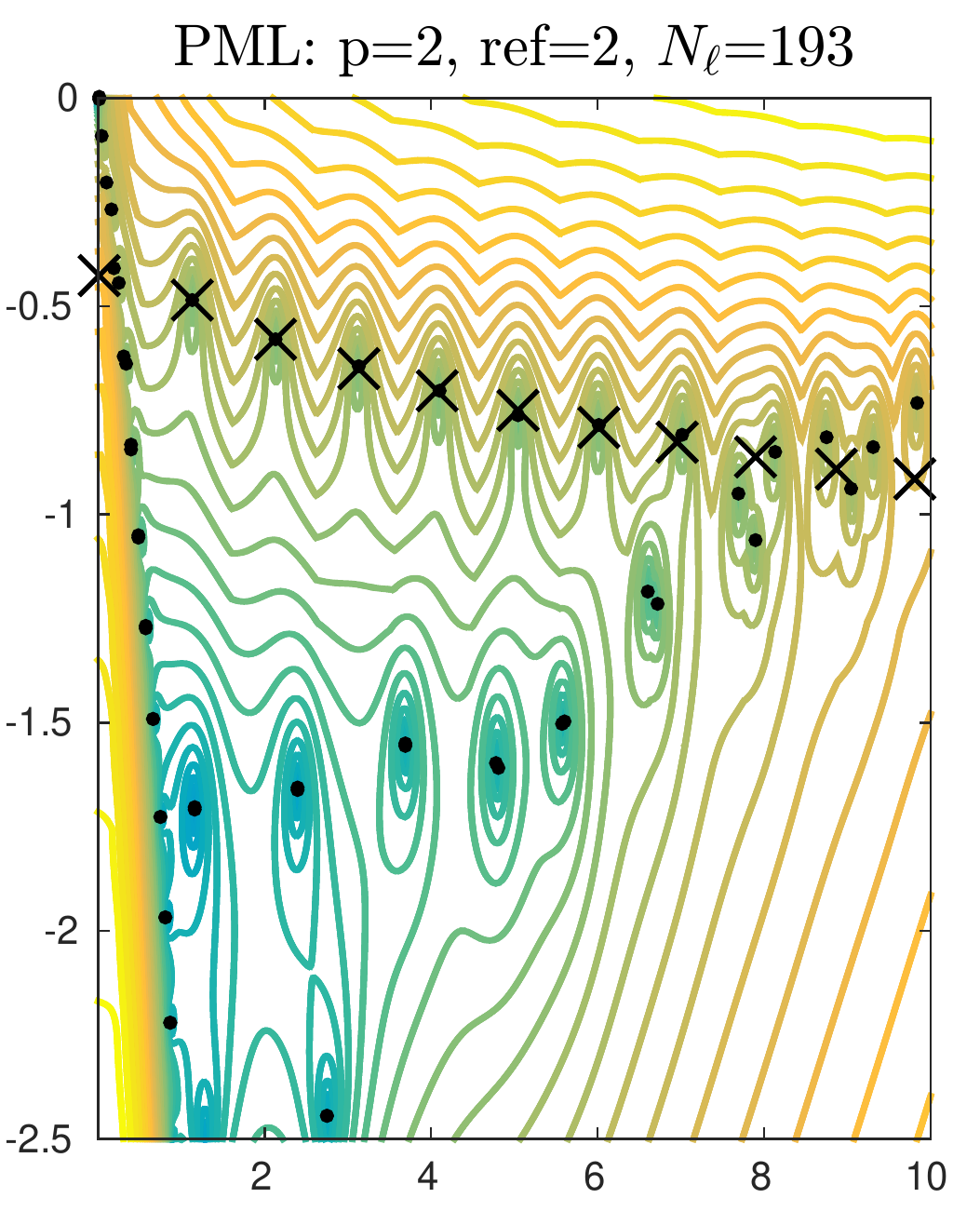} };
		
		\draw( 0.00,  0.0) node { \includegraphics[scale=0.43]{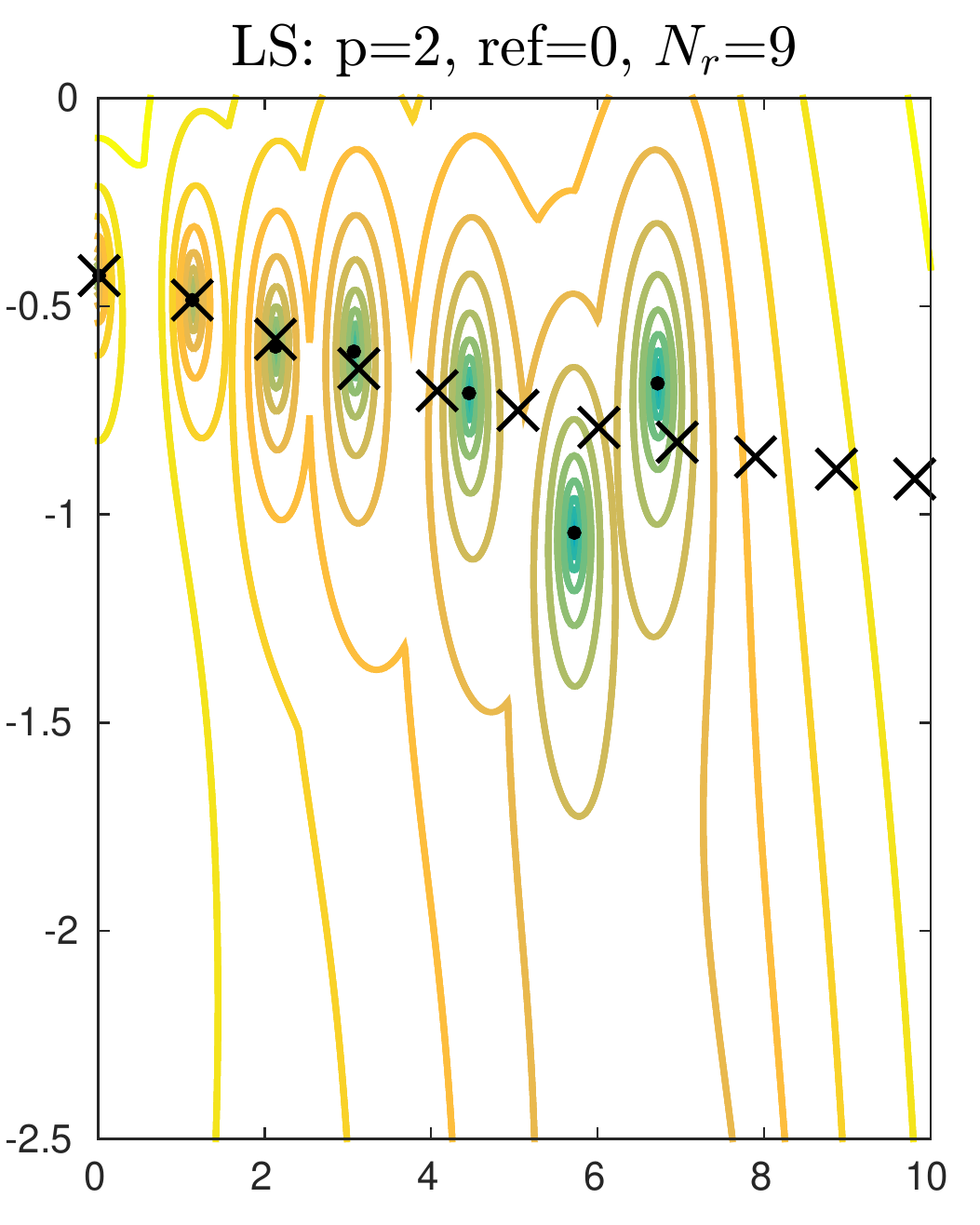} };
		\draw( 4.70,  0.0) node { \includegraphics[scale=0.43]{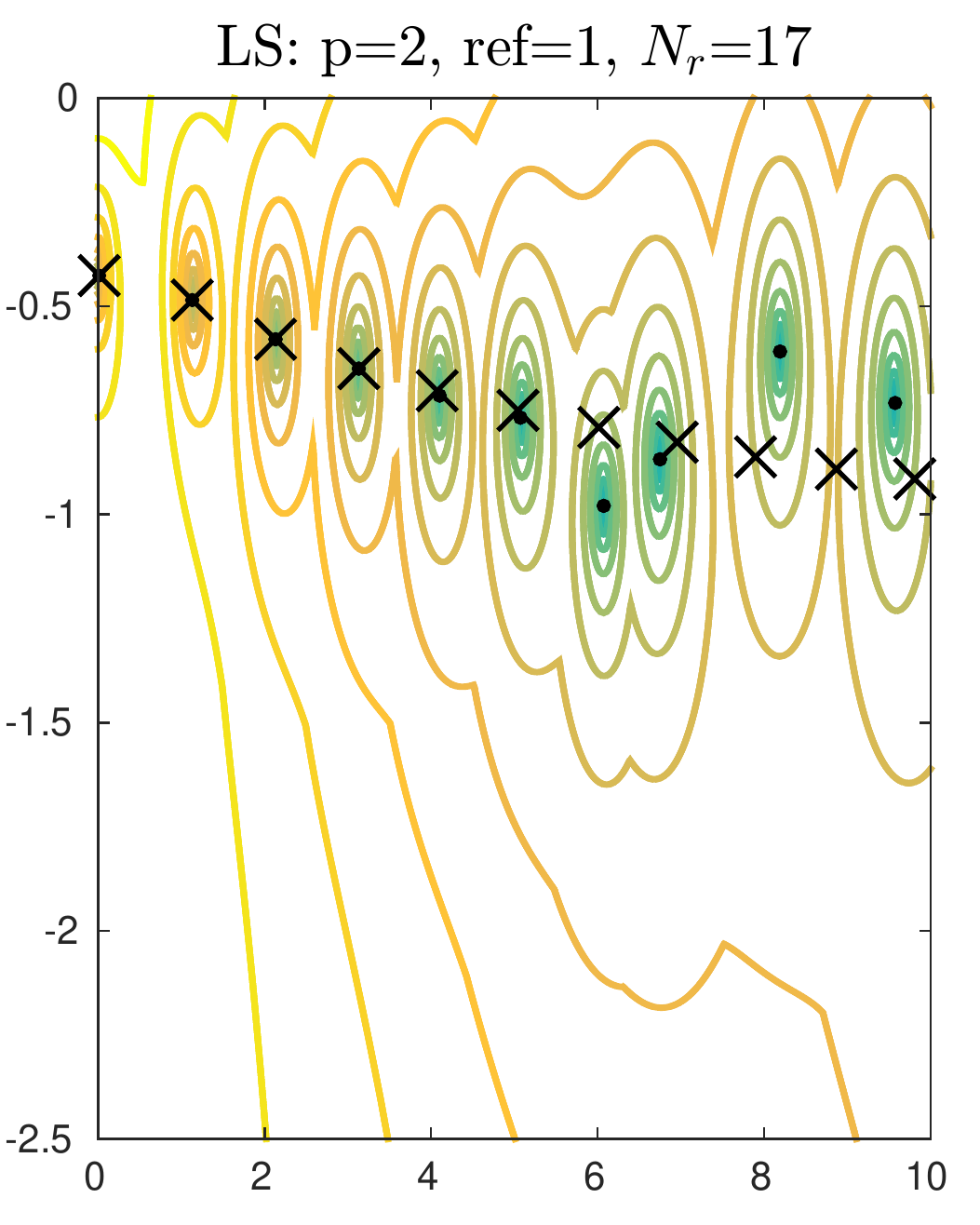} };
		\draw( 9.40,  0.0) node { \includegraphics[scale=0.43]{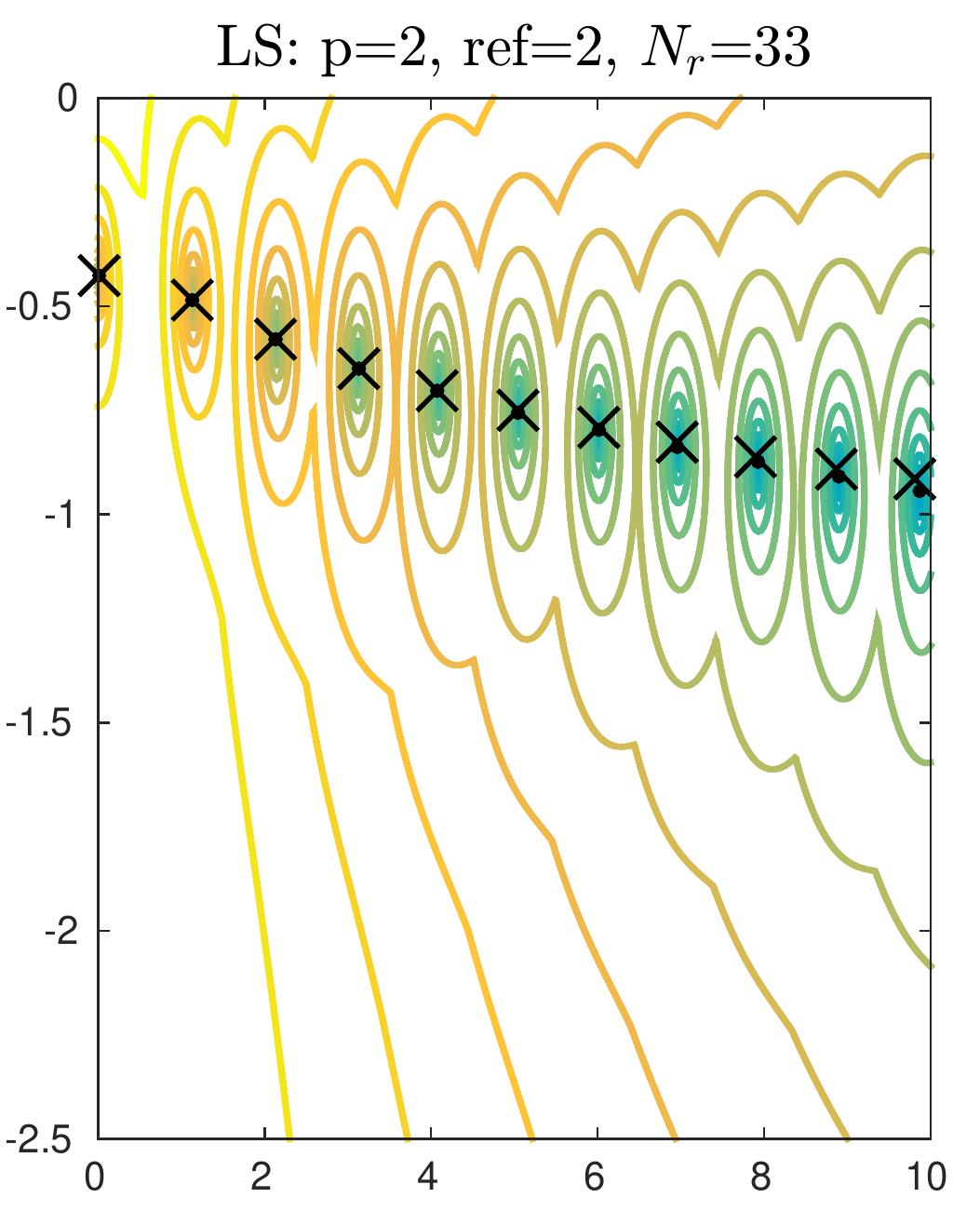} };
		
		
		
	\end{tikzpicture}
	\vspace*{-15mm}
	\caption{\emph{Pseudospectrum for the Bump problem:  Panels on the top are computed using the DtN based formulation, in the middle with the PML based formulation, and in the bottom using the Lippmann-Schwinger formulation. For reference, we mark with crosses $(\times)$ the reference eigenvalues listed in Table \ref{tab:bump_reference}.}}
	\label{fig:psudo_bump}
\end{figure}

\textit{The PML based formulation:}
Fig. \ref{fig:i_check_PML_air_filled} shows eigenvalues $k_j^\fem$ computed using equation \eqref{eq:mat_pml}
with $\sigma_0=5$. Spurious solutions are also present in the PML based formulation. However, in the PML formulation we also have eigenvalues close to the critical line \eqref{eq:CriticalLine} that are not approximations to the resonances $k_j$. Furthermore, we observe that $k_0$ cannot be approximated with any PML discretizaion because $k_0$ lies under the critical line with $\re k_0=0$. Similarly we see that $k_1$ is located very close to the critical line and the value $\epsilon_1$ does not get below $10^{-3}$ for any discretization. The approximation $k_j^\fem$ corresponding to $k_2$ has $\epsilon_j\approx 0$ for finer discretizations.

As discussed above, the number of eigenvalues computed with the Lippmann-Schwinger formulation [Fig. \ref{fig:psudo_air_filled} $a1)$ and $b1)$] equal the number of exact resonances. In the same region, [Fig. \ref{fig:i_check_PML_air_filled} $a)$ and $b)$] the PML formulation shows spurious eigenvalues.

\begin{figure} 
	\begin{tikzpicture}[thick,scale=0.9, every node/.style={scale=0.9}]
		\draw( 0.00, 0.0)  node {
			\includegraphics[scale=0.45]{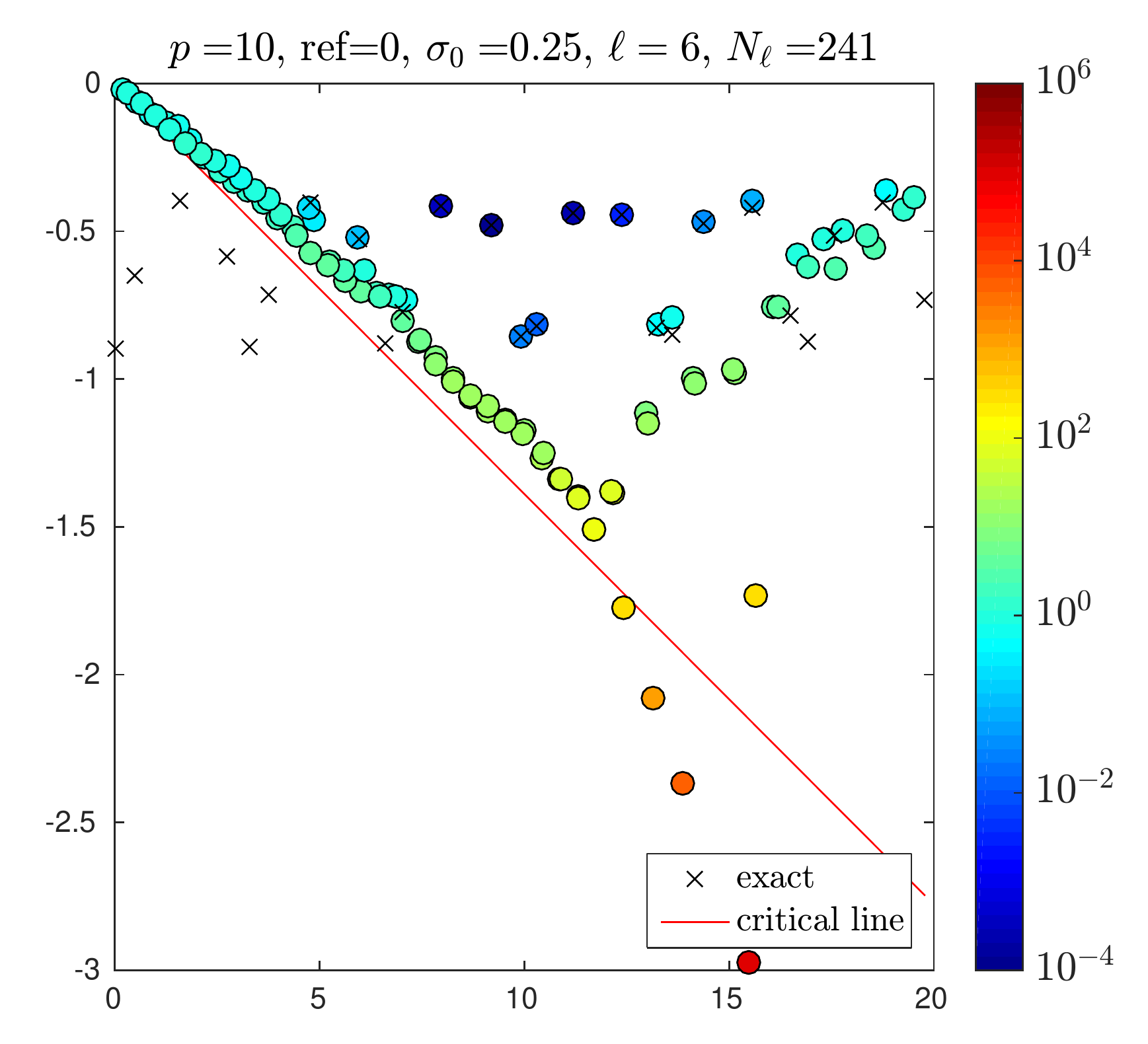} };
		\draw( 8.00, 0.0)  node {
			\includegraphics[scale=0.45]{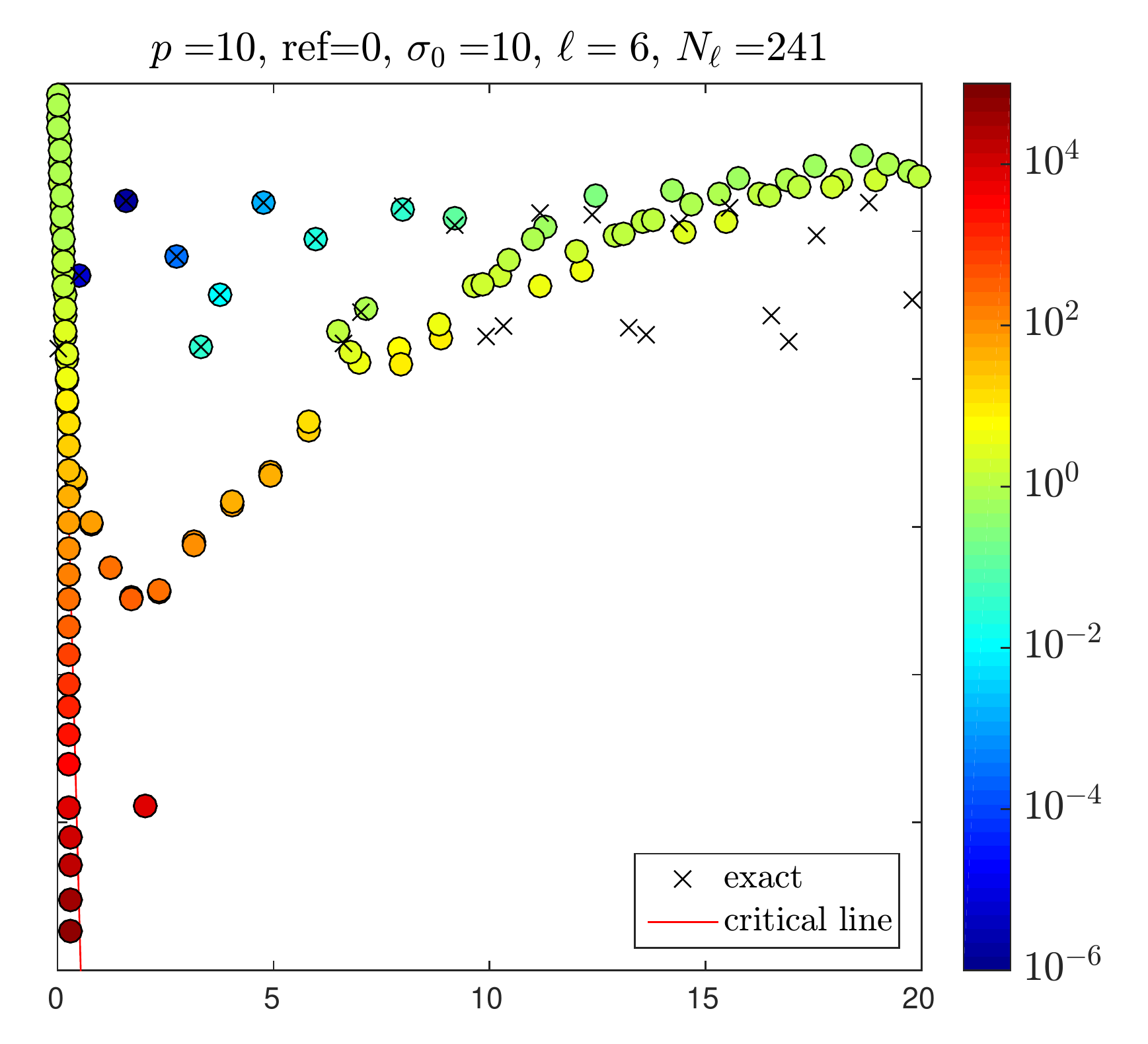} };
		
		\node[text width=3cm] at (-1.70, 3.30) {$a)$};
		\node[text width=3cm] at ( 6.00, 3.30) {$b)$};
		
	\end{tikzpicture}
	\vspace*{-5mm}
	\caption{\emph{Comparison of computed eigenvalues $k_j^\fem$ of \eqref{eq:mat_pml} corresponding to the air-filled-cavity problem described in Sect. \ref{sec:air_filled} with different $\sigma_0$. In colours we give $\epsilon_j$ for each $k_j^\fem$.} }
	\label{fig:PML_air_filled_sigma}
\end{figure}

In the DtN formulation, the finite element space used in Fig. \ref{fig:i_check_DtN_air_filled} d) was sufficient to clear the chosen region from spurious eigenvalues. The same discretization of the resonator region was used in the PML formulation but Fig. \ref{fig:i_check_PML_air_filled} d) shows several spurious eigenvalues in the feasible region. However, 
an increase of the polynomial degree to $p>22$ results in no spurious eigenvalues in the shown region bounded by the critical line \eqref{eq:CriticalLine}. By comparing Fig. \ref{fig:i_check_DtN_air_filled} and Fig. \ref{fig:i_check_PML_air_filled} for equivalent discretizations it is clear that the pairs from the DtN-map have smaller $\epsilon_j$ than those of the finite PML for the chosen $\pml$.
\newline
Fig. \ref{fig:PML_air_filled_sigma} shows eigenvalue approximations corresponding to $\sigma_0=1/4$ and to $\sigma_0=10$. The same FE discretization was used for both $\sigma_0$ and these computations verify that the smallest $\epsilon$ such that $k_j^\nu\in\sigma_{\epsilon}(T^\fem)$ depends critically on $\sigma_0$. 

\subsubsection{Results for the Bump problem}\label{sec:res_bump} 
We use the refractive index given in Sect. \eqref{sec:bump}, and compute eigenpairs and pseudospectrum with the DtN, PML and Lippmann-Schwinger (LS) based formulations, correspondingly from \eqref{eq:mat_dtn}, \eqref{eq:mat_pml} and \eqref{eq:lipp_collocation}. 
In the Bump problem there are no exact eigenpairs but we list in Table \ref{tab:bump_reference} reference values computed with DtN-FEM, with $p=20$, $h=0.125$, and $N_\dtn=481$.

The reference solutions, are used in the following experiment: Computations are performed for all formulations, polynomial order fixed to $p=2$ and show results for three consecutively refined meshes. 
Figure \ref{fig:psudo_bump}, depicts reference eigenvalues with ($\times$), computed eigenvalues with (\tikz\draw[black,fill=black] (0,0) circle (.5ex);), and pseudospectrum in colored contours.
The results follow the lines discussed in sections \ref{sec:res_single_slab} and \ref{sec:res_air_filled}. 
In the regions shown, the number of eigenvalues for the DtN and PML formulations is greater than the number of reference solutions, while the LS has equal number of eigenvalues compared to the number of reference solutions. An exception is the extremely coarse discretization with $N_r=9$, where in the selected region the LS results in less number of numerical eigenvalues than reference eigenvalues. Moreover, we see that the resolvent norm of the LS computations is large only close to the corresponding spectrum, while the resolvent norm for the DtN and PML formulations are large even away from the spectrum.

The pseudospectrum of the discretized Lippmann-Schwinger operator in  Fig. \ref{fig:psudo_air_filled} and Fig. \ref{fig:psudo_bump} shows that it is robust also for a course discretization. The test \eqref{eq:eps_check} is therefore successful in identifying spurious eigenvalues from the DtN-FEM and PML-FEM.

\section{Conclusions and Outlook}\label{sec:Con}

In this paper, we have discussed the approximation of resonances of the Helmholtz problem in open domains. Particularly, we give a new characterization of the spurious eigenvalues arising from the truncated PML problem by introducing a DtN for the finite PML in one dimension. This formulation is then used to derive a new error estimate and reference solutions. Furthermore, we propose a method to detect spurious solutions from DtN and PML computations, when the finite element method is used to approximate resonances.
In our numerical experiments with the DtN map and with the truncated PML, spurious solutions were in all cases present for moderately large finite element spaces but never for finite element spaces with very good approximation properties. In our computations optimal convergence rates are reached for solutions of the discrete problems \eqref{eq:mat_dtn} and \eqref{eq:mat_pml}. Our method to detect spurious solutions could numerically distinguish between approximating eigenpairs and spurious solutions even when true resonance and spurious solutions are mixed. All computations were in one dimension but the PML formulation is similar 
in higher dimensions. Future work includes the efficient evaluation of the Lippman-Schwinger equation in higher dimensions and efficient calculations of eigenvalues with the DtN formulation of the resonance problem.


\section*{Acknowledgements} 
This work is founded by the Swedish Research Council under Grant No.\ $621$-$2012$-$3863$.

\bibliography{extracted}
\bibliographystyle{siam}

\end{document}